\documentclass[a4paper]{amsart}

\usepackage[mathscr]{eucal}
\usepackage{amsmath,amsfonts,amssymb}
\usepackage{mathrsfs}
\usepackage{bm}
\usepackage{multirow}
\usepackage[dvips]{graphicx}
\graphicspath{{./graphs/}}
\usepackage{fancyhdr}

\newcommand{\mathsym}[1]{{}}

\theoremstyle{plain}

\newtheorem{conjecture0}{Conjecture}

\theoremstyle{definition}

\begin{document}
\quad\vspace{-0.1in}

\title[On the zeros of certain modular functions for the normalizers of congruence subgroups]{On the zeros of certain modular functions for \\ the normalizers of congruence subgroups of low levels \\ II}
\author{Junichi Shigezumi}

\maketitle \vspace{-0.1in}
\begin{center}
Graduate School of Mathematics Kyushu University\\
Hakozaki 6-10-1 Higashi-ku, Fukuoka, 812-8581 Japan\\
{\it E-mail address} : j.shigezumi@math.kyushu-u.ac.jp \vspace{-0.05in}
\end{center} \quad

\begin{quote}
{\small\bfseries Abstract.}
We research the location of the zeros of the Eisenstein series and the modular functions from the Hecke type Faber polynomials associated with the normalizers of congruence subgroups which are of genus zero and of level at most twelve.

In Part II, we will observe the location of the zeros of the above functions by numerical calculation.\\  \vspace{-0.15in}

\vspace{-0.15in}

 \vspace{0.15in}
\end{quote}

\section*{Introduction}

The motive of this research is to decide the location of the zeros of modular functions. The Eisenstein series and the Hecke type Faber polynomials are the most interesting and important modular forms.

F. K. C. Rankin and H. P. F. Swinnerton-Dyer considered the problem of locating the zeros of the Eisenstein series $E_k(z)$ in the standard fundamental domain $\mathbb{F}$ (See \cite{RSD}). They proved that all of the zeros of $E_k(z)$ in $\mathbb{F}$ lie on the unit circle. They also stated towards the end of their study that ``This method can equally well be applied to Eisenstein series associated with subgroups of the modular group.'' However, it seems unclear how widely this claim holds. 

Subsequently, T. Miezaki, H. Nozaki, and the present author considered the same problem for the Fricke group $\Gamma_0^{*}(p)$ (see \cite{Kr}, \cite{Q}), and proved that all of the zeros of the Eisenstein series $E_{k, p}^{*}(z)$ in a certain fundamental domain lie on a circle whose radius is equal to $1 / \sqrt{p}$,  $p = 2, 3$ (see \cite{MNS}). Furthermore, we also proved that almost all the zeros of the Eisenstein series in a certain fundamental domain lie on circles whose radius are equal to $1 / \sqrt{p}$ or $1 / (2 \sqrt{p})$,  $p = 5, 7$ (see \cite{SJ2}).

Let $\Gamma$ be a discrete subgroup of $\text{\upshape SL}_2(\mathbb{R})$, and let $h$ be the width of $\Gamma$, then we define
\begin{equation}
\mathbb{F}_{0,\Gamma} := \left\{z \in \mathbb{H} \: ; \: - h / 2 < Re(z) < h / 2 \: , \: |c z + d| > 1 \: \text{for} \: \forall \gamma = \left(\begin{smallmatrix} a & b \\ c & d \end{smallmatrix}\right) \in \Gamma \: \text{s.t.} \: c \ne 0 \right\} \label{def-fd*}.
\end{equation}
We have a fundamental domain $\mathbb{F}_{\Gamma}$ such that $\mathbb{F}_{0,\Gamma} \subset \mathbb{F}_{\Gamma} \subset \overline{\mathbb{F}_{0,\Gamma}}$. Let $\mathbb{F}_{\Gamma}$ be such a fundamental domain.

For the modular group $\text{SL}_2(\mathbb{Z})$ and the Fricke groups $\Gamma_0^{*}(p)$ ($p = 2, 3$), all the zeros of the Eisenstein series for the cusp $\infty$ lie on the arcs on the boundary of their certain fundamental domains.

H. Hahn considered that the location of the zeros of the Eisenstein series for the cusp $\infty$ for every genus zero Fucksian group $\Gamma$ of the first kind with $\infty$ as a cusp which satisfies that its hauptmodul $J_{\Gamma}$ takes real value on $\partial \mathbb{F}_{\Gamma}$, and proved that almost all the zeros of the Eisenstein series for the cusp $\infty$ for $\Gamma$ lie on $\partial \mathbb{F}_{\Gamma}$ under some more assumption (see \cite{H}).

Also, T. Asai, M. Kaneko, and H. Ninomiya considered the problem of locating the zeros of modular functions $F_m(z)$ for $\text{SL}_2(\mathbb{Z})$ which correspond to the Hecke type Faber polynomial $P_m$, that is, $F_m(z) = P_m (J(z))$ (See \cite{AKN}). They proved that all of the zeros of $F_m(z)$ in $\mathbb{F}$ lie on the unit circle for each $m \geqslant 1$. After that, E. Bannai, K. Kojima, and T. Miezaki considered the same problem for the normalizers of congruence subgroups which correspond the conjugacy classes of the Monster group (See \cite{BKM}). They observed the location of the zeros by numerical calculation, then almost all of the zeros of the modular functions from Hecke type Faber polynomial lie on the lower arcs when the group satisfy the same assumption of the theorem of H. Hahn. In particular, T. Miezaki proved that all of the zeros of the modular functions from the Hecke type Faber polynomials for the Fricke group $\Gamma_0^{*}(2)$ lie on the lower arcs of its fundamental domain in their paper.

Now, we have the following conjectures:
\begin{conjecture0}\label{conj0}
Let $\Gamma$ be a genus zero Fucksian group of the first kind with $\infty$ as a cusp. If the hauptmodul $J_{\Gamma}$ takes real value on $\partial \mathbb{F}_{\Gamma}$, all of the zeros of the Eisenstein series for the cusp $\infty$ for $\Gamma$ in $\mathbb{F}_{\Gamma}$ lie on the arcs
\begin{equation*}
\partial \mathbb{F}_{\Gamma} \setminus \{z \in \mathbb{H} \: ; \: Re(z) = \pm h / 2\}.
\end{equation*}
\end{conjecture0}

\begin{conjecture0}\label{conj01}
Let $\Gamma$ be a genus zero Fucksian group of the first kind with $\infty$ as a cusp. If the hauptmodul $J_{\Gamma}$ takes real value on $\partial \mathbb{F}_{\Gamma}$, all but at most $c_h(\Gamma)$ of the zeros of modular function from the Hecke type Faber polynomial of degree $m$ for $\Gamma$ in $\mathbb{F}_{\Gamma}$ lie on the arcs
\begin{equation*}
\partial \mathbb{F}_{\Gamma} \setminus \{z \in \mathbb{H} \: ; \: Re(z) = \pm h / 2\}
\end{equation*}
for all but finite number of $m$ and for the constant number $c_h(\Gamma)$ which does not depend on $m$.
\end{conjecture0}

In this paper, we will observe the location of the zeros of the Eisenstein series and the modular functions from Hecke type Faber polynomials for the normalizers of congruence subgroups, as a first step of a challenge for the above conjectures.\\

The normalizers of congruence subgroups of level at most $12$ which satisfies the assumption of above conjectures are 
\begin{align*}
&\text{SL}_2(\mathbb{Z}), \; \Gamma_0^{*}(2), \; \Gamma_0(2), \; \Gamma_0^{*}(3), \; \Gamma_0(3), \; \Gamma_0^{*}(4), \; \Gamma_0(4), \; \Gamma_0^{*}(5), \; \Gamma_0(6)+, \; \Gamma_0^{*}(6), \; \Gamma_0(6)+3, \; \Gamma_0(6),\\
&\Gamma_0^{*}(7), \; \Gamma_0^{*}(8), \; \Gamma_0(8), \; \Gamma_0^{*}(9), \; \Gamma_0(10)+, \; \Gamma_0^{*}(10), \; \Gamma_0(10)+5, \; \Gamma_0(12)+, \; \Gamma_0^{*}(12),\\
&\Gamma_0(12)+4, \; \text{and} \; \Gamma_0(12).
\end{align*}

For the Conjecture \ref{conj0}, $\text{SL}_2(\mathbb{Z})$, $\Gamma_0^{*}(2)$, and $\Gamma_0^{*}(3)$ verify  Conjecture \ref{conj0}. For the other cases, we can prove by numerical calculation for the Eisenstein series of weight $k \leqslant 500$.

For the Conjecture \ref{conj01}, $\text{SL}_2(\mathbb{Z})$ and $\Gamma_0^{*}(2)$ verify  Conjecture \ref{conj01} for every degree $m$, where we have $c_h(\Gamma) = 0$ for each case. Furthermore, for $\Gamma_0(2)$, $\Gamma_0^{*}(3)$, $\Gamma_0(3)$, $\Gamma_0^{*}(4)$, $\Gamma_0(4)$, $\Gamma_0(6)+$, $\Gamma_0(6)+3$, $\Gamma_0(6)$, $\Gamma_0(8)$, $\Gamma_0^{*}(9)$, $\Gamma_0(10)+$, $\Gamma_0(10)+5$, $\Gamma_0(12)+$, $\Gamma_0(12)+4$, and $\Gamma_0(12)$, we can prove all of the zeros of the modular function from the Hecke type Faber polynomial of every degee $m \leqslant 200$ in each fundamental domain lie on the lower arcs  by numerical calculation.

On the other hand, for $\Gamma_0^{*}(5)$ and $\Gamma_0^{*}(7)$, we can prove by numerical calculation for the modular function from the Hecke type Faber polynomial of every degee $m = 1$ and $3 \leqslant m \leqslant 200$, where we have $c_h(\Gamma) = 0$ for each case. When $m = 2$, there is just one zero which is on the boundary of its fundamental domain but not on the lower arcs for the each group.

For $\Gamma_0^{*}(6)$ and $\Gamma_0^{*}(8)$, we can prove by numerical calculation for the modular function from the Hecke type Faber polynomial of every degee $m \leqslant 200$ which satisfy $m \not\equiv 0 \pmod{2}$ and $m \not\equiv 2 \pmod{4}$, respectively. For the remaining degrees, there is just one zero which is on the boundary of its fundamental domain but not on the lower arcs for the each group, that is, $c_h(\Gamma) = 1$.

Finally, for $\Gamma_0^{*}(10)$ and $\Gamma_0^{*}(12)$, we have just two zeros which are not on the boundary of each fundamental domain for degrees $m =7, 9, 11$ and $m = 3, 6, 12, 13, 15$, respectively. Furthermore, there is just one zero which is on the boundary of its fundamental domain but not on the lower arcs for the case $m \equiv 0 \pmod{2}$ and $m \equiv 2, 4 \pmod{6}$, respectively. For the other cases, we can prove that all of the zeros are on the lower arcs of each fundamental domain by numerical calculation.

\begin{table}[htbp]
{\small \begin{center}
\begin{tabular}{ccc}
\hline
\hspace{-0.5in}$\Gamma$ & \hspace{-0.3in}Eisenstein series $(k \leqslant 500)$ & Hecke type Faber polynomial $(m \leqslant 200)$\\
\hline
\begin{tabular}{l}
$\text{SL}_2(\mathbb{Z})$, $\Gamma_0^{*}(2)$, $\Gamma_0(2)$, $\Gamma_0^{*}(3)$, $\Gamma_0(3)$,\\ $\Gamma_0^{*}(4)$, $\Gamma_0(4)$, $\Gamma_0(6)+$, $\Gamma_0(6)+3$,\\
$\Gamma_0(6)$, $\Gamma_0(8)$, $\Gamma_0^{*}(9)$, $\Gamma_0(10)+$, $\Gamma_0(10)+5$,\\
$\Gamma_0(12)+$, $\Gamma_0(12)+4$, $\Gamma_0(12)$.
\end{tabular} & \multirow{6}{*}{\hspace{-0.15in}$\bigcirc$} & $\bigcirc$\\
$\Gamma_0^{*}(5)$, $\Gamma_0^{*}(7)$ && $m = 2$, \quad $\langle 1 \rangle$\\
$\Gamma_0^{*}(6)$ && $m$ : even, \quad $\langle 1 \rangle$\\
$\Gamma_0^{*}(8)$ && $m \equiv 0 \pmod{4}$, \quad $\langle 1 \rangle$\\
$\Gamma_0^{*}(10)$ && $m = 7, 9, 11, \; [2], \quad m$ : even, $\langle 1 \rangle$\\
$\Gamma_0^{*}(12)$ && $m = 3, 6, 12, 13, 15, \; [2] \quad m \equiv 2, 4 \pmod{6}, \; \langle 1 \rangle$\\
\hline
\end{tabular}
\end{center}
\begin{flushleft}
\hspace{0.51in}`$\bigcirc$': all of the zeros lie on lower arcs.\\
\hspace{0.5in}$\langle \ \cdot \ \rangle$: the number of zeros which are on $\partial \mathcal{F}$ but not on lower arcs.\\
\hspace{0.51in}$[ \ \cdot \ ]$\hspace{0.02in}: the number of zeros which are not on $\partial \mathcal{F}$.
\end{flushleft}}
\caption{Result by numerical calculation}
\end{table}

\begin{figure}[hbtp]
\begin{center}
{{\small Lower arcs of $\partial \mathbb{F}_{6+2}$}\includegraphics[width=2in]{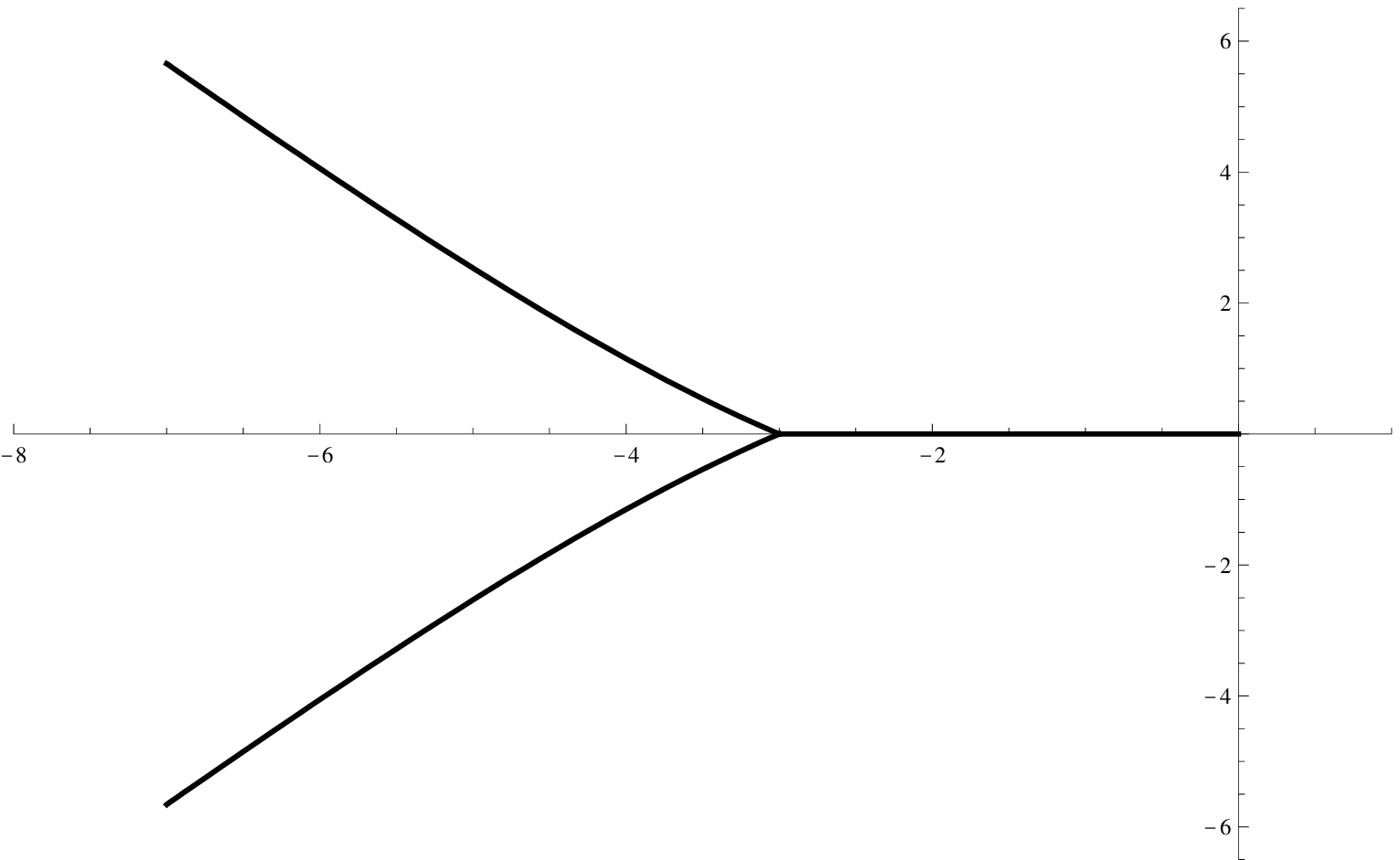}}
\end{center}
\caption{Image by $J_{6+2}$ ($\Gamma_0(6)+2$)}\label{Im-J6D-intro0}
\end{figure}

If the hauptmodul $J_{\Gamma}$ does not take real value on $\partial \mathbb{F}_{\Gamma}$ (cf. Figure \ref{Im-J6D-intro0}), it seems to be not similar. Such cases are followings;
\begin{equation*}
\Gamma_0(5), \; \Gamma_0(6)+2, \; \Gamma_0(7), \; \Gamma_0(9), \; \Gamma_0(10)+2, \; \Gamma_0(10), \; \Gamma_0^{*}(11), \; \text{and} \; \Gamma_0(12)+3.
\end{equation*}

\begin{figure}[hbtp]
\begin{center}
{{$\begin{matrix}\text{\small The zeros of $E_{k, 6+2}^{\infty}$} \\ \text{\small for $4 \leqslant k \leqslant 40$}\end{matrix}$}\includegraphics[width=4.2in]{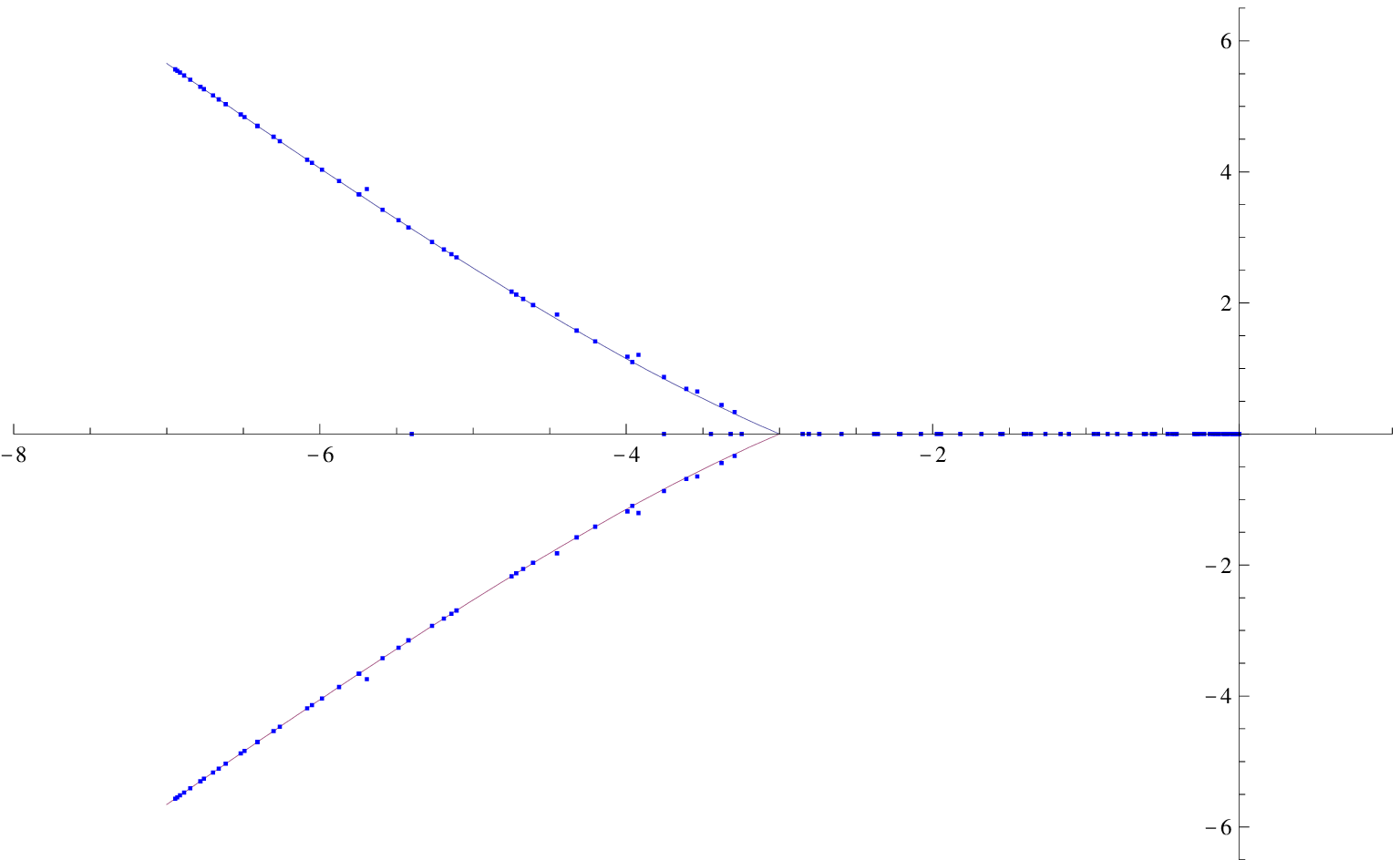}}\\
{{\small The zeros of $E_{1000, 6+2}^{\infty}$}\includegraphics[width=4.2in]{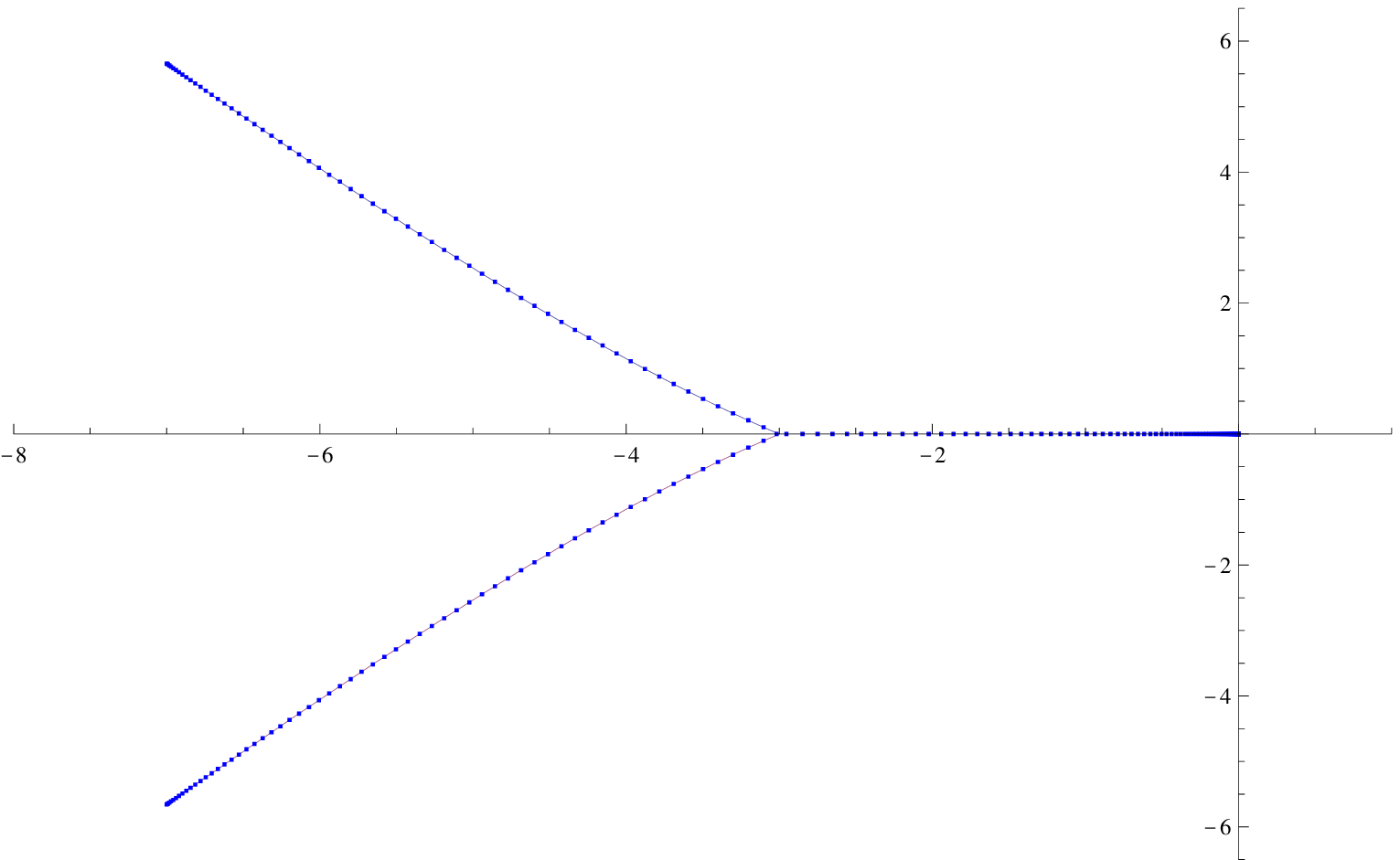}}
\end{center}
\caption{Image by $J_{6+2}$ ($\Gamma_0(6)+2$)}\label{Im-J6D-intro1}
\end{figure}

For $\Gamma_0(5)$, $\Gamma_0(6)+2$, $\Gamma_0(7)$, $\Gamma_0(10)+2$, $\Gamma_0(10)$, and $\Gamma_0^{*}(11)$, we can observe that the zeros of the Eisenstein series for cusp $\infty$ do not lie on the lower arcs of their fundamental domains by numerical calculation. However, when the weight of Eisenstein series increases, then the location of the zeros seems to approach to lower arcs. (See Figure \ref{Im-J6D-intro1})

\begin{figure}[hbtp]
\begin{center}
{{$\begin{matrix}\text{\small The zeros of $F_{m, 6+2}$} \\ \text{\small for $1 \leqslant m \leqslant 40$}\end{matrix}$}\includegraphics[width=4.2in]{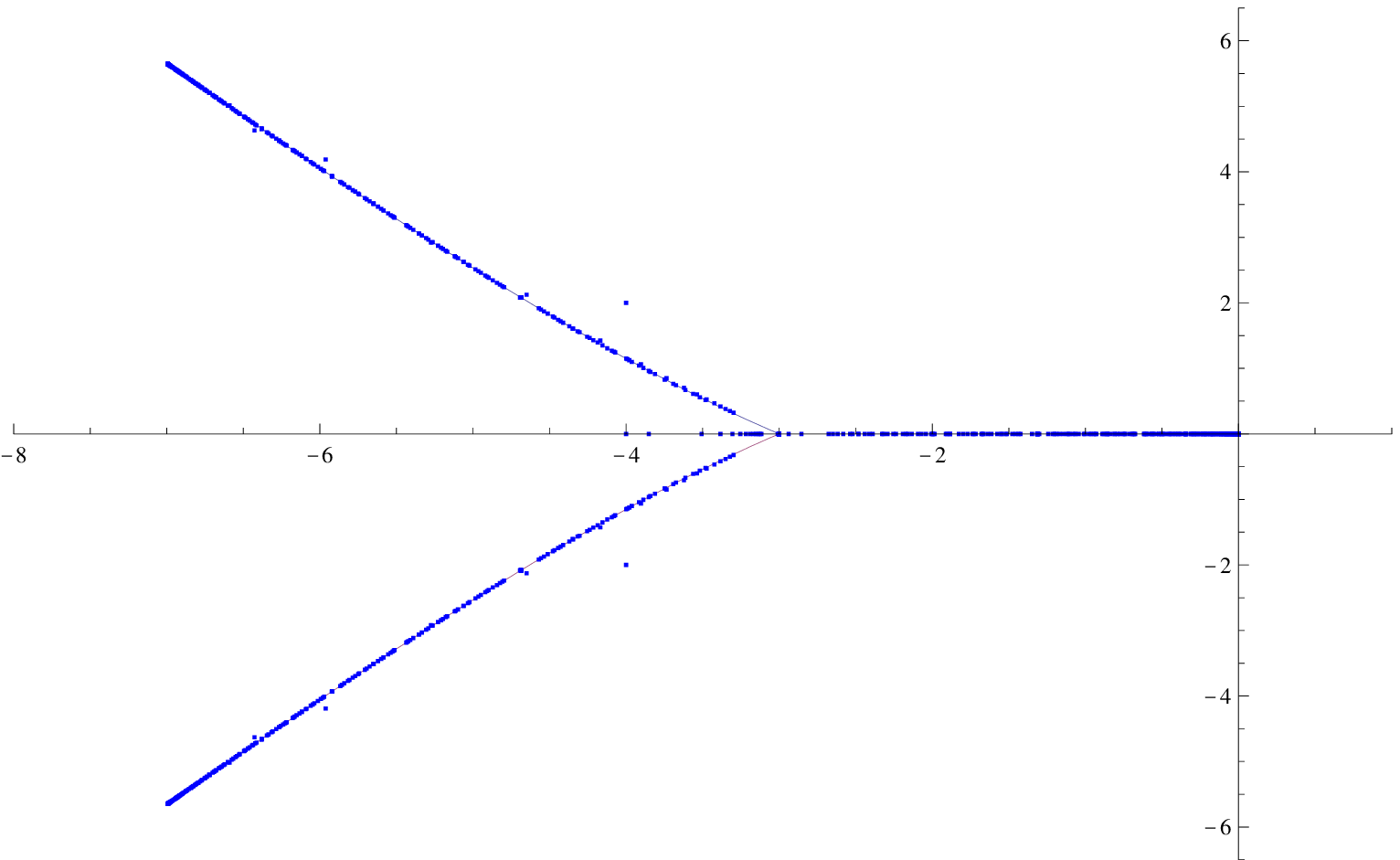}}\\
{{\small The zeros of $F_{200, 6+2}$}\includegraphics[width=4.2in]{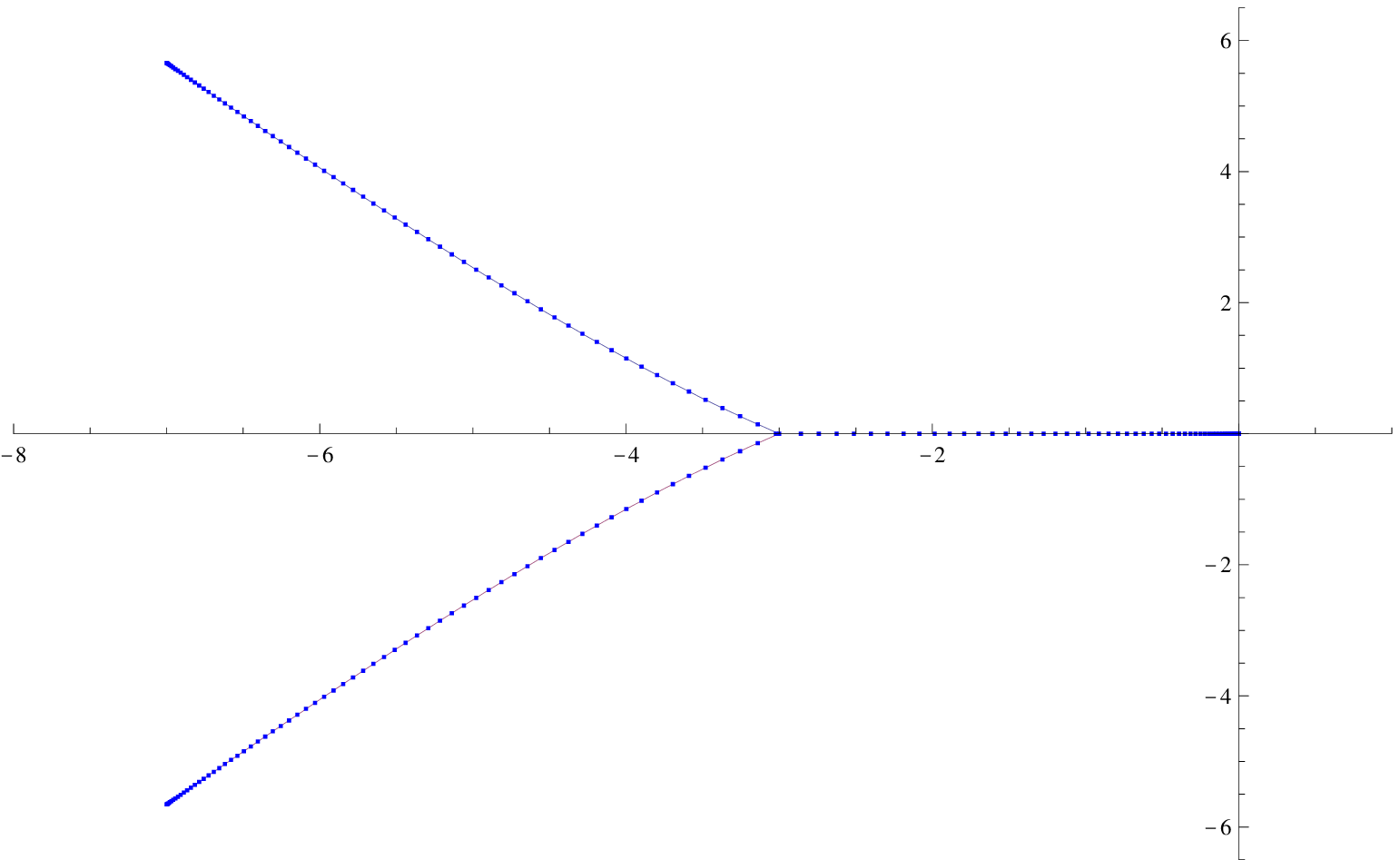}}
\end{center}
\caption{Image by $J_{6+2}$ ($\Gamma_0(6)+2$)}\label{Im-J6D-intro2}
\end{figure}

Also, for the zeros of the modular functions from the Hecke type Faber polynomials, we can observe that there are some zeros which do not lie on the lower arcs of their fundamental domains by numerical calculation. Furthermore, when the degree $m$ increases, then the location of the zeros seems to approach to lower arcs. (See Figure \ref{Im-J6D-intro2})

\begin{figure}[hbtp]
\begin{center}
{{\small $\Gamma_0(9)$}\includegraphics[width=2in]{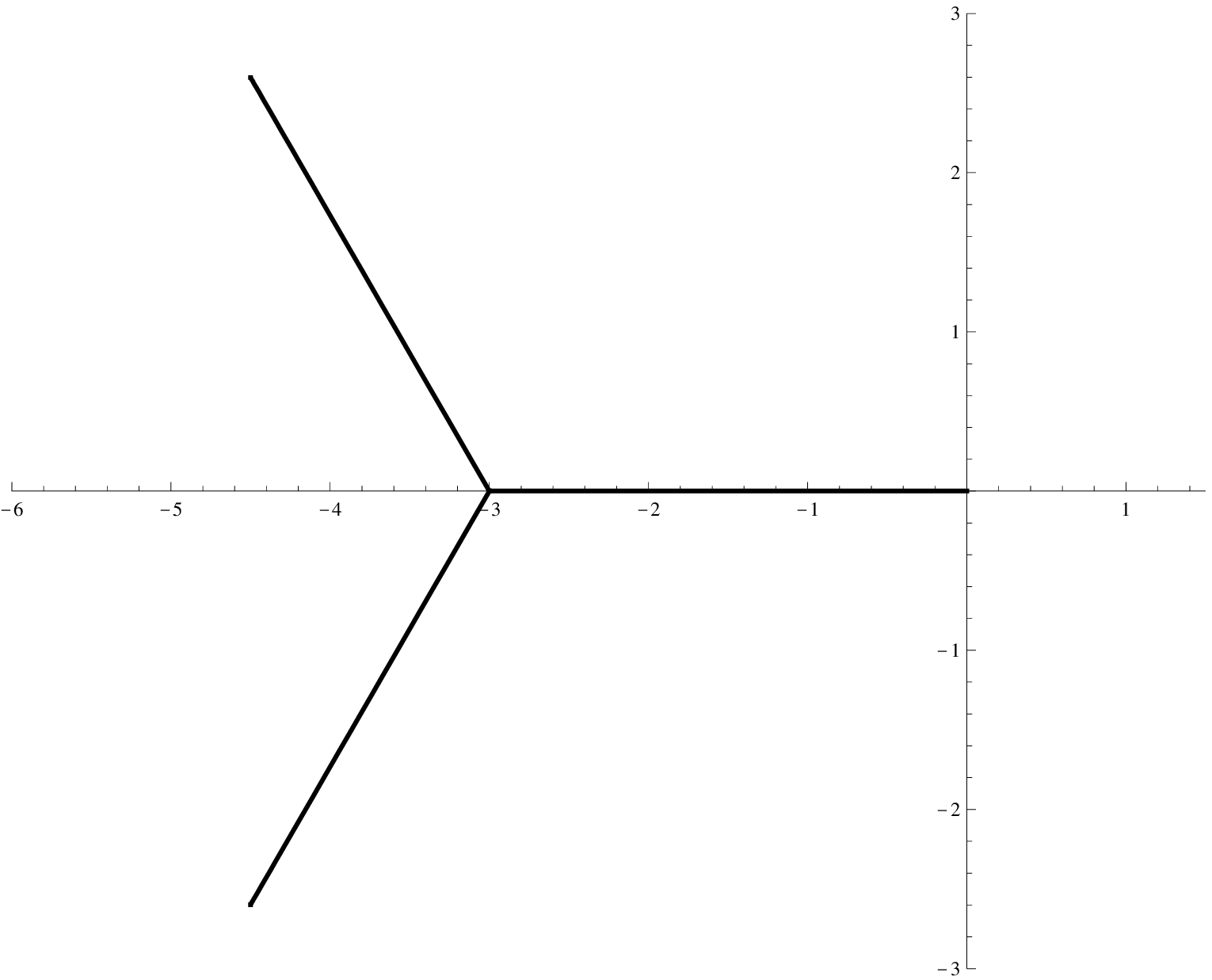}} \;
{{\small $\Gamma_0(12)+3$}\includegraphics[width=2.6in]{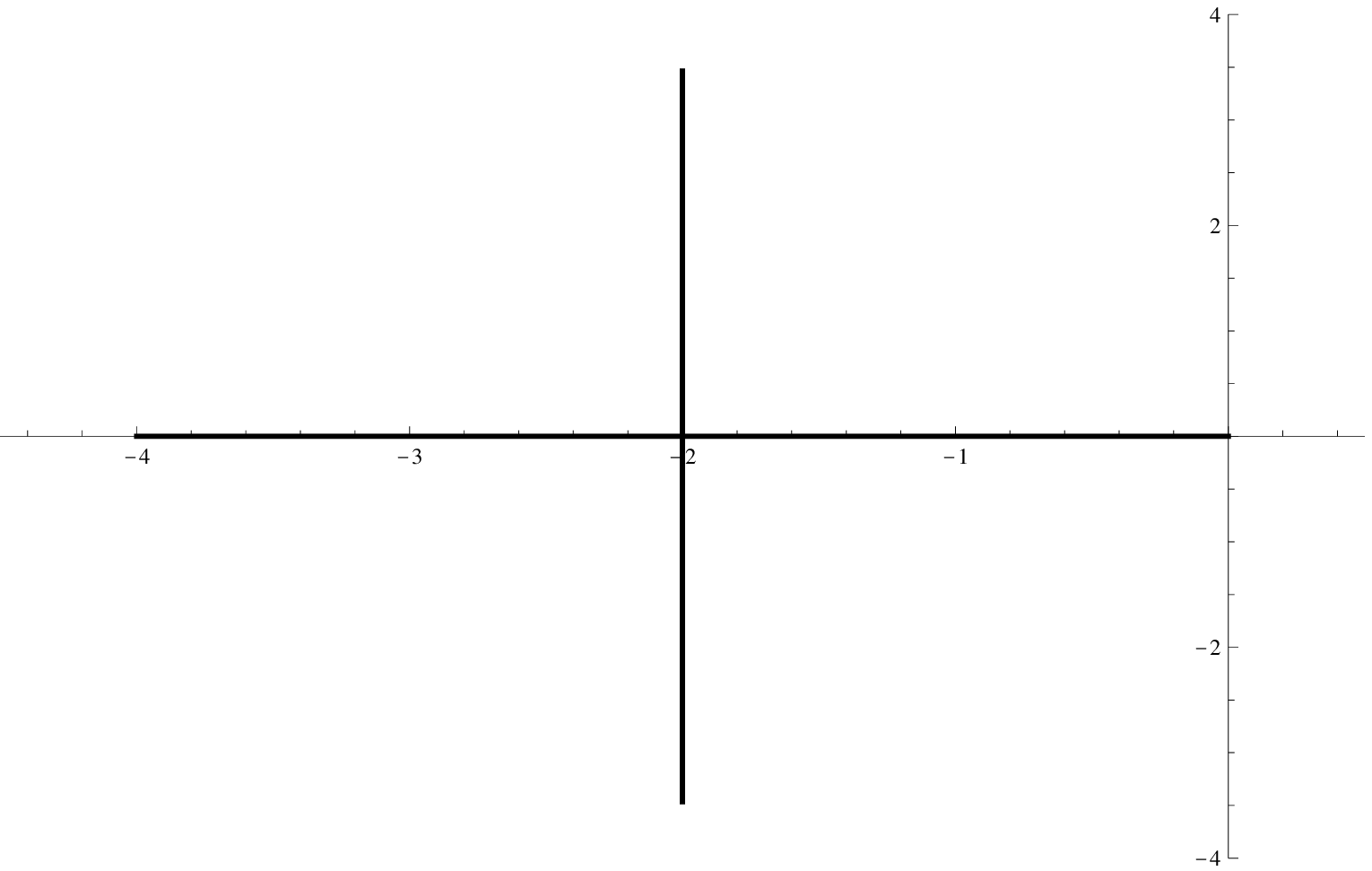}}
\end{center}
\caption{Image of the lower arcs of the fundamental domains by hauptmoduls}\label{Im-SP-intro}
\end{figure}

On the other hand, $\Gamma_0(9)$ and $\Gamma_0(12)+3$ seem to show the special cases. We can prove that all of the zeros of the Eisenstein series of weight $k \leqslant 500$ lie on the lower arcs of their fundamental domains by numerical calculation. Also, we can prove that all of the zeros of the modular function from the Hecke type Faber polynomial of degee $m \leqslant 200$ lie on the lower arcs by numerical calculation. On the other hand, they do not satisfy the assumption of Conjecture \ref{conj0} and \ref{conj01}. However, the image of lower arcs by its hauptmodul draw a interesting figure. (Figure \ref{Im-SP-intro})

\newpage
\quad

We refer to \cite{MNS}, \cite{SJ1}, and \cite{SJ2} for some groups. However, note that definitions in this paper are sometimes different from that in it.\\

In `Part I', we will consider the general theory of modular functions for the normalizers of the congruence subgroups $\Gamma_0(N)$ of level $N \leqslant 12$. And in `Part II', we will observe the location of the zeros of the Eisenstein series and the the modular functions from Hecke type Faber polynomials for the normalizers in Part I by numerical calculation.

\clearpage

\markboth{\sc On the zeros of Eisenstein Series for the normalizers of congruence subgroups}{\sc Level \thesection}

\section{Level $1$}

\subsection{$\text{SL}_2(\mathbb{Z})$}

We have $\text{SL}_2(\mathbb{Z}) = \langle \left( \begin{smallmatrix} 1 & 1 \\ 0 & 1 \end{smallmatrix} \right), \: \left( \begin{smallmatrix} 0 & -1 \\ 1 & 0 \end{smallmatrix} \right) \rangle$.\\

\paragraph{\bf Location of the zeros of the Eisenstein series}
F. K. C. Rankin and H. P. F. Swinnerton-Dyer proved that all of the zeros of $E_k$ lie on the lower arcs of $\partial \mathbb{F}$. (See \cite{RSD})
\begin{figure}[hbtp]
\begin{center}
\includegraphics[width=1.5in]{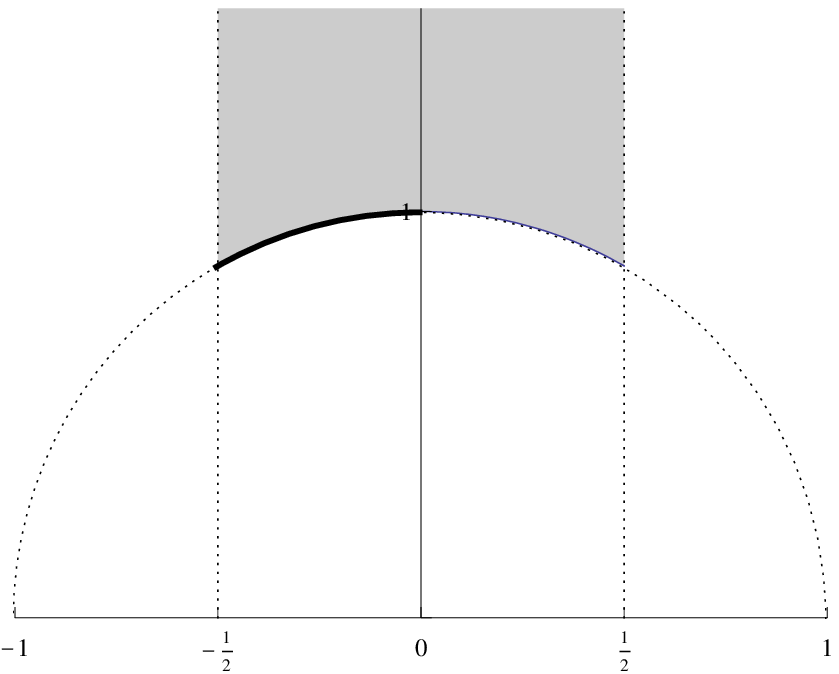}
\end{center}
\caption{Location of the zeros of the Eisenstein series}
\end{figure}

\paragraph{\bf Location of the zeros of Hecke type Faber Polynomial}
T. Asai, M. Kaneko, and H. Ninomiya proved that all of the zeros of $F_m$ lie on the lower arcs of $\partial \mathbb{F}$. (See \cite{AKN})

\clearpage

\section{Level $2$}

We have $\Gamma_0(2)+=\Gamma_0^{*}(2)$ and $\Gamma_0(2)-=\Gamma_0(2)$. We have $W_2 = \left(\begin{smallmatrix}0&-1 / \sqrt{2}\\ \sqrt{2}&0\end{smallmatrix}\right)$, \\

\subsection{$\Gamma_0^{*}(2)$} 

We have $\Gamma_0^{*}(2) = \langle \left( \begin{smallmatrix} 1 & 1 \\ 0 & 1 \end{smallmatrix} \right), \: W_2 \rangle$.\\

\paragraph{\bf Location of the zeros of the Eisenstein series}
T. Miezaki, H. Nozaki, and the present author proved that all of the zeros of $E_{k, 2+}$ lie on the lower arcs of $\partial \mathbb{F}_{2+}$. (See \cite{MNS})
\begin{figure}[hbtp]
\begin{center}
\includegraphics[width=1.5in]{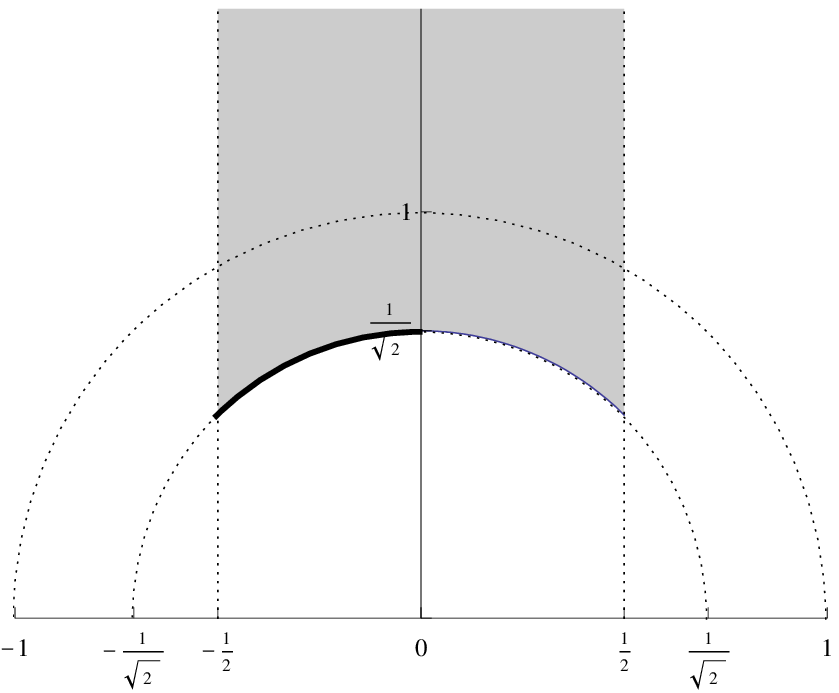}
\end{center}
\caption{Location of the zeros of the Eisenstein series}
\end{figure}

\paragraph{\bf Location of the zeros of Hecke type Faber Polynomial}
T. Miezaki proved that all of the zeros of $F_{m, 2+}$ lie on the lower arcs of $\partial \mathbb{F}_{2+}$. (See \cite{BKM})\\

\subsection{$\Gamma_0(2)$}

We have $\Gamma_0(2) = \langle \left( \begin{smallmatrix} 1 & 1 \\ 0 & 1 \end{smallmatrix} \right), \: \left( \begin{smallmatrix} 1 & 0 \\ 2 & 1 \end{smallmatrix} \right) \rangle$ and $\gamma_0 = W_2$.\\

\paragraph{\bf Location of the zeros of the Eisenstein series}
Since $W_2^{- 1} \Gamma_0(2) W_2 = \Gamma_0(2)$, we have
\begin{equation}
E_{k, 2}^0(W_2 z) = (\sqrt{2} z)^k E_{k, 2}^{\infty}(z).
\end{equation}
Furthermore, we have
\begin{equation*}
E_{k, 2}^0 (- 1/2 + i / (2 \tan\theta/2)) = ((e^{i \theta} - 1) / \sqrt{2})^k E_{k, 2}^{\infty}((e^{i \theta} - 1) / 2).
\end{equation*}
Then, if we have the zeros of $E_{k, 2}^{\infty}$ in the lower arcs of $\partial \mathbb{F}_2$, then we have the zeros of $E_{k, 2}^0$ in $\{ z \: ; \: Re(z) = - 1/2 \}$. (See the below figure)

For $k \leqslant 1000$, we can prove that all of the zeros of $E_{k, 2}^{\infty}$ lie on the lower arcs of $\partial \mathbb{F}_2$ by numerical calculation.

\begin{figure}[hbtp]
\begin{center}
{{$E_{k, 2}^{\infty}$}\includegraphics[width=1.5in]{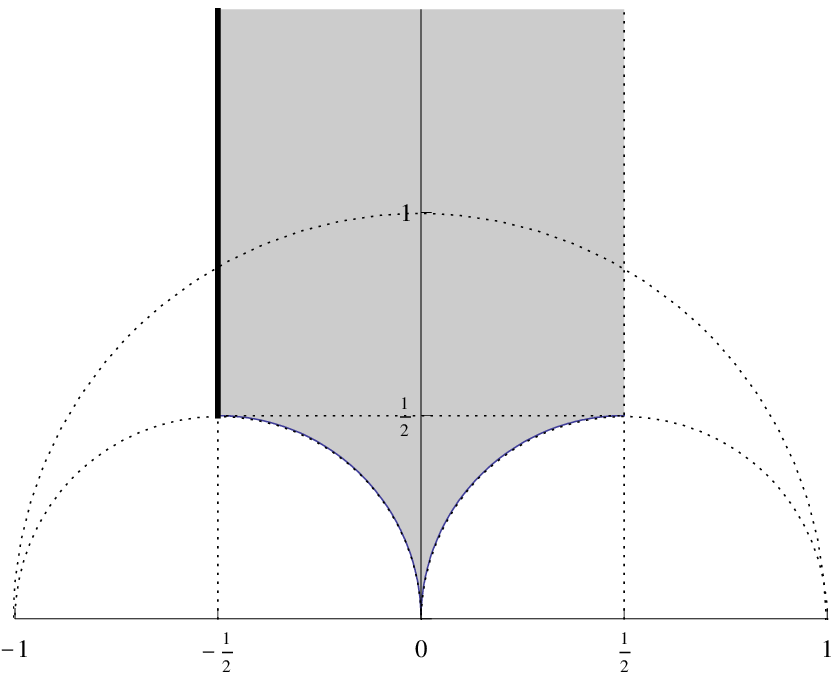}}
\quad \quad
{{$E_{k, 2}^{0}$}\includegraphics[width=1.5in]{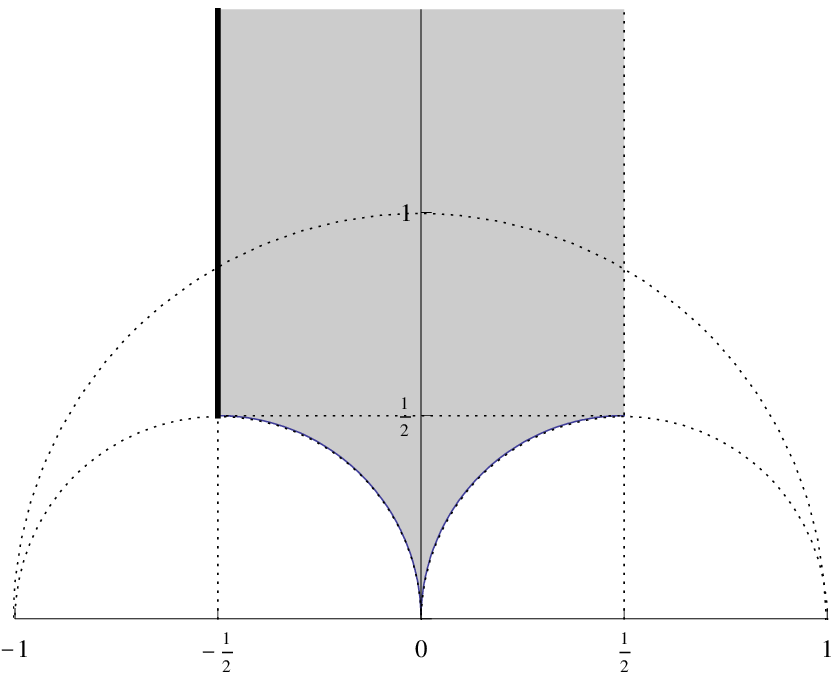}}
\end{center}
\caption{Location of the zeros of the Eisenstein series}
\end{figure}

\paragraph{\bf Location of the zeros of Hecke type Faber Polynomial}
For $m \leqslant 200$, we can prove that all of the zeros of $F_{m, 2}$ lie on the lower arcs of $\partial \mathbb{F}_2$ by numerical calculation.

\clearpage

\section{Level $3$}

We have $\Gamma_0(3)+=\Gamma_0^{*}(3)$ and $\Gamma_0(3)-=\Gamma_0(3)$. We have $W_3 = \left(\begin{smallmatrix}0&-1 / \sqrt{3}\\ \sqrt{3}&0\end{smallmatrix}\right)$.\\

\subsection{$\Gamma_0^{*}(3)$}

We have $\Gamma_0^{*}(3) = \langle \left( \begin{smallmatrix} 1 & 1 \\ 0 & 1 \end{smallmatrix} \right), \: W_3 \rangle$.\\

\paragraph{\bf Location of the zeros of the Eisenstein series}
T. Miezaki, H. Nozaki, and the present author proved that all of the zeros of $E_{k, 3+}$ lie on the lower arcs of $\partial \mathbb{F}_{3+}$.
\begin{figure}[hbtp]
\begin{center}
\includegraphics[width=1.5in]{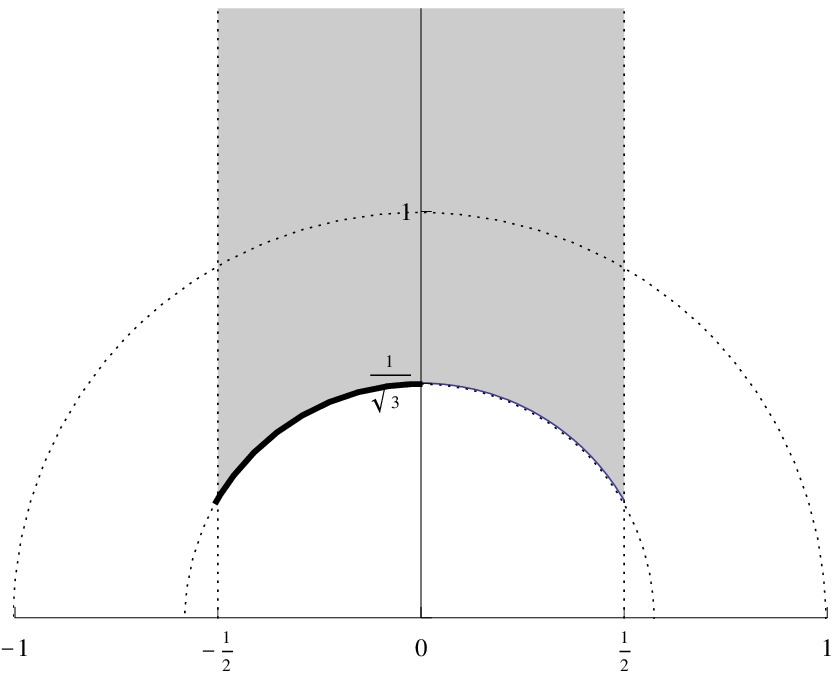}
\end{center}
\caption{Location of the zeros of the Eisenstein series}
\end{figure}

\paragraph{\bf Location of the zeros of Hecke type Faber Polynomial}
For $m \leqslant 200$, we can prove that all of the zeros of $F_{m, 3+}$ lie on the lower arcs of $\partial \mathbb{F}_{3+}$ by numerical calculation.\\

\subsection{$\Gamma_0(3)$}

We have $\Gamma_0(3) = \langle \left( \begin{smallmatrix} 1 & 1 \\ 0 & 1 \end{smallmatrix} \right), \: - \left( \begin{smallmatrix} 1 & 0 \\ 3 & 1 \end{smallmatrix} \right) \rangle$ and $\gamma_0 = W_3$.\\

\paragraph{\bf Location of the zeros of the Eisenstein series}
Since $W_3^{- 1} \Gamma_0(3) W_3 = \Gamma_0(3)$, we have
\begin{equation}
E_{k, 3}^0(W_3 z) = (\sqrt{3} z)^k E_{k, 3}^{\infty}(z).
\end{equation}
Furthermore, we have
\begin{equation*}
E_{k, 3}^0 (- 1/2 + i / (2 \tan\theta/2)) = ((e^{i \theta} - 1) / \sqrt{3})^k E_{k, 3}^{\infty}((e^{i \theta} - 1) / 3).
\end{equation*}

For $k \leqslant 1000$, we can prove that all of the zeros of $E_{k, 3}^{\infty}$ lie on the lower arcs of $\partial \mathbb{F}_3$ by numerical calculation.

\begin{figure}[hbtp]
\begin{center}
{{$E_{k, 3}^{\infty}$}\includegraphics[width=1.5in]{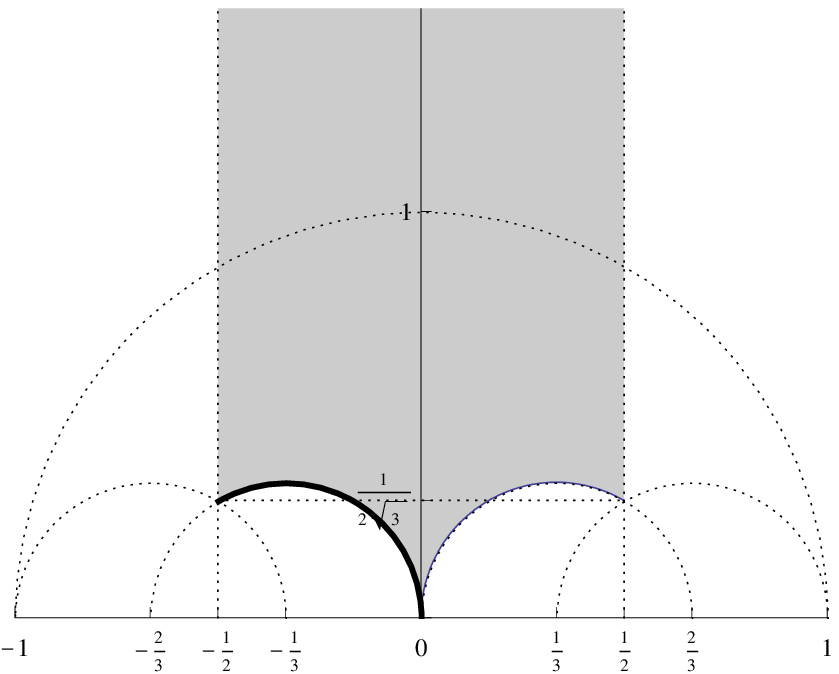}}
\quad \quad
{{$E_{k, 3}^{0}$}\includegraphics[width=1.5in]{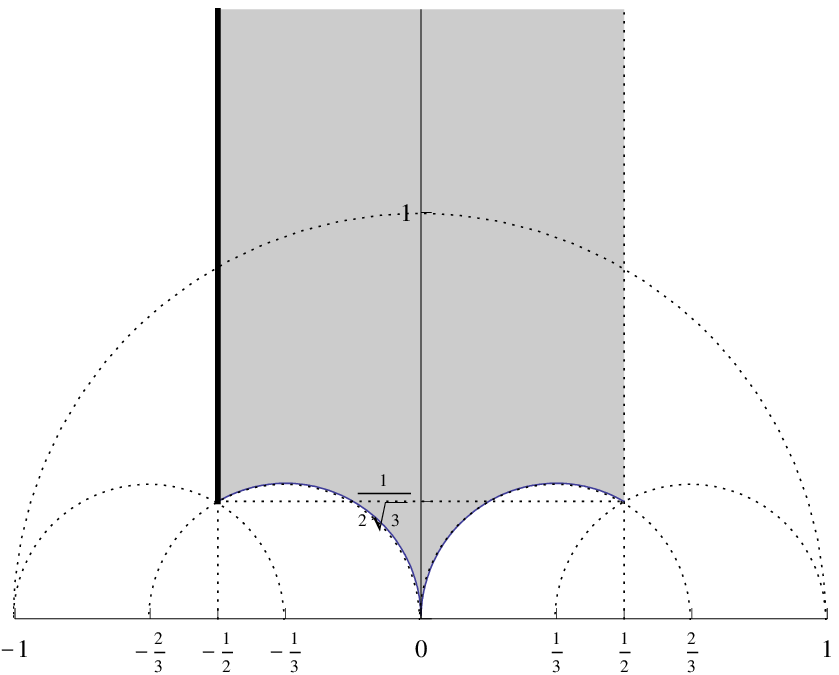}}
\end{center}
\caption{Location of the zeros of the Eisenstein series}
\end{figure}

\paragraph{\bf Location of the zeros of Hecke type Faber Polynomial}
For $m \leqslant 200$, we can prove that all of the zeros of $F_{m, 3}$ lie on the lower arcs of $\partial \mathbb{F}_3$ by numerical calculation.\\

\clearpage

\section{Level $4$}

We have $\Gamma_0(4)+=\Gamma_0(4)+4=\Gamma_0^{*}(4)$ and $\Gamma_0(4)-=\Gamma_0(4)$. We have $W_4 = \left(\begin{smallmatrix} 0 & -1/2 \\ 2 & 0 \end{smallmatrix}\right)$ and define $W_{4-, 2} = \left(\begin{smallmatrix} -1 & -1 \\ 2 & 1 \end{smallmatrix}\right)$ and $W_{4+, 2} = \left(\begin{smallmatrix} -1/\sqrt{2} & -3/(2 \sqrt{2}) \\ \sqrt{2} & 1/\sqrt{2} \end{smallmatrix}\right)$.\\

\subsection{$\Gamma_0^{*}(4)$}

We have $\Gamma_0^{*}(4) = T_{1/2}^{-1} \: \Gamma_0(2) \: T_{1/2}$ and $\Gamma_0^{*}(4) = \langle \left( \begin{smallmatrix} 1 & 1 \\ 0 & 1 \end{smallmatrix} \right), \: W_4 \rangle$.\\

\paragraph{\bf Location of the zeros of the Eisenstein series}
Since $(W_{4+, 2})^{- 1} \Gamma_0^{*}(4) W_{4+, 2} = \Gamma_0^{*}(4)$, we have
\begin{equation}
E_{k, 4+4}^{-1/2}(\gamma_{-1/2} z) = (\sqrt{2} z + 1/\sqrt{2})^k E_{k, 4+4}^{\infty}(z).
\end{equation}
Furthermore, we have
\begin{equation*}
E_{k, 4+4}^{-1/2} (i \tan(\theta /2) / 2) = ((e^{i \theta} + 1) / \sqrt{2})^k E_{k, 4+4}^{\infty}(e^{i \theta} / 2).
\end{equation*}

Now, recall that $\Gamma_0^{*}(4) = T_{1/2}^{-1} \: \Gamma_0(2) \: T_{1/2}$. Then, for $k \leqslant 1000$, since we can prove that all of the zeros of $E_{k, 2}^{\infty}$ lie on the lower arcs of $\partial \mathbb{F}_2$ by numerical calculation, we have all of the zeros of $E_{k, 4+4}^{\infty}$ in the lower arcs of $\partial \mathbb{F}_{4+4}$.

\begin{figure}[hbtp]
\begin{center}
{{$E_{k, 4+4}^{\infty}$}\includegraphics[width=1.5in]{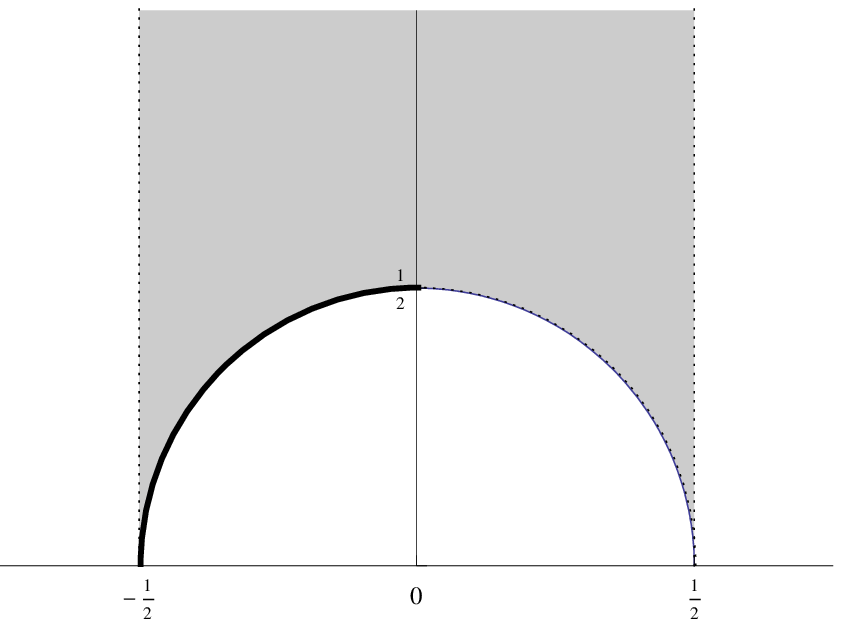}}
\quad \quad
{{$E_{k, 4+4}^{-1/2}$}\includegraphics[width=1.5in]{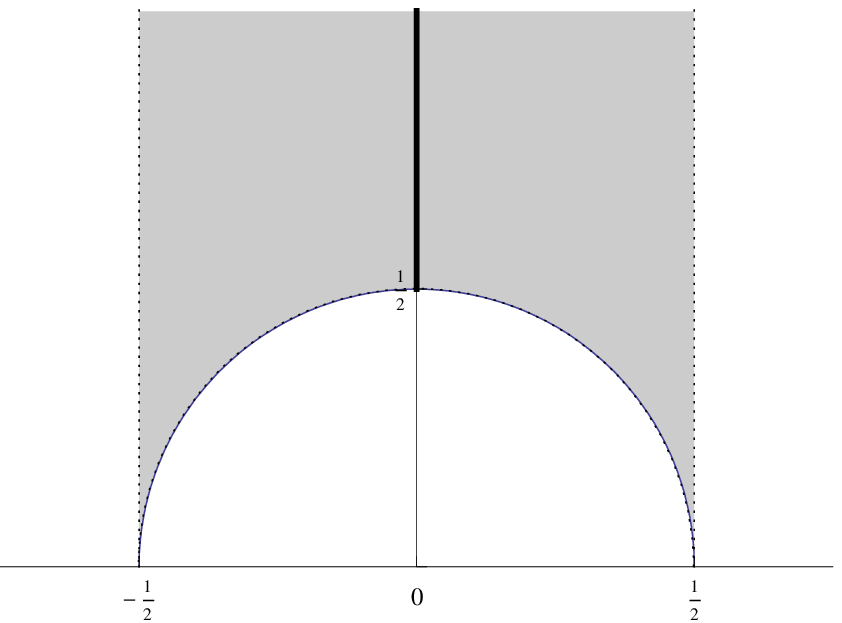}}
\end{center}
\caption{Location of the zeros of the Eisenstein series}
\end{figure}

\paragraph{\bf Location of the zeros of Hecke type Faber Polynomial}
For $m \leqslant 200$, we can prove that all of the zeros of $F_{m, 4+4}$ lie on the lower arcs of $\partial \mathbb{F}_{4+4}$ by numerical calculation.\\

\subsection{$\Gamma_0(4)$}

We have $\Gamma_0(4) = V_2^{-1} \Gamma(2) V_2$ and $\Gamma_0(4) = \langle -I, \: \left( \begin{smallmatrix} 1 & 1 \\ 0 & 1 \end{smallmatrix} \right), \: \left( \begin{smallmatrix} 1 & 0 \\ 4 & 1 \end{smallmatrix} \right) \rangle$. Furthermore, we have $\gamma_0 = W_4$ and $\gamma_{-1/2} = W_{4-, 2}$.\\

\paragraph{\bf Location of the zeros of the Eisenstein series}
Since $W_4^{- 1} \Gamma_0(4) W_4 = \gamma_{-1/2}^{- 1} \Gamma_0(4) \gamma_{-1/2} = \Gamma_0(4)$, we have
\begin{gather}
E_{k, 4}^0(W_4 z) = (2 z)^k E_{k, 4}^{\infty}(z),\\
E_{k, 4}^{-1/2}(\gamma_{-1/2} z) = (2 z + 1)^k E_{k, 4}^{\infty}(z).
\end{gather}
Furthermore, we have
\begin{gather*}
E_{k, 4}^0 (- 1/2 + i \tan(\theta /2) / 2) = ((e^{i \theta} - 1) / 2)^k E_{k, 4}^{\infty}((e^{i \theta} - 1) / 4),\\
E_{k, 4}^{-1/2} (i \tan(\theta /2) / 2) = ((e^{i \theta} + 1) / 2)^k E_{k, 4}^{\infty}((e^{i \theta} - 1) / 4).
\end{gather*}

Now, recall that $E_{k, 4}^{\infty}(z) = E_{k, 2}^{\infty}(2 z)$, then $E_{k, 4}^{\infty}(z)$ has $\lfloor k / 4 \rfloor - 1$ zeros in $\{ |z| = 1/4, \; -1/4 < Re(z) < 0 \}$, and $v_{-1/4 + i/4}(E_{k, 4}^{\infty}) = 1$ for $k \equiv 2 \pmod{4}$. Moreover, by the transformation with $W_{4-, 2}$ for $E_{k, 2}^{\infty}$, we have
\begin{equation*}
E_{k, 4}^{\infty}((e^{i \theta} - 1) / 4) = E_{k, 2}^{\infty}((e^{i \theta} - 1) / 2)
 = e^{i k (\pi - \theta)} E_{k, 2}^{\infty}((e^{i (\pi - \theta)} - 1) / 2) = e^{i k (\pi - \theta)} E_{k, 4}^{\infty}((e^{i (\pi - \theta)} - 1) / 4).
\end{equation*}

For $k \leqslant 1000$, we can prove that all of the zeros of $E_{k, 2}^{\infty}$ lie on the lower arcs of $\partial \mathbb{F}_2$ by numerical calculation, then we have all of the zeros of $E_{k, 4}^{\infty}$ in the lower arcs of $\partial \mathbb{F}_4$.

\begin{figure}[hbtp]
\begin{center}
{{$E_{k, 4}^{\infty}$}\includegraphics[width=1.5in]{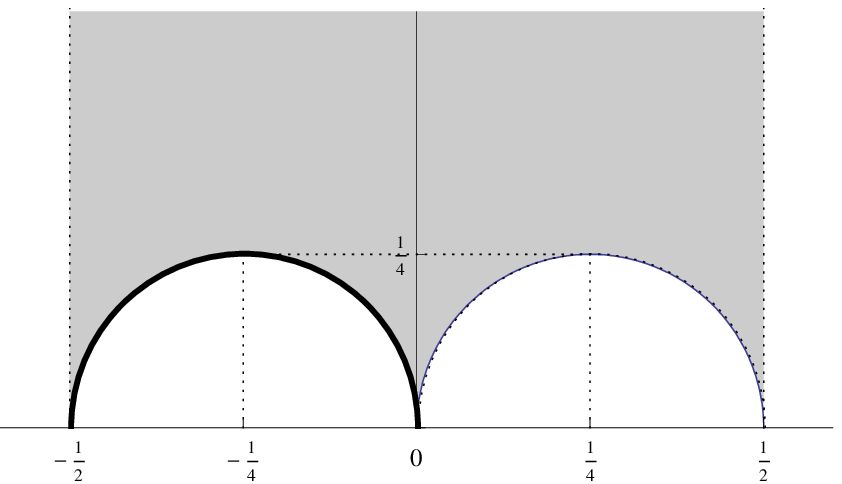}}
\quad \quad
{{$E_{k, 4}^{0}$}\includegraphics[width=1.5in]{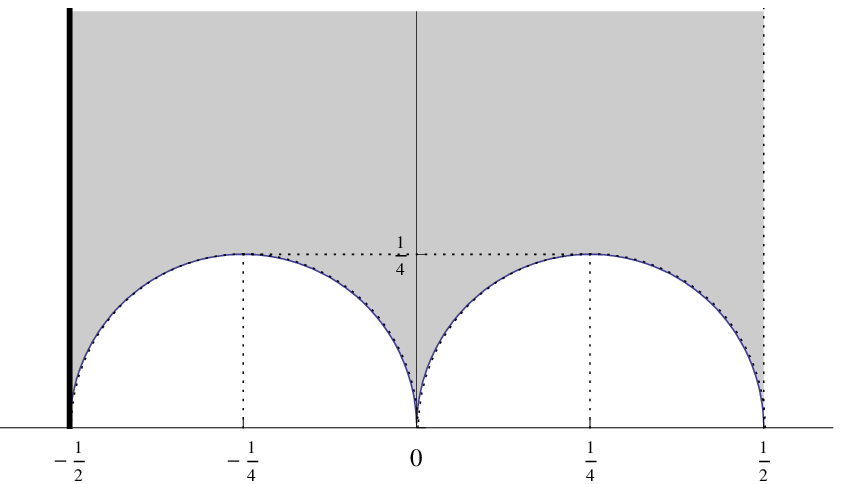}}
\quad \quad
{{$E_{k, 4}^{-1/2}$}\includegraphics[width=1.5in]{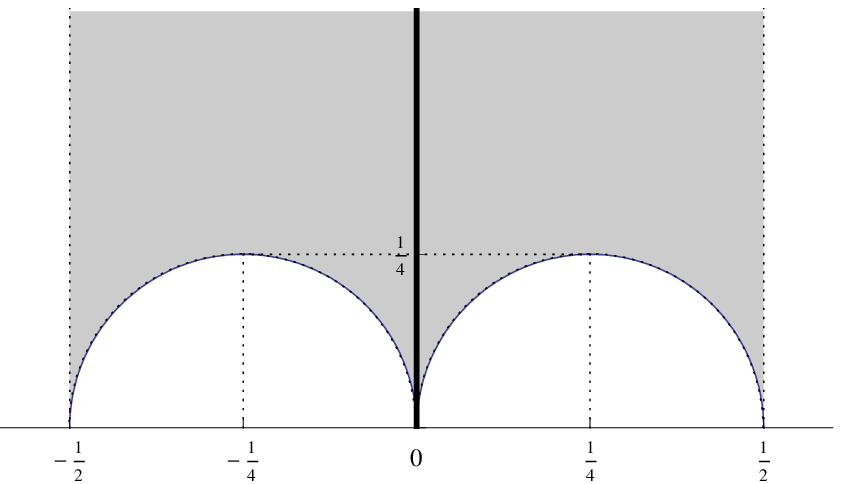}}
\end{center}
\caption{Location of the zeros of the Eisenstein series}
\end{figure}
\newpage

\paragraph{\bf Location of the zeros of Hecke type Faber Polynomial}
For $m \leqslant 200$, we can prove that all of the zeros of $F_{m, 4}$ lie on the lower arcs of $\partial \mathbb{F}_4$ by numerical calculation.\\

\clearpage

\section{Level $5$}

We have $\Gamma_0(5)+=\Gamma_0^{*}(5)$ and $\Gamma_0(5)-=\Gamma_0(5)$. We have $W_5 = \left(\begin{smallmatrix}0&-1 / \sqrt{5}\\ \sqrt{5}&0\end{smallmatrix}\right)$.\\

\subsection{$\Gamma_0^{*}(5)$}

We have $\Gamma_0^{*}(5) = \langle \left( \begin{smallmatrix} 1 & 1 \\ 0 & 1 \end{smallmatrix} \right), \: W_5, \: \left( \begin{smallmatrix} 3 & 1 \\ 5 & 2 \end{smallmatrix} \right) \rangle$.\\

\paragraph{\bf Location of the zeros of the Eisenstein series}
In \cite{SJ2}, the present author proved that all of the zeros of $E_{k, 5+}$ lie on the lower arcs of $\partial \mathbb{F}_{5+}$ if $4 \mid k$, and we prove all but at most one of the zeros of $E_{k, 5+}$ lie there if $4 \nmid k$. Furthermore, let $\alpha_5 \in [0, \pi]$ be the angle which satisfies $\tan\alpha_5 = 2$, and let $\alpha_{5, k} \in [0, \pi]$ be the angle which satisfies $\alpha_{5, k} \equiv k (\pi / 2 + \alpha_5) / 2 \pmod{\pi}$. We prove that all of the zeros of $E_{k, 5+}(z)$ in $\mathbb{F}^{*}(5)$ are on the lower arcs of $\partial \mathbb{F}_{5+}$ for $4 \mid k$ if $\alpha_{5, k} < (116/180) \pi$ or $(117/180) \pi < \alpha_{5, k}$.

In addition, for $k \leqslant 2500$, we can prove that all of the zeros of $E_{k, 5+}$ lie on the lower arcs of $\partial \mathbb{F}_{5+}$ by numerical calculation.

\begin{figure}[hbtp]
\begin{center}
\includegraphics[width=1.5in]{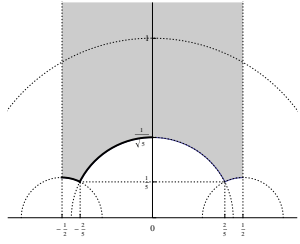}
\end{center}
\caption{Location of the zeros of the Eisenstein series}
\end{figure}

\paragraph{\bf Location of the zeros of Hecke type Faber Polynomial}
For $m = 1$ and $3 \leqslant m \leqslant 200$, we can prove that all of the zeros of $F_{m, 5+}$ lie on the lower arcs of $\partial \mathbb{F}_{5+}$ by numerical calculation. On the other hand, by numerical calculation, we can prove that all but one of the zeros of $F_{2, 5+}$ lie on the lower arcs of $\partial \mathbb{F}_{5+}$, and one of the zeros of $F_{2, 5+}$ lies on $\partial \mathbb{F}_{5+}$ but does not on the lower arcs.\\

\subsection{$\Gamma_0(5)$}

We have $\Gamma_0(5) = \langle \left( \begin{smallmatrix} 1 & 1 \\ 0 & 1 \end{smallmatrix} \right), \: \left( \begin{smallmatrix} 1 & 0 \\ 5 & 1 \end{smallmatrix} \right), \: \left( \begin{smallmatrix} 3 & 1 \\ 5 & 2 \end{smallmatrix} \right) \rangle$ and $\gamma_0 = W_5$.\\

\paragraph{\bf Location of the zeros of the Eisenstein series}
Since $W_5^{- 1} \Gamma_0(5) W_5 = \Gamma_0(5)$, we have
\begin{equation}
E_{k, 5}^0(W_5 z) = (\sqrt{5} z)^k E_{k, 5}^{\infty}(z).
\end{equation}
Furthermore, we have
\begin{align*}
E_{k, 5}^0 (- 1/2 + i / (2 \tan\theta/2)) &= ((e^{i \theta} - 1) / \sqrt{5})^k E_{k, 5}^{\infty}((e^{i \theta} - 1) / 5),\\
E_{k, 5}^0 ((e^{i \theta'} + 2) / 3) &= ((e^{i \theta} - 2) / \sqrt{5})^k E_{k, 5}^{\infty}((e^{i \theta} - 2) / 5),\\
E_{k, 5}^0 ((e^{i (\pi - \theta')} - 2) / 3) &= ((2 e^{i \theta} - 1) / \sqrt{5})^k E_{k, 5}^{\infty}((e^{i \theta} + 2) / 5),
\end{align*}
where $e^{i \theta'} = (4 - 5 \cos\theta + 3 i \sin\theta) / (5 - 4 \cos\theta)$.
\begin{figure}[hbtp]
\begin{center}
{{$E_{k, 5}^{\infty}$}\includegraphics[width=1.5in]{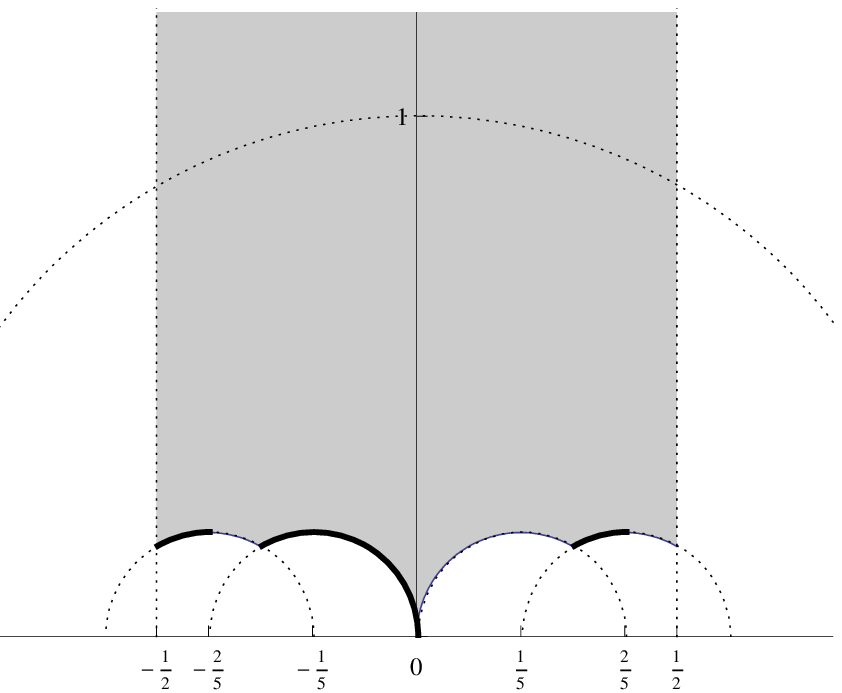}}
\quad \quad
{{$E_{k, 5}^{0}$}\includegraphics[width=1.5in]{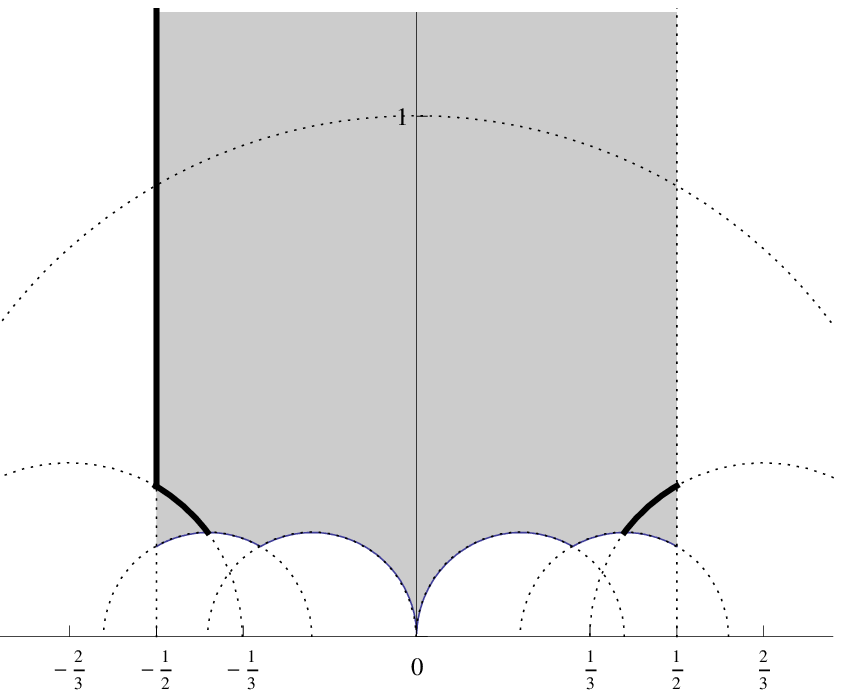}}
\end{center}
\caption{Neighborhood of location of the zeros of the Eisenstein series}
\end{figure}

We can verify whether the zeros lie on $\partial \mathbb{F}_5$ if $J_5$ takes real value there. However, $J_5$ does not take real value, then all we can do is to observe the graphs.

\begin{figure}[hbtp]
\begin{center}
{{\small Lower arcs of $\partial \mathbb{F}_5$}\includegraphics[width=2.5in]{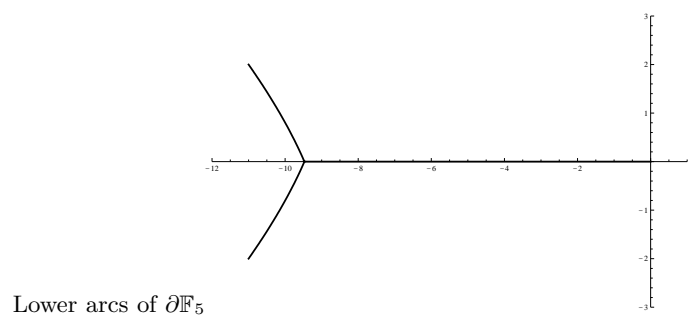}}
\end{center}
\caption{Image by $J_5$}\label{Im-J5B}
\end{figure}

\begin{figure}[hbtp]
\begin{center}
\includegraphics[width=6.3in]{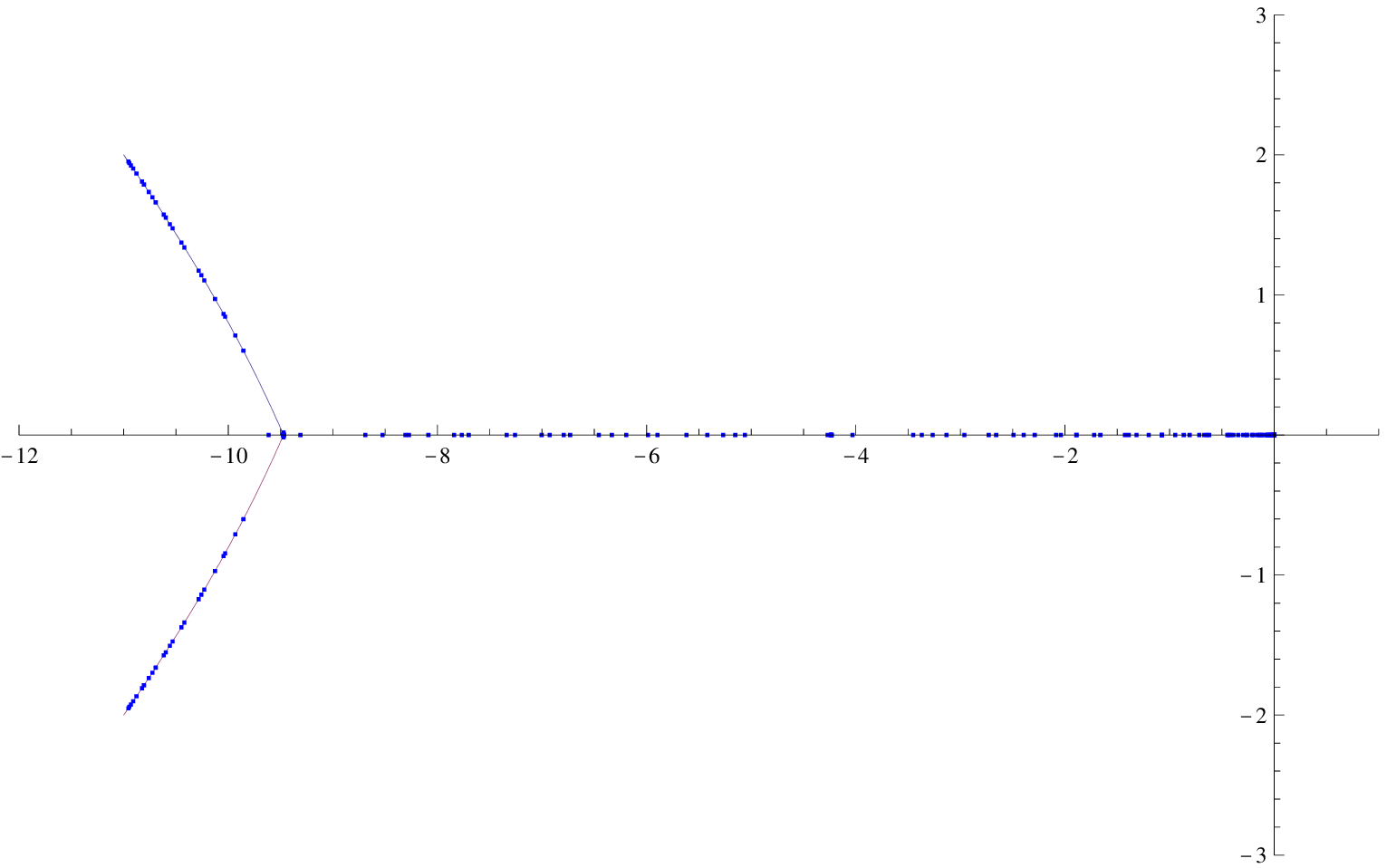}\\
The zeros of $E_{k, 5}^{\infty}$ for $4 \leqslant k \leqslant 40$\\
\includegraphics[width=6.3in]{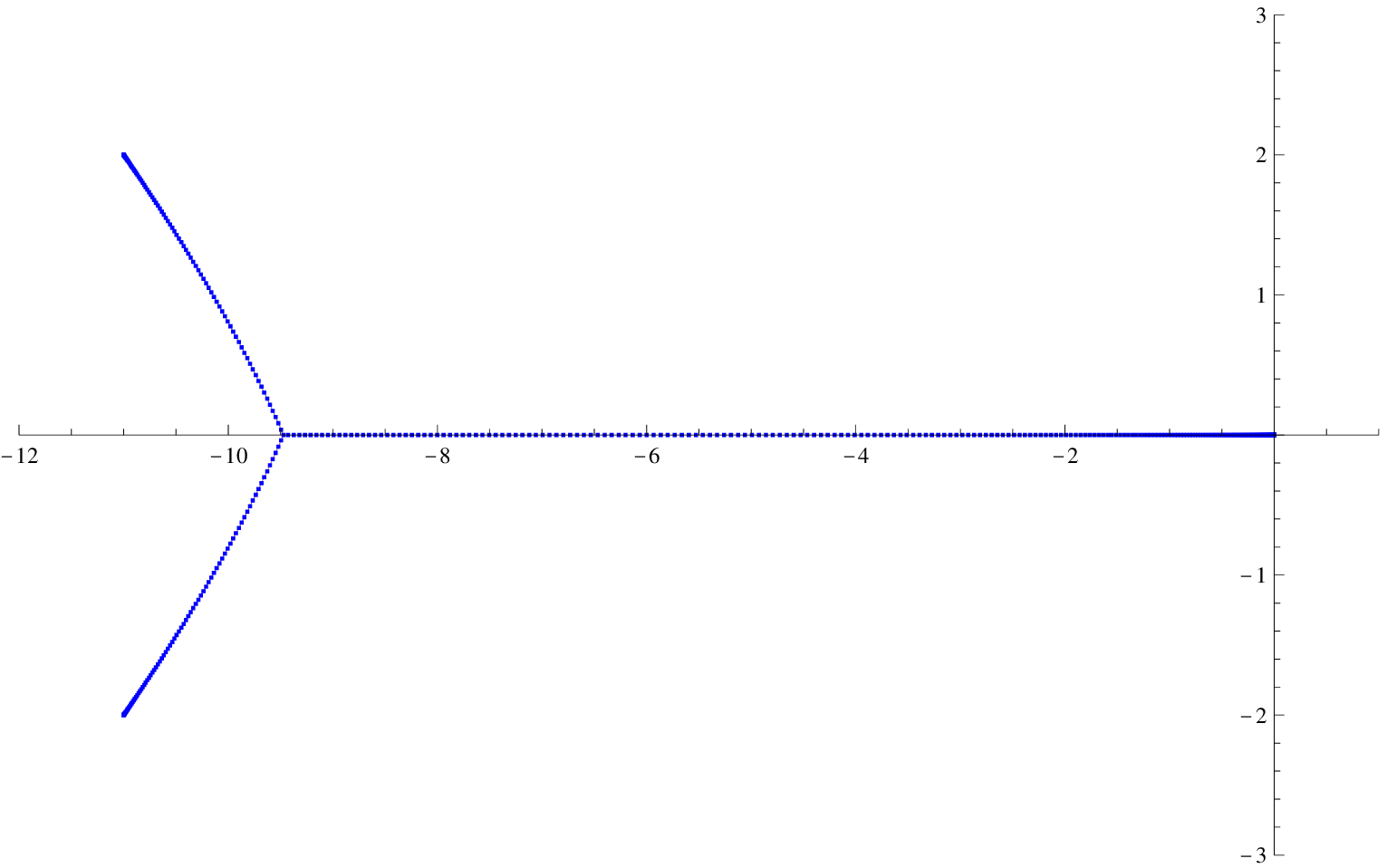}\\
The zeros of $E_{1000, 5}^{\infty}$
\end{center}
\caption{Image by $J_5$}\label{Im-J5Bz}
\end{figure}

Now, we can observe that some zeros of $E_{k, 5}^{\infty}$ do not lie on the lower arcs of $\partial \mathbb{F}_5$ for small weight $k$ by numerical calculation, but they seems to lie on $\partial \mathbb{F}_5$ except for lower arcs. However, when the weight $k$ increases, then the location of the zeros seems to approach to lower arcs of $\partial \mathbb{F}_5$. (see Figure \ref{Im-J5Bz})\\

\begin{figure}[hbtp]
\begin{center}
\includegraphics[width=6.3in]{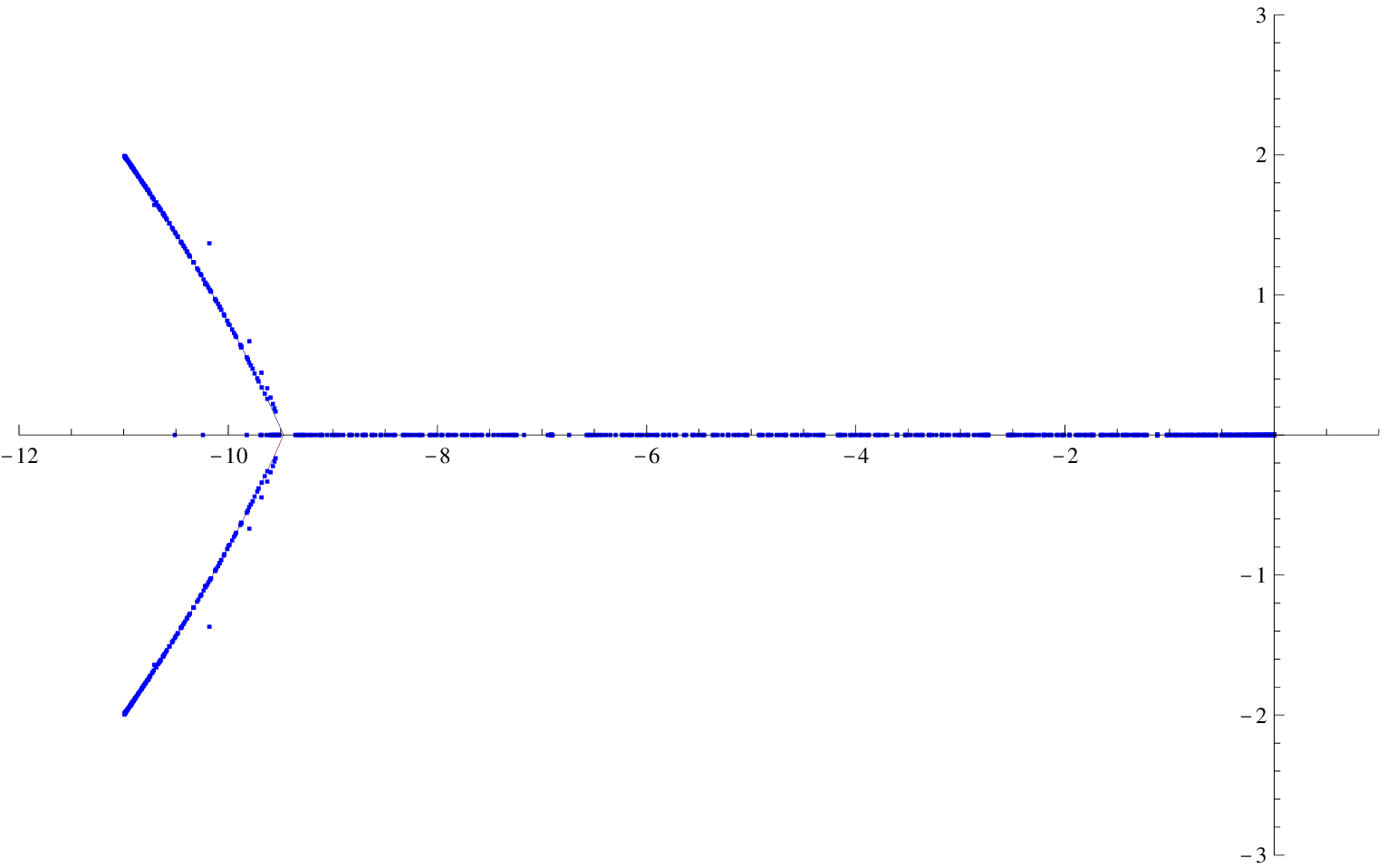}\\
The zeros of $F_{m, 5}$ for $m \leqslant 40$\\
\includegraphics[width=6.3in]{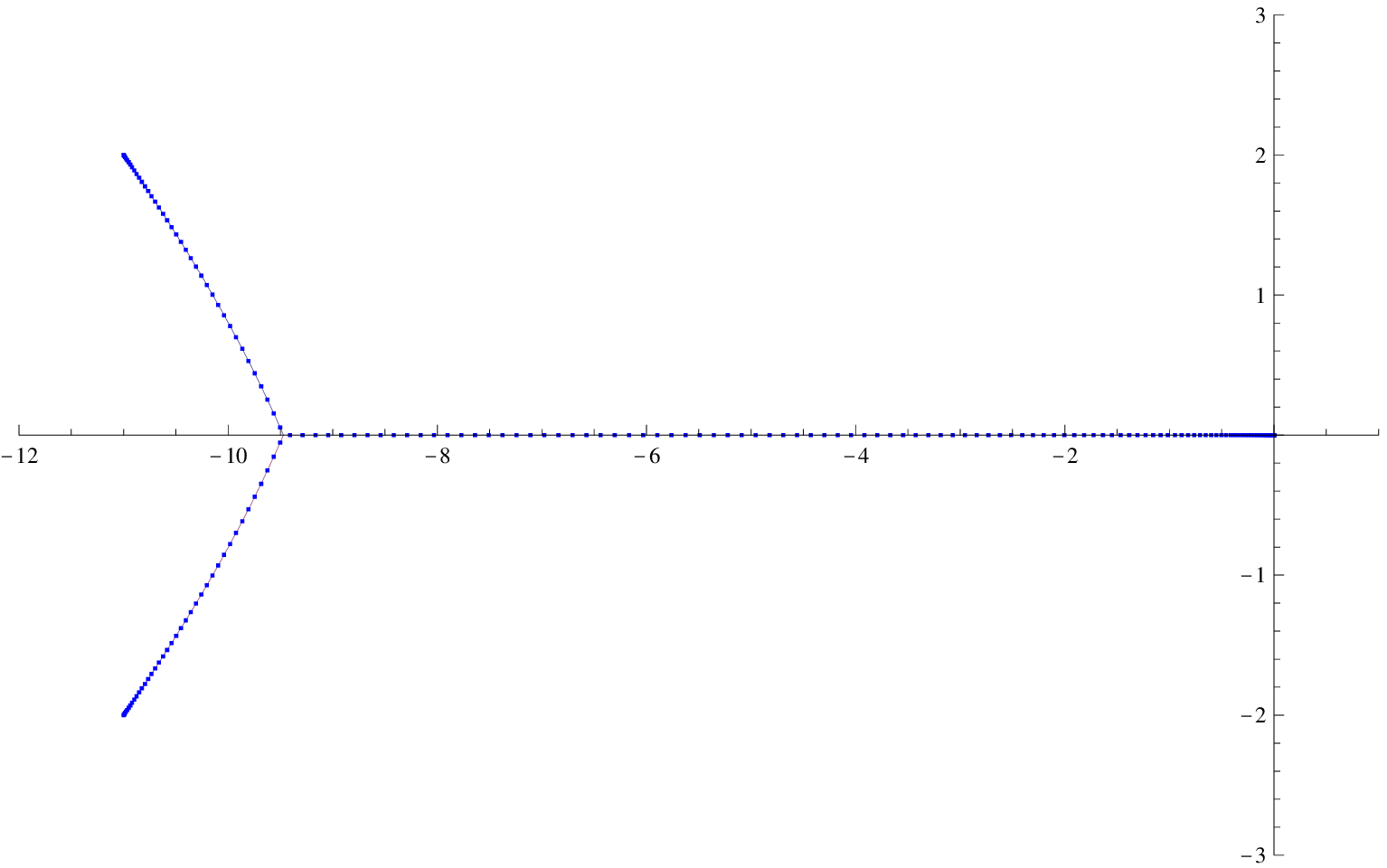}\\
The zeros of $F_{200, 5}$
\end{center}
\caption{Image by $J_5$}\label{Im-J5Bhz}
\end{figure}

\paragraph{\bf Location of the zeros of Hecke type Faber Polynomial}
Similarly to the Eisenstein series, we can observe that some zeros of $F_{m, 5}$ do not lie on the lower arcs of $\partial \mathbb{F}_5$ for small weight $m$ by numerical calculation. However, when the weight $m$ increases, then the location of the zeros seems to approach to lower arcs of $\partial \mathbb{F}_5$. (see Figure \ref{Im-J5Bhz})
 \clearpage

\section{Level $6$}

We have $\Gamma_0(6)+$, $\Gamma_0(6)+6=\Gamma_0^{*}(6)$, $\Gamma_0(6)+3$, $\Gamma_0(6)+2$, and $\Gamma_0(6)-=\Gamma_0(6)$. We have $W_6 = \left(\begin{smallmatrix}0 & - 1 / \sqrt{6}\\ \sqrt{6} & 0\end{smallmatrix}\right)$, $W_{6, 2} := \left(\begin{smallmatrix} -\sqrt{2} & -1/\sqrt{2}\\ 3\sqrt{2} & \sqrt{2}\end{smallmatrix}\right)$, and $W_{6, 3} := \left(\begin{smallmatrix} -\sqrt{3} & -2/\sqrt{3}\\ 2\sqrt{3} & \sqrt{3}\end{smallmatrix}\right)$.\\

\subsection{$\Gamma_0(6)+$}

We have $\Gamma_0(6)+ = \langle \left( \begin{smallmatrix} 1 & 1 \\ 0 & 1 \end{smallmatrix} \right), \: W_6, \: W_{6, 3} \rangle$.\\

\paragraph{\bf Location of the zeros of the Eisenstein series}
For $k \leqslant 1000$, we can prove that all of the zeros of $E_{k, 6+}$ lie on the lower arcs of $\partial \mathbb{F}_{6+}$ by numerical calculation.
\begin{figure}[hbtp]
\begin{center}
\includegraphics[width=1.5in]{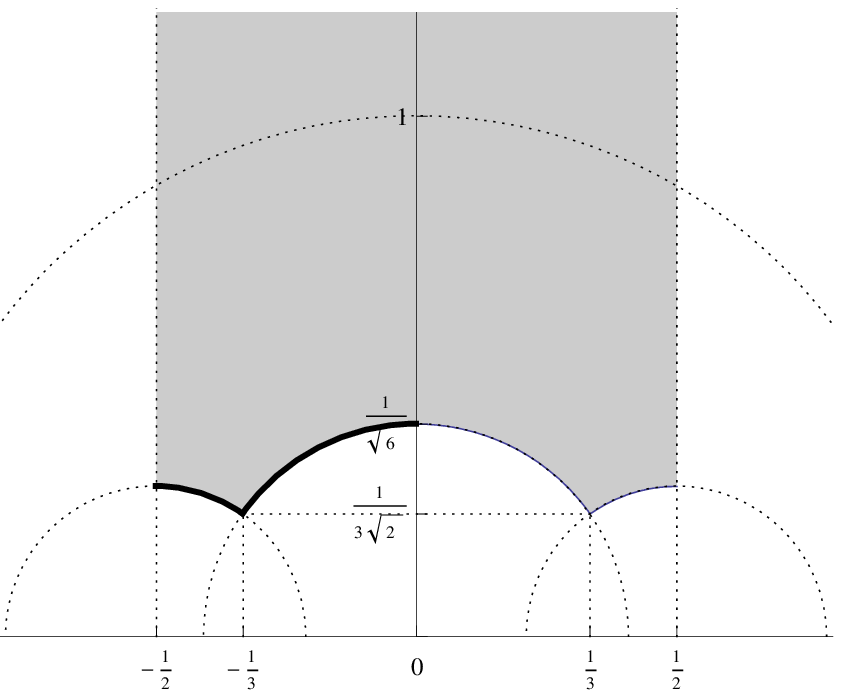}
\end{center}
\caption{Location of the zeros of the Eisenstein series}
\end{figure}

\paragraph{\bf Location of the zeros of Hecke type Faber Polynomial}
For $m \leqslant 200$, we can prove that all of the zeros of $F_{m, 6+}$ lie on the lower arcs of $\partial \mathbb{F}_{6+}$ by numerical calculation.\\

\subsection{$\Gamma_0(6)+6 = \Gamma_0^{*}(6)$}

We have $\Gamma_0^{*}(6) = \langle \left( \begin{smallmatrix} 1 & 1 \\ 0 & 1 \end{smallmatrix} \right), \: W_6, \: \left( \begin{smallmatrix} 5 & 2 \\ 12 & 5 \end{smallmatrix} \right) \rangle$ and $\gamma_{-1/2} = W_{6, 3}$\\

\paragraph{\bf Location of the zeros of the Eisenstein series}
Since $W_{6, 3}^{- 1} \Gamma_0^{*}(6) W_{6, 3} = \Gamma_0^{*}(6)$, we have
\begin{equation}
E_{k, 6+6}^{-1/2}(W_{6, 3} z) = (2 \sqrt{3} z + \sqrt{3})^k E_{k, 6+6}^{\infty}(z).
\end{equation}
Furthermore, we have
\begin{align*}
E_{k, 6+6}^{-1/2} (i \tan(\theta/2) / 2) &= ((e^{i \theta} + 1) / (2 \sqrt{3}))^k E_{k, 6+6}^{\infty}((e^{i \theta} - 5) / 12),\\
E_{k, 6+6}^{-1/2} (e^{i \theta'} / \sqrt{6}) &= (\sqrt{3} e^{i \theta} + \sqrt{2})^k E_{k, 6+6}^{\infty}(e^{i \theta} / \sqrt{6}),
\end{align*}
where $e^{i \theta'} = (- 2 \sqrt{6} - 5 \cos\theta + i \sin\theta) / (5 + 2 \sqrt{6} \cos\theta)$.

For $k \leqslant 750$, we can prove that all of the zeros of $E_{k, 6+6}^{\infty}$ lie on the lower arcs of $\partial \mathbb{F}_{6+6}$ by numerical calculation.

\begin{figure}[htbp]
\begin{center}
{{$E_{k, 6+6}^{\infty}$}\includegraphics[width=1.5in]{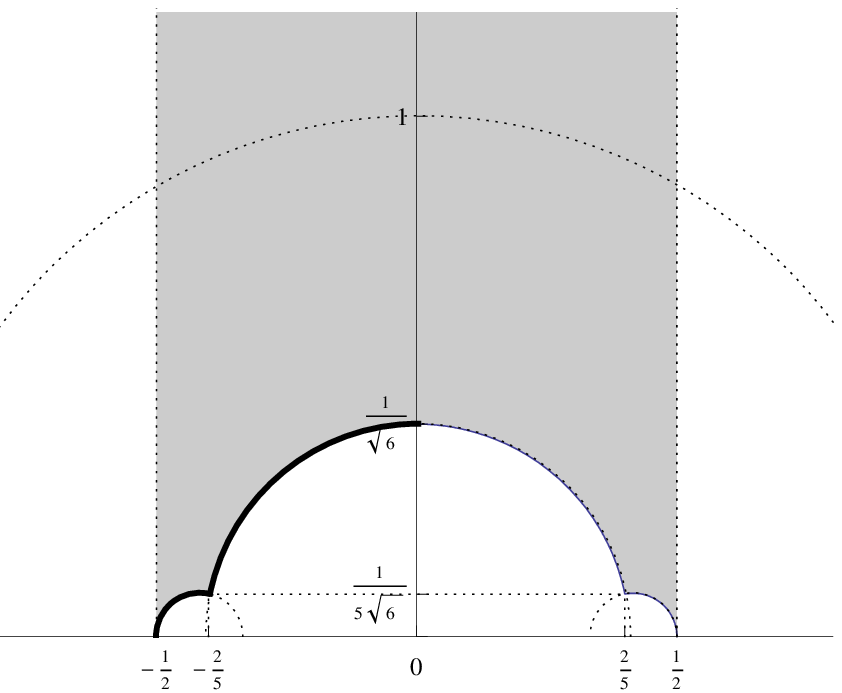}}
\quad \quad
{{$E_{k, 6+6}^{-1/2}$}\includegraphics[width=1.5in]{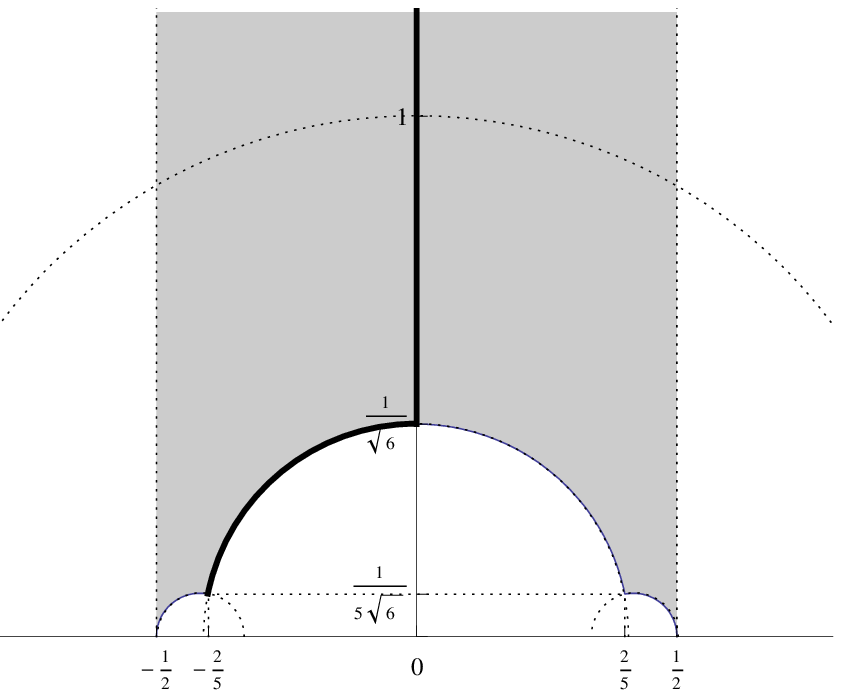}}
\end{center}
\caption{Location of the zeros of the Eisenstein series}
\end{figure}

\paragraph{\bf Location of the zeros of Hecke type Faber Polynomial}
For every odd integer $m \leqslant 200$, we can prove that all of the zeros of $F_{m, 6+6}$ lie on the lower arcs of $\partial \mathbb{F}_{6+6}$ by numerical calculation. On the other hand, by numerical calculation, for every even integer $m \leqslant 200$, we can prove that all but one of the zeros of $F_{m, 6+6}$ lie on the lower arcs of $\partial \mathbb{F}_{6+6}$, and one of the zeros of $F_{m, 6+6}$ lies on $\partial \mathbb{F}_{6+6}$ but does not on the lower arcs.

\newpage

\subsection{$\Gamma_0(6)+3$}

We have $\Gamma_0(6)+3 = \langle \left( \begin{smallmatrix} 1 & 1 \\ 0 & 1 \end{smallmatrix} \right), \: \left( \begin{smallmatrix} 1 & 0 \\ 6 & 1 \end{smallmatrix} \right), \: W_{6, 3} \rangle$ and $\gamma_0 = W_6$.\\

\paragraph{\bf Location of the zeros of the Eisenstein series}
Since $W_6^{- 1} (\Gamma_0(6)+3) W_6 = \Gamma_0(6)+3$, we have
\begin{equation}
E_{k, 6+3}^0(W_6 z) = (\sqrt{6} z)^k E_{k, 6+3}^{\infty}(z).
\end{equation}
Furthermore, we have
\begin{align*}
E_{k, 6+3}^0 (-1/2 + i / (2 \tan(\theta/2))) &= ((e^{i \theta} - 1) / \sqrt{6})^k E_{k, 6+3}^{\infty}((e^{i \theta} - 1) / 6),\\
E_{k, 6+3}^0 (e^{i \theta'} / (2 \sqrt{3}) - 1/2) &= ((\sqrt{3} e^{i \theta} -1) / \sqrt{2})^k E_{k, 6+3}^{\infty}(e^{i \theta} / (2 \sqrt{3}) - 1/2),
\end{align*}
where $e^{i \theta'} = (\sqrt{3} - 2 \cos\theta + i \sin\theta) / (2 - \sqrt{3} \cos\theta)$.

For $k \leqslant 600$, we can prove that all of the zeros of $E_{k, 6+3}^{\infty}$ lie on the lower arcs of $\partial \mathbb{F}_{6+3}$ by numerical calculation.

\begin{figure}[hbtp]
\begin{center}
{{$E_{k, 6+3}^{\infty}$}\includegraphics[width=1.5in]{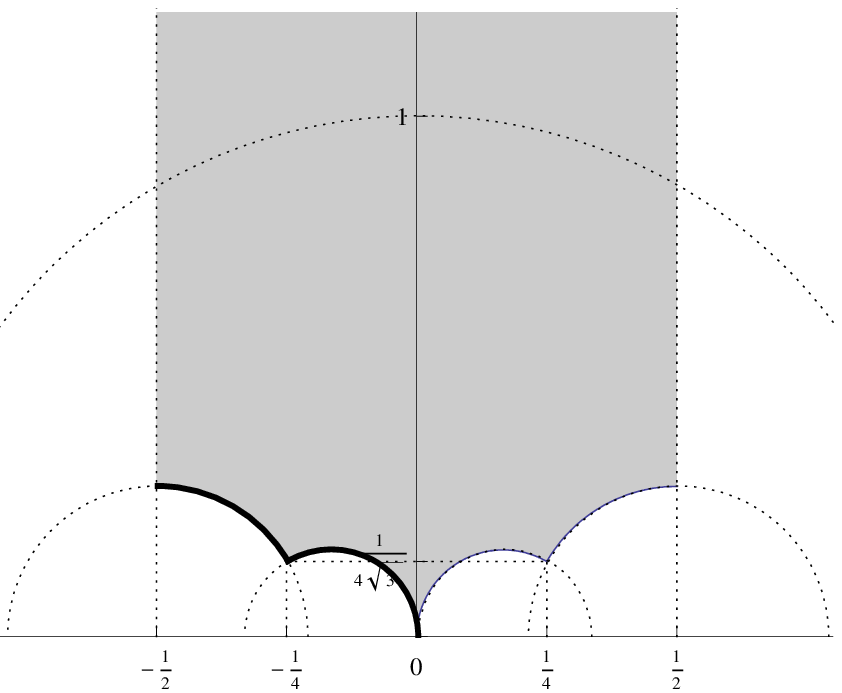}}
\quad \quad
{{$E_{k, 6+3}^0$}\includegraphics[width=1.5in]{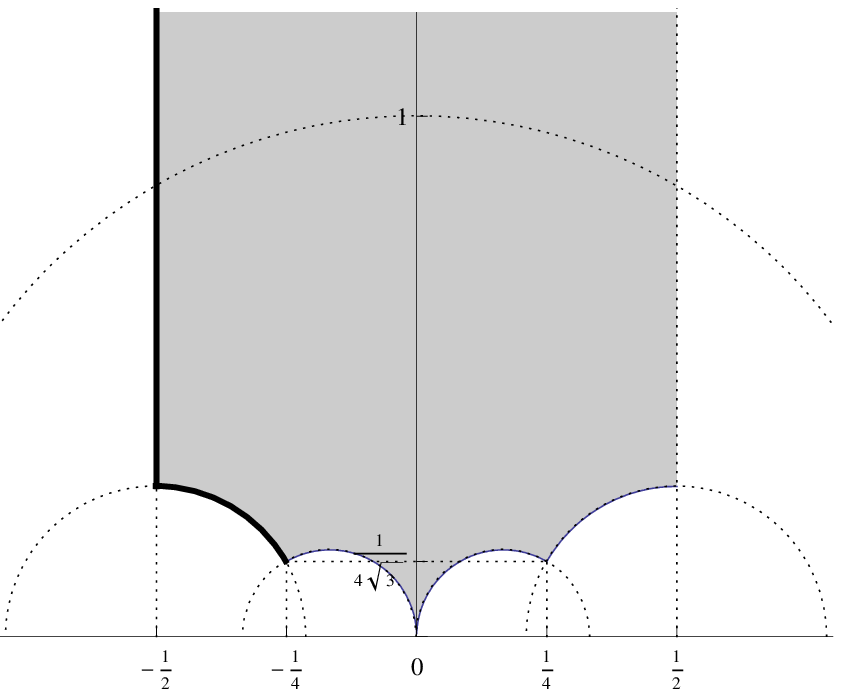}}
\end{center}
\caption{Location of the zeros of the Eisenstein series}
\end{figure}

\paragraph{\bf Location of the zeros of Hecke type Faber Polynomial}
For $m \leqslant 200$, we can prove that all of the zeros of $F_{m, 6+3}$ lie on the lower arcs of $\partial \mathbb{F}_{6+3}$ by numerical calculation.\\

\subsection{$\Gamma_0(6)+2$}

We have $\Gamma_0(6)+2 = \langle \left( \begin{smallmatrix} 1 & 1 \\ 0 & 1 \end{smallmatrix} \right), \: \left( \begin{smallmatrix} 1 & 0 \\ 6 & 1 \end{smallmatrix} \right), \: W_{6, 2} \rangle$ and $\gamma_0 = W_6$.\\

\paragraph{\bf Location of the zeros of Eisenstein series}

Since $W_6^{- 1} (\Gamma_0(6)+2) W_6 = \Gamma_0(6)+2$, we have
\begin{equation}
E_{k, 6+2}^0(W_6 z) = (\sqrt{6} z)^k E_{k, 6+2}^{\infty}(z).
\end{equation}
Furthermore, we have
\begin{align*}
E_{k, 6+2}^0 (- 1/2 + i / (2 \tan\theta/2)) &= ((e^{i \theta} - 1) / \sqrt{6})^k E_{k, 6+2}^{\infty}((e^{i \theta} - 1) / 6),\\
E_{k, 6+2}^0 (e^{i \theta'} / \sqrt{2} + 1) &= ((e^{i \theta} - \sqrt{2}) / \sqrt{3})^k E_{k, 6+2}^{\infty}(e^{i \theta} / (3 \sqrt{2}) - 1/3),\\
E_{k, 6+2}^0 (e^{i \theta''} / \sqrt{2} - 1) &= ((e^{i \theta} + \sqrt{2}) / \sqrt{3})^k E_{k, 6+2}^{\infty}(e^{i \theta} / (3 \sqrt{2}) + 1/3),
\end{align*}
where $e^{i \theta'} = (3 \cos\theta - 2 \sqrt{2} + i \sin\theta) / (2 \sqrt{2} - 3 \cos\theta)$ and $e^{i \theta''} = (3 \cos\theta + 2 \sqrt{2} + i \sin\theta) / (2 \sqrt{2} + 3 \cos\theta)$.
\begin{figure}[hbtp]
\begin{center}
{{$E_{k, 6+2}^{\infty}$}\includegraphics[width=1.5in]{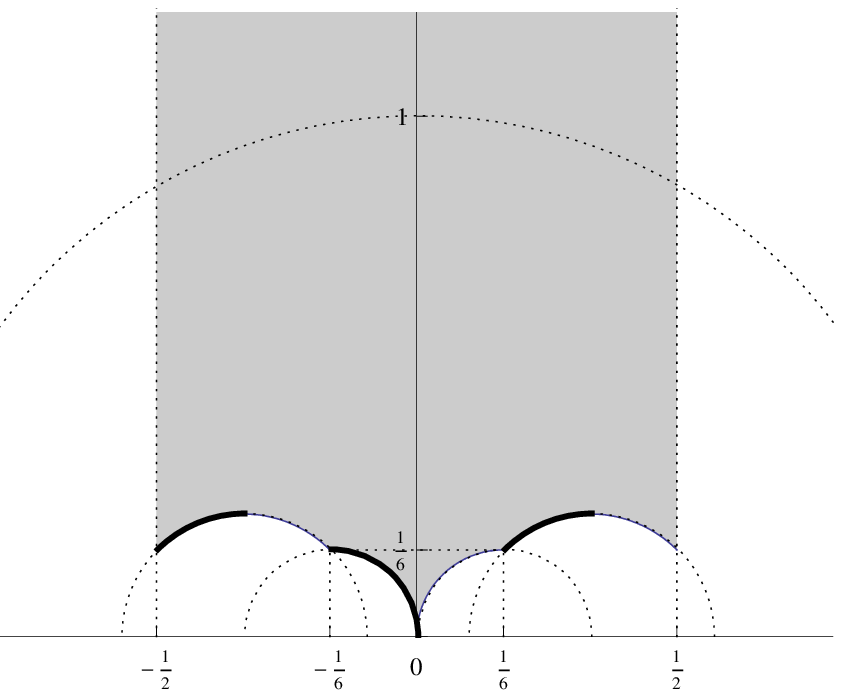}}
\quad \quad
{{$E_{k, 6+2}^{0}$}\includegraphics[width=1.5in]{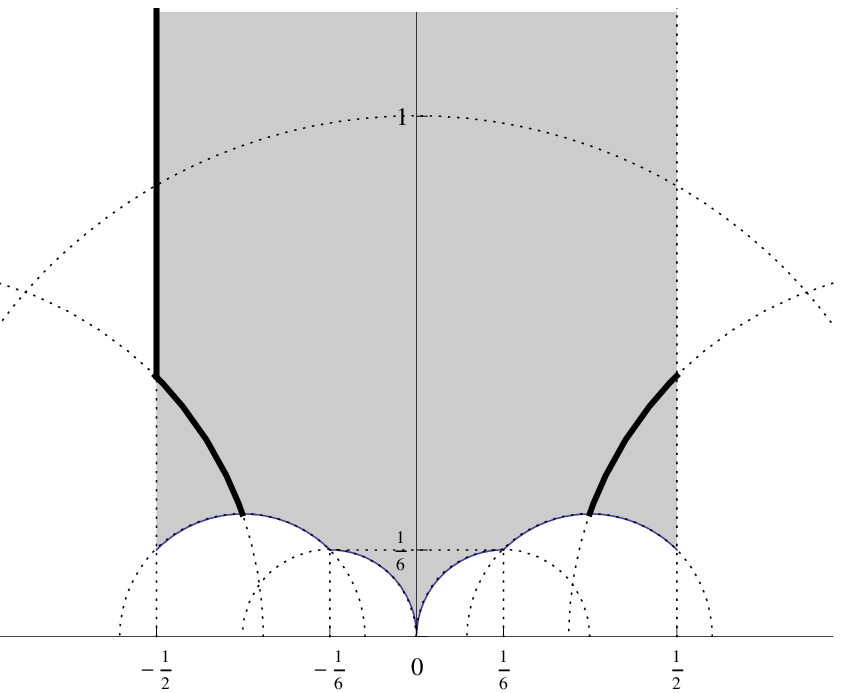}}
\end{center}
\caption{Neighborhood of location of the zeros of the Eisenstein series}
\end{figure}

\begin{figure}[hbtp]
\begin{center}
{{\small Lower arcs of $\partial \mathbb{F}_{6+2}$}\includegraphics[width=2.5in]{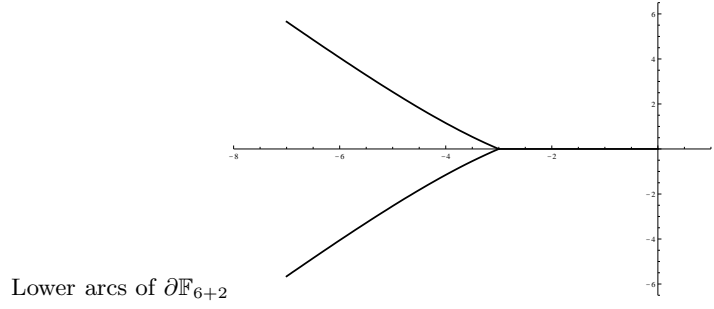}}
\end{center}
\caption{Image by $J_{6+2}$}\label{Im-J6D}
\end{figure}

\begin{figure}[hbtp]
\begin{center}
\includegraphics[width=6.3in]{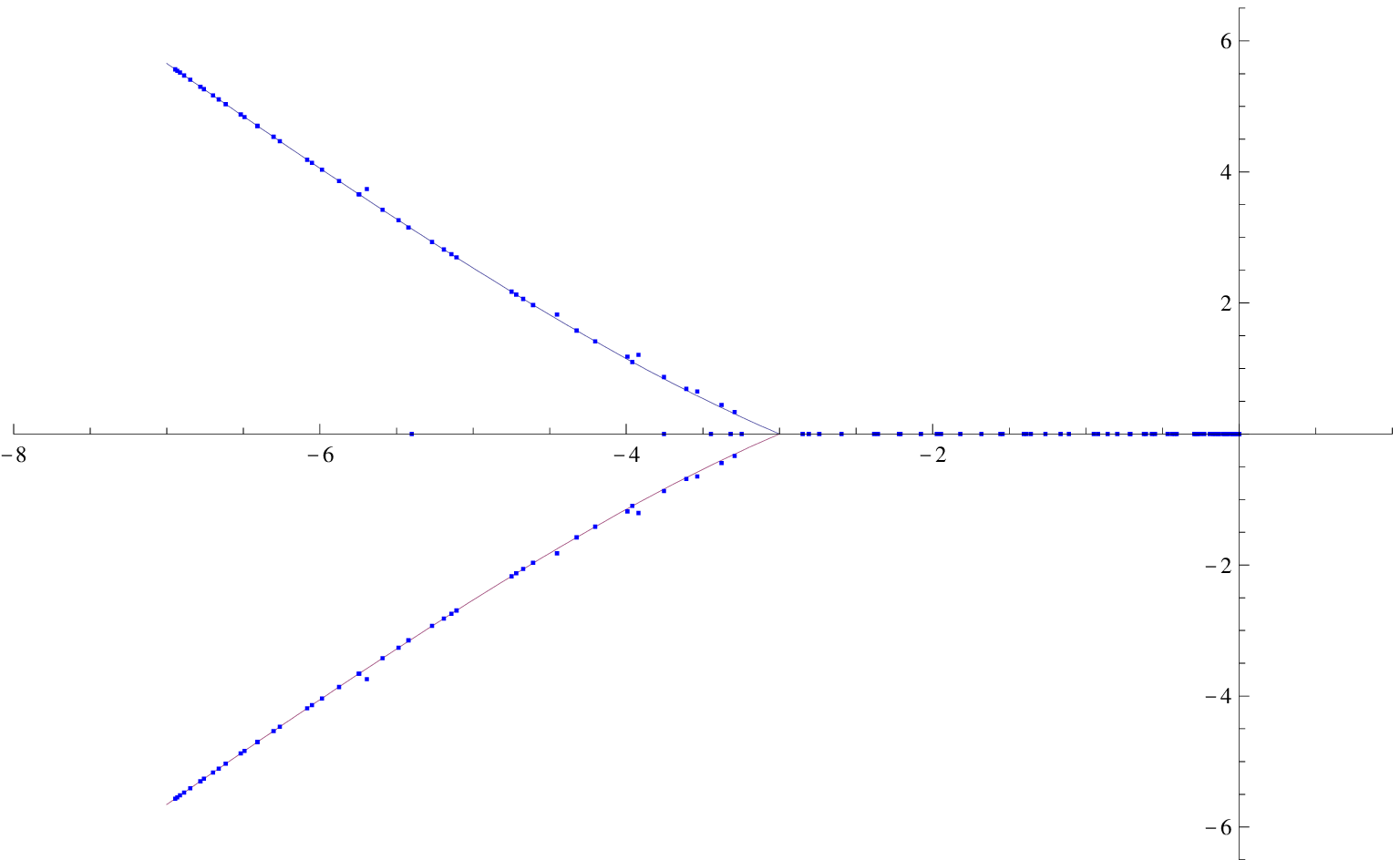}\\
The zeros of $E_{k, 6+2}^{\infty}$ for $4 \leqslant k \leqslant 40$\\
\includegraphics[width=6.3in]{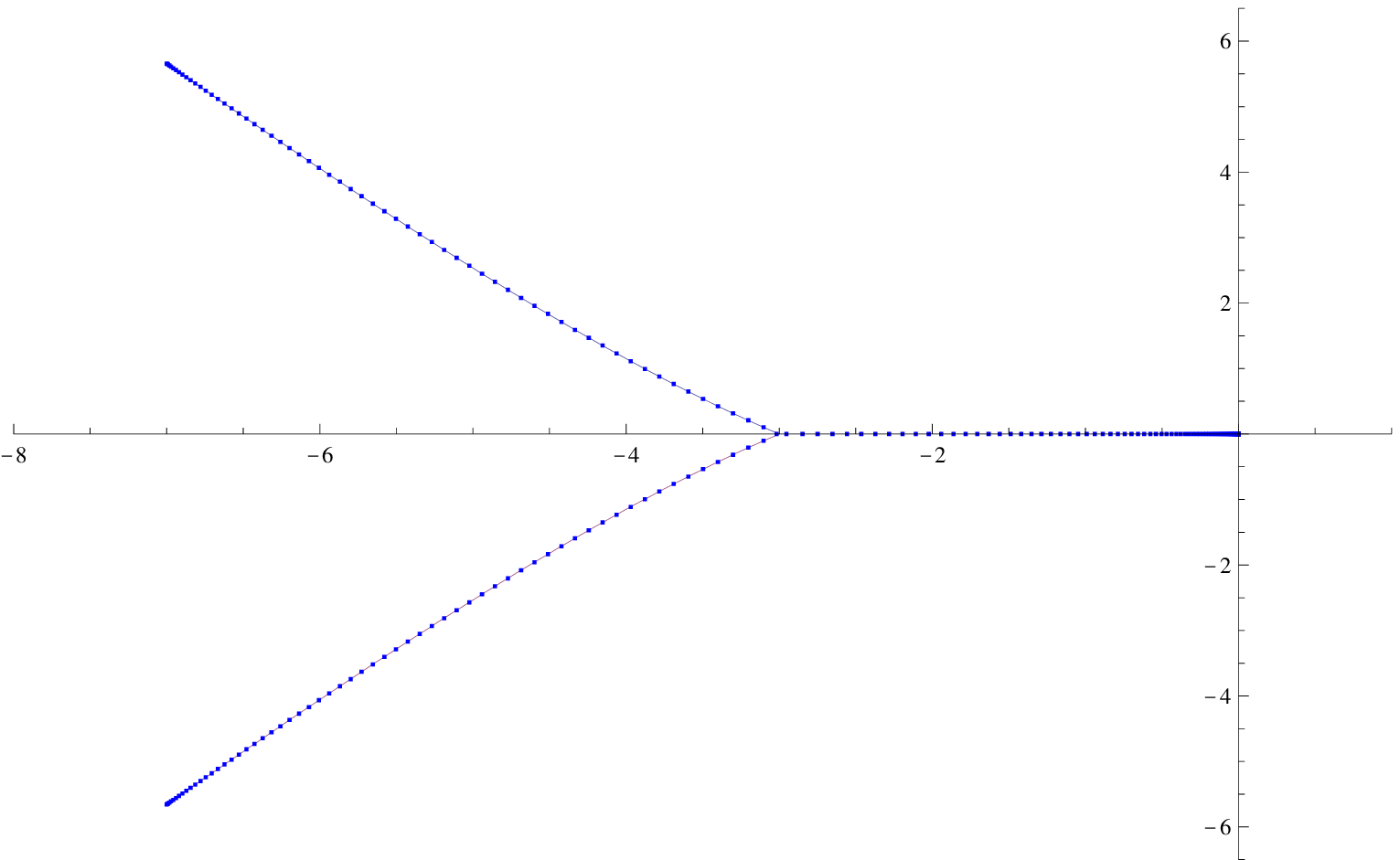}\\
The zeros of $E_{1000, 6+2}^{\infty}$
\end{center}
\caption{Image by $J_{6+2}$}\label{Im-J6Dz}
\end{figure}

Now, we can observe that the zeros of $E_{k, 6+2}^{\infty}$ do not lie on the arcs of $\partial \mathbb{F}_{6+2}$ for small weight $k$ by numerical calculation. However, when the weight $k$ increases, then the location of the zeros seems to approach to lower arcs of $\partial \mathbb{F}_{6+2}$. (See Figure \ref{Im-J6Dz})\\

\begin{figure}[hbtp]
\begin{center}
\includegraphics[width=6.3in]{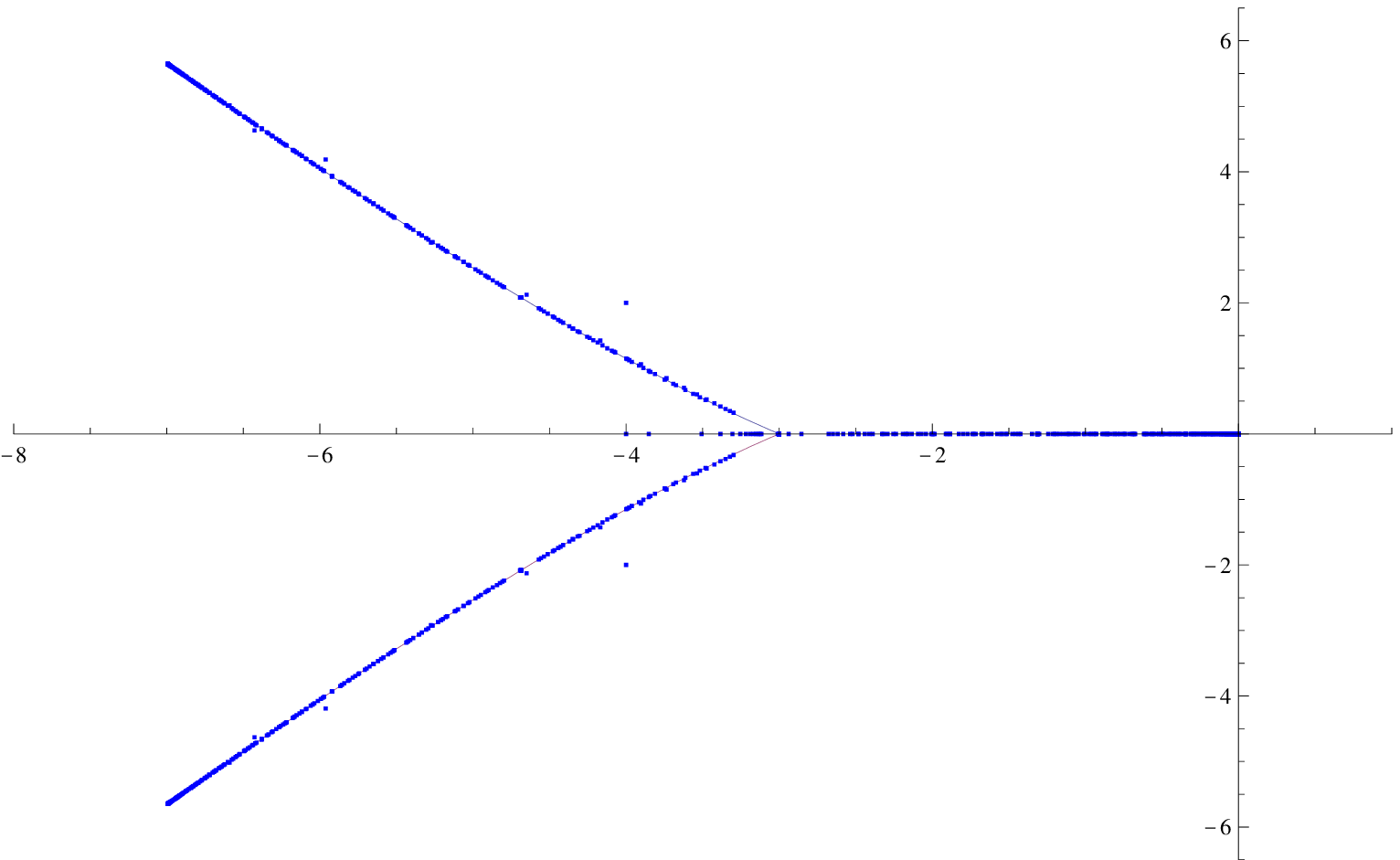}\\
The zeros of $F_{m, 6+2}$ for $m \leqslant 40$\\
\includegraphics[width=6.3in]{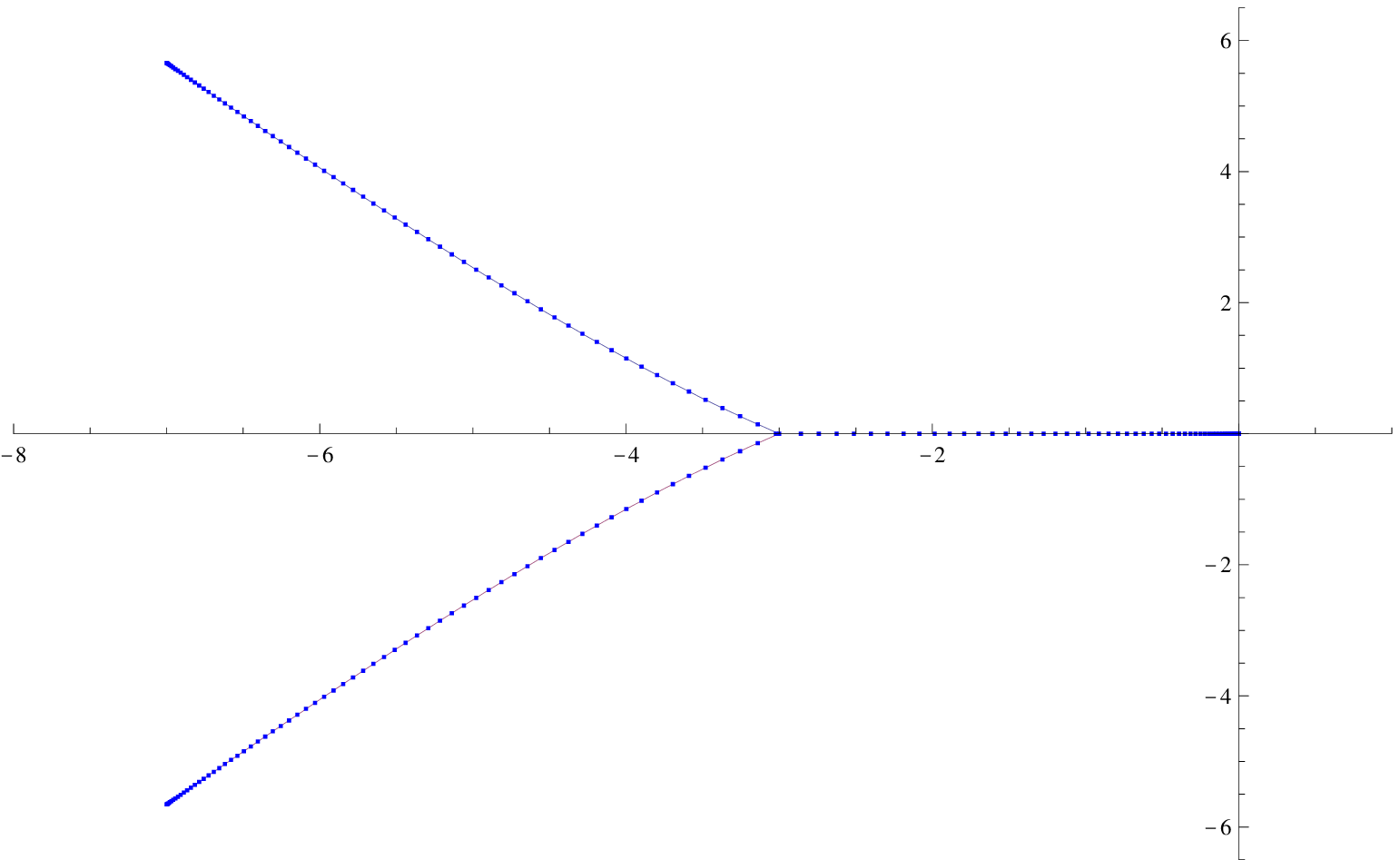}\\
The zeros of $F_{200, 6+2}$
\end{center}
\caption{Image by $J_{6+2}$}\label{Im-J6Dhz}
\end{figure}

\paragraph{\bf Location of the zeros of Hecke type Faber Polynomial}
We can observe that some zeros of $F_{m, 6+2}$ do not lie on the lower arcs of $\partial \mathbb{F}_{6+2}$ for small weight $m$ by numerical calculation. However, when the weight $m$ increases, then the location of the zeros seems to approach to lower arcs of $\partial \mathbb{F}_{6+2}$. (see Figure \ref{Im-J6Dhz})\\

\subsection{$\Gamma_0(6)$}

We have $\Gamma_0(6) = \langle - I, \: \left( \begin{smallmatrix} 1 & 1 \\ 0 & 1 \end{smallmatrix} \right), \: \left( \begin{smallmatrix} 1 & 0 \\ 6 & 1 \end{smallmatrix} \right), \: \left( \begin{smallmatrix} 5 & 2 \\ 12 & 5 \end{smallmatrix} \right) \rangle$, $\gamma_{-1/2} = W_{6, 3}$, $\gamma_0 = W_6$, $\gamma_{-1/2} = W_{6, 3}$, and $\gamma_{-1/3} = W_{6, 2}$.\\

\paragraph{\bf Location of the zeros of the Eisenstein series}
Since $W_6^{- 1} \Gamma_0(6) W_6 = W_{6, 3}^{- 1} \Gamma_0(6) W_{6, 3} = W_{6, 2}^{- 1} \Gamma_0(6) W_{6, 2} = \Gamma_0(6)$, we have
\begin{equation*}
(\sqrt{6} z)^{-k} E_{k, 6}^0(W_6 z) = (2 \sqrt{3} z + \sqrt{3})^{-k} E_{k, 6}^{-1/2}(W_{6, 3} z)
 = (3 \sqrt{2} z + \sqrt{2})^{-k} E_{k, 6}^{-1/3}(W_{6, 2} z) = E_{k, 6}^{\infty}(z).
\end{equation*}
Furthermore, we have
\begin{align*}
E_{k, 6}^0 (-1/2 + i / (2 \tan(\theta/2))) &= ((e^{i \theta} - 1) / \sqrt{6})^k E_{k, 6}^{\infty}((e^{i \theta} - 1) / 6),\\
E_{k, 6}^0 ((e^{i \theta'} - 5) / 12) &= ((5 e^{i \theta} -1) / (2 \sqrt{6}))^k E_{k, 6}^{\infty}((e^{i \theta} - 5) / 12),\\
E_{k, 6}^{-1/2} ((e^{i \theta''} - 1) / 6) &= ((e^{i \theta} + 2) / (2 \sqrt{3}))^k E_{k, 6}^{\infty}((e^{i \theta} - 1) / 6),\\
E_{k, 6}^{-1/2} (i \tan(\theta/2)) &= ((e^{i \theta} + 1) / (2 \sqrt{3}))^k E_{k, 6}^{\infty}((e^{i \theta} - 5) / 12),\\
E_{k, 6}^{-1/3} (-1/2+i \tan(\theta/2) / 6) &= ((e^{i \theta} + 1) / \sqrt{2})^k E_{k, 6}^{\infty}((e^{i \theta} - 1) / 6),\\
E_{k, 6}^{-1/3} (i / (3 \tan(\theta/2))) &= ((e^{i \theta} - 1) / (2 \sqrt{2}))^k E_{k, 6}^{\infty}((e^{i \theta} - 5) / 12),
\end{align*}
where $e^{i \theta'} = (5 - 13 \cos\theta + 12 i \sin\theta) / (13 - 5 \cos\theta)$ and $e^{i \theta''} = (- 4 - 5 \cos\theta + 3 i \sin\theta) / (5 + 4 \cos\theta)$.

For $k \leqslant 500$, we can prove that all of the zeros of $E_{k, 6}^{\infty}$ lie on the lower arcs of $\partial \mathbb{F}_6$ by numerical calculation.

\begin{figure}[hbtp]
\begin{center}
{{$E_{k, 6}^{\infty}$}\includegraphics[width=1.5in]{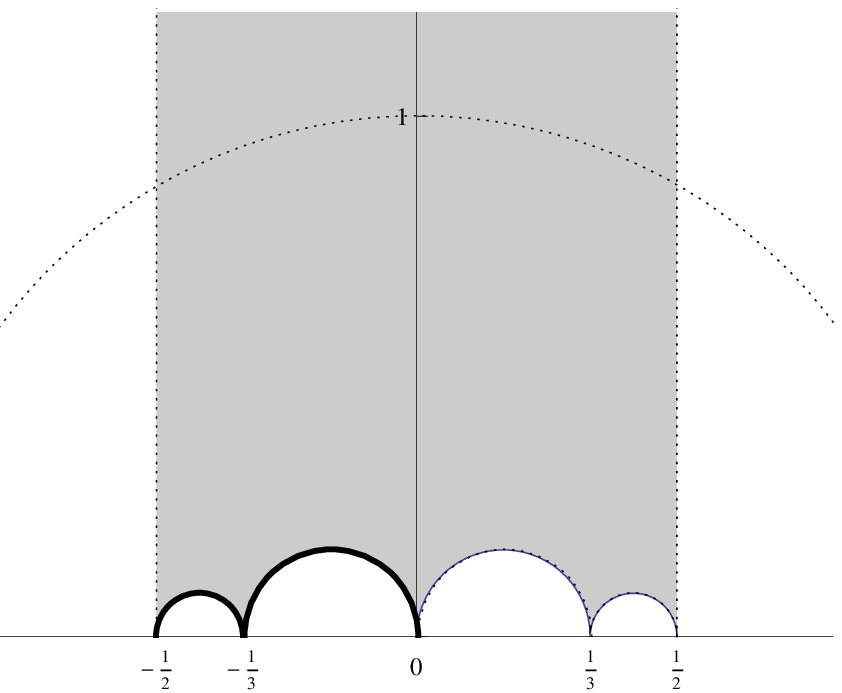}}
\quad \quad
{{$E_{k, 6}^0$}\includegraphics[width=1.5in]{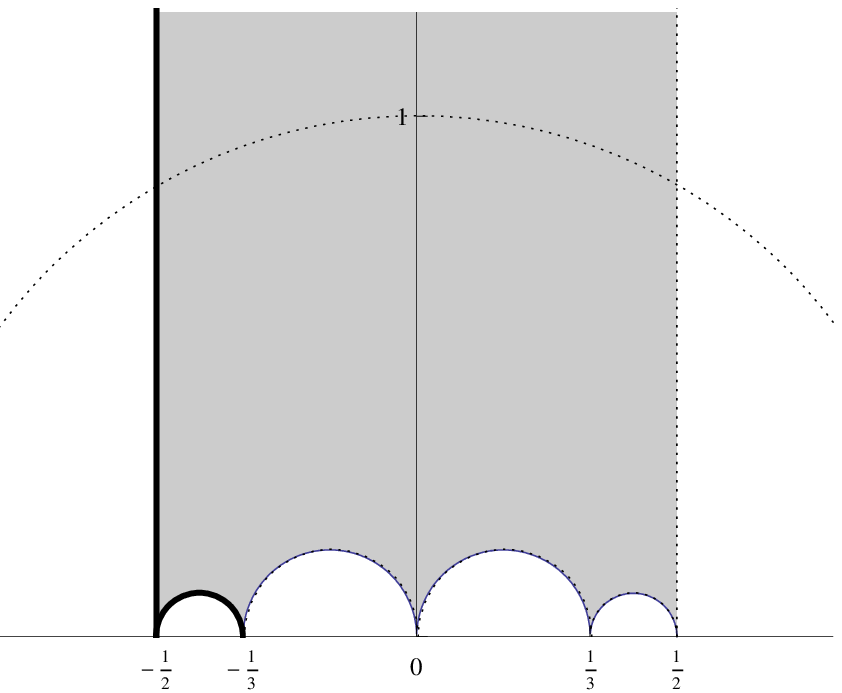}}\\

{{$E_{k, 6}^{-1/2}$}\includegraphics[width=1.5in]{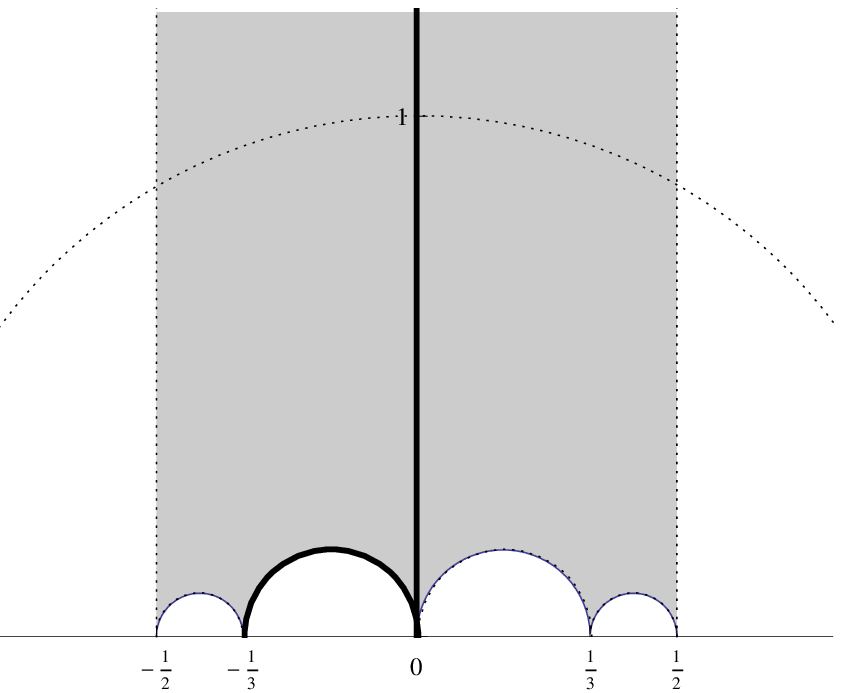}}
\quad \quad
{{$E_{k, 6}^{-1/3}$}\includegraphics[width=1.5in]{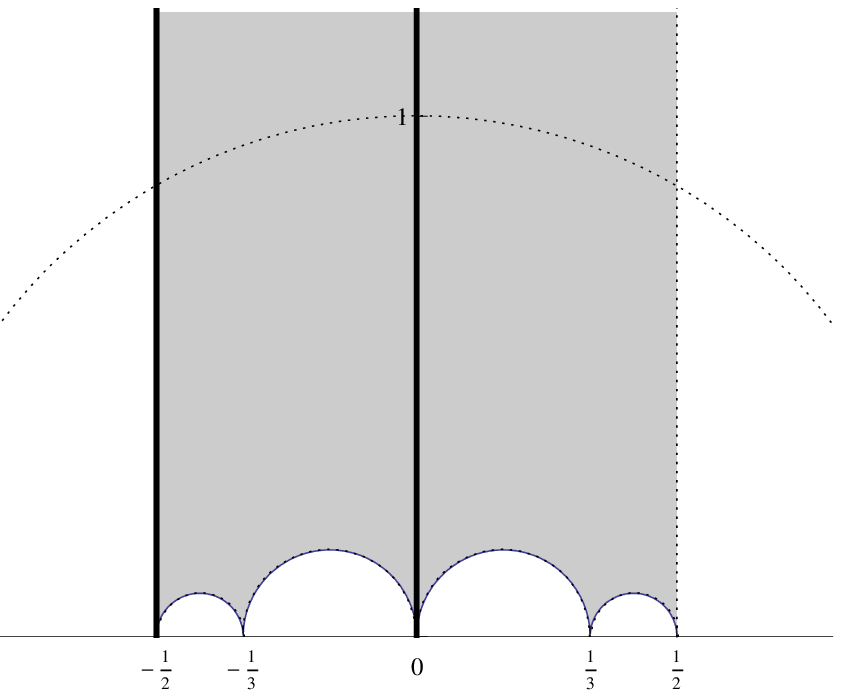}}
\end{center}
\caption{Location of the zeros of the Eisenstein series}
\end{figure}

\newpage

\paragraph{\bf Location of the zeros of Hecke type Faber Polynomial}
For $m \leqslant 200$, we can prove that all of the zeros of $F_{m, 6}$ lie on the lower arcs of $\partial \mathbb{F}_{6}$ by numerical calculation.

\clearpage

\section{Level $7$}

We have $\Gamma_0(7)+=\Gamma_0^{*}(7)$ and $\Gamma_0(7)-=\Gamma_0(7)$. We have $W_7 = \left(\begin{smallmatrix}0&-1 / \sqrt{7}\\ \sqrt{7}&0\end{smallmatrix}\right)$.\\

\subsection{$\Gamma_0^{*}(7)$}

We have $\Gamma_0^{*}(7) = \langle \left( \begin{smallmatrix} 1 & 1 \\ 0 & 1 \end{smallmatrix} \right), \: W_7, \: \left( \begin{smallmatrix} 3 & 1 \\ 7 & 2 \end{smallmatrix} \right) \rangle$.\\

\paragraph{\bf Location of the zeros of the Eisenstein series}
In \cite{SJ2}, the present author proved that all of the zeros of $E_{k, 7+}$ lie on the lower arcs of $\partial \mathbb{F}_{7+}$ if $6 \mid k$, and we prove all but at most one of the zeros of $E_{k, 7+}$ lie there if $6 \nmid k$. Furthermore, let $\alpha_7 \in [0, \pi]$ be the angle which satisfies $\tan\alpha_7 = 5 / \sqrt{3}$, and let $\alpha_{7, k} \in [0, \pi]$ be the angle which satisfies $\alpha_{7, k} \equiv k (\pi / 2 + \alpha_7) / 2 \pmod{\pi}$. We prove that all of the zeros of $E_{k, 7+}(z)$ in $\mathbb{F}^{*}(7)$ are on the lower arcs of $\partial \mathbb{F}_{7+}$ if ``$\alpha_{7, k} < (127.68/180) \pi$ or $(128.68/180) \pi < \alpha_{7, k}$ for $k \equiv 2 \pmod{6}$'' or ``$\alpha_{7, k} < (108.5/180) \pi$ or $(109.5/180) \pi < \alpha_{7, k}$ for $k \equiv 4 \pmod{6}$''.

In addition, for $k \leqslant 3000$, we can prove that all of the zeros of $E_{k, 7+}$ lie on the lower arcs of $\partial \mathbb{F}_{7+}$ by numerical calculation.

\begin{figure}[hbtp]
\begin{center}
\includegraphics[width=1.5in]{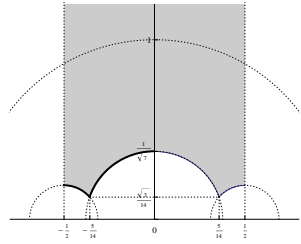}
\end{center}
\caption{Location of the zeros of the Eisenstein series}
\end{figure}

\paragraph{\bf Location of the zeros of Hecke type Faber Polynomial}
For $m = 1$ and $3 \leqslant m \leqslant 200$, we can prove that all of the zeros of $F_{m, 7+}$ lie on the lower arcs of $\partial \mathbb{F}_{7+}$ by numerical calculation. On the other hand, by numerical calculation, we can prove that all but one of the zeros of $F_{2, 7+}$ lie on the lower arcs of $\partial \mathbb{F}_{7+}$, and one of the zeros of $F_{2, 7+}$ lies on $\partial \mathbb{F}_{7+}$ but does not on the lower arcs.\\

\subsection{$\Gamma_0(7)$}

We have $\Gamma_0(7) = \langle \left( \begin{smallmatrix} 1 & 1 \\ 0 & 1 \end{smallmatrix} \right), \: - \left( \begin{smallmatrix} 1 & 0 \\ 7 & 1 \end{smallmatrix} \right), \: \left( \begin{smallmatrix} 4 & 1 \\ 7 & 2 \end{smallmatrix} \right) \rangle$ and $\gamma_0 = W_7$.\\

\paragraph{\bf Location of the zeros of the Eisenstein series}
Since $W_7^{- 1} \Gamma_0(7) W_7 = \Gamma_0(7)$, we have
\begin{equation}
E_{k, 7}^0(W_7 z) = (\sqrt{7} z)^k E_{k, 7}^{\infty}(z).
\end{equation}
Furthermore, we have
\begin{align*}
E_{k, 7}^0 (- 1/2 + i / (2 \tan\theta/2)) &= ((e^{i \theta} - 1) / \sqrt{7})^k E_{k, 7}^{\infty}((e^{i \theta} - 1) / 7),\\
E_{k, 7}^0 ((e^{i \theta'} + 2) / 3) &= ((- 2 e^{i \theta} - 1) / \sqrt{7})^k E_{k, 7}^{\infty}((e^{i \theta} - 3) / 7),\\
E_{k, 7}^0 ((e^{i (\pi - \theta')} - 2) / 3) &= ((e^{i \theta} + 2) / \sqrt{7})^k E_{k, 7}^{\infty}((e^{i \theta} + 2) / 7),
\end{align*}
where $e^{i \theta'} = (- 4 - 5 \cos\theta + 3 i \sin\theta) / (5 + 4 \cos\theta)$.
\begin{figure}[hbtp]
\begin{center}
{{$E_{k, 7}^{\infty}$}\includegraphics[width=1.5in]{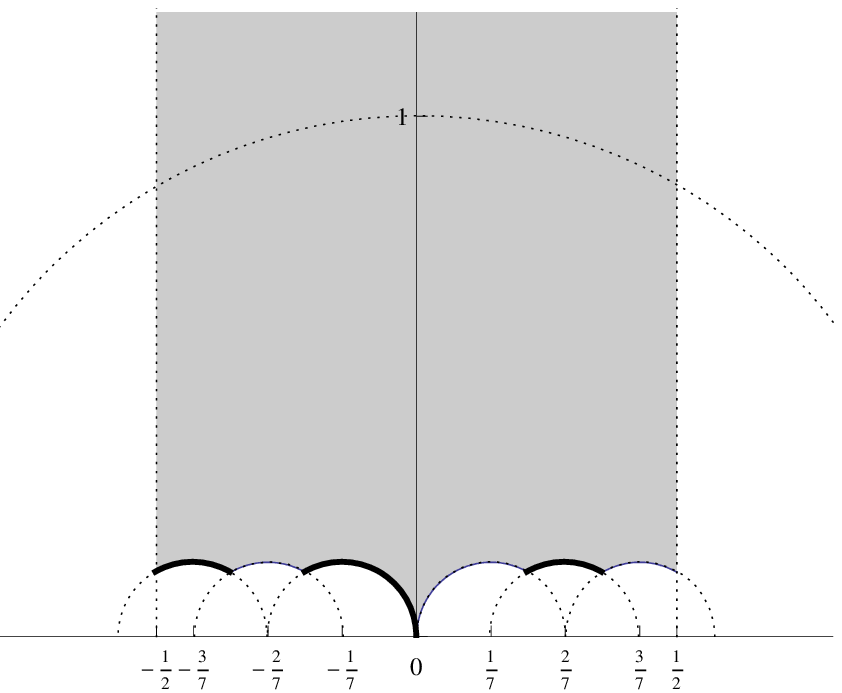}}
\quad \quad
{{$E_{k, 7}^{0}$}\includegraphics[width=1.5in]{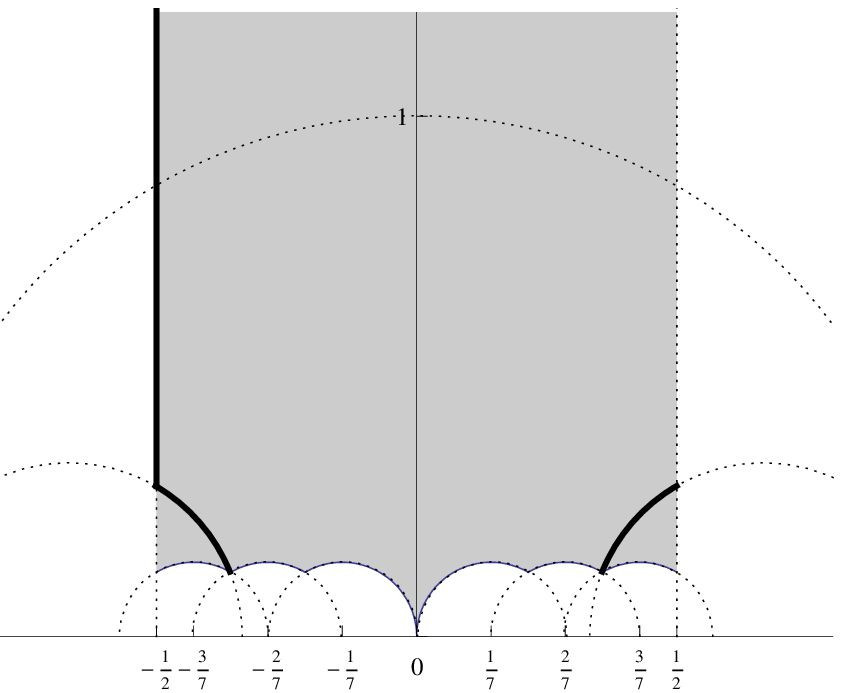}}
\end{center}
\caption{Neighborhood of location of the zeros of the Eisenstein series}
\end{figure}

\begin{figure}[hbtp]
\begin{center}
{{\small Lower arcs of $\partial \mathbb{F}_7$}\includegraphics[width=2.5in]{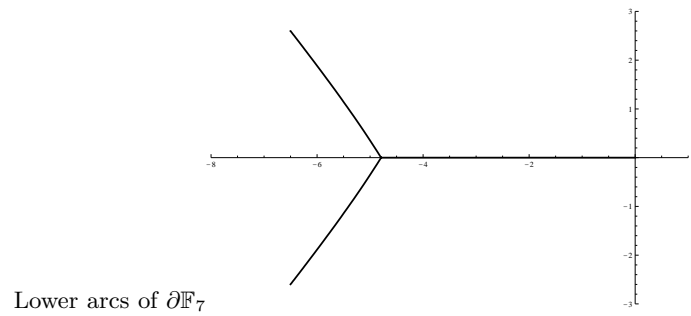}}
\end{center}
\caption{Image by $J_7$}\label{Im-J7B}
\end{figure}

\begin{figure}[hbtp]
\begin{center}
\includegraphics[width=6.3in]{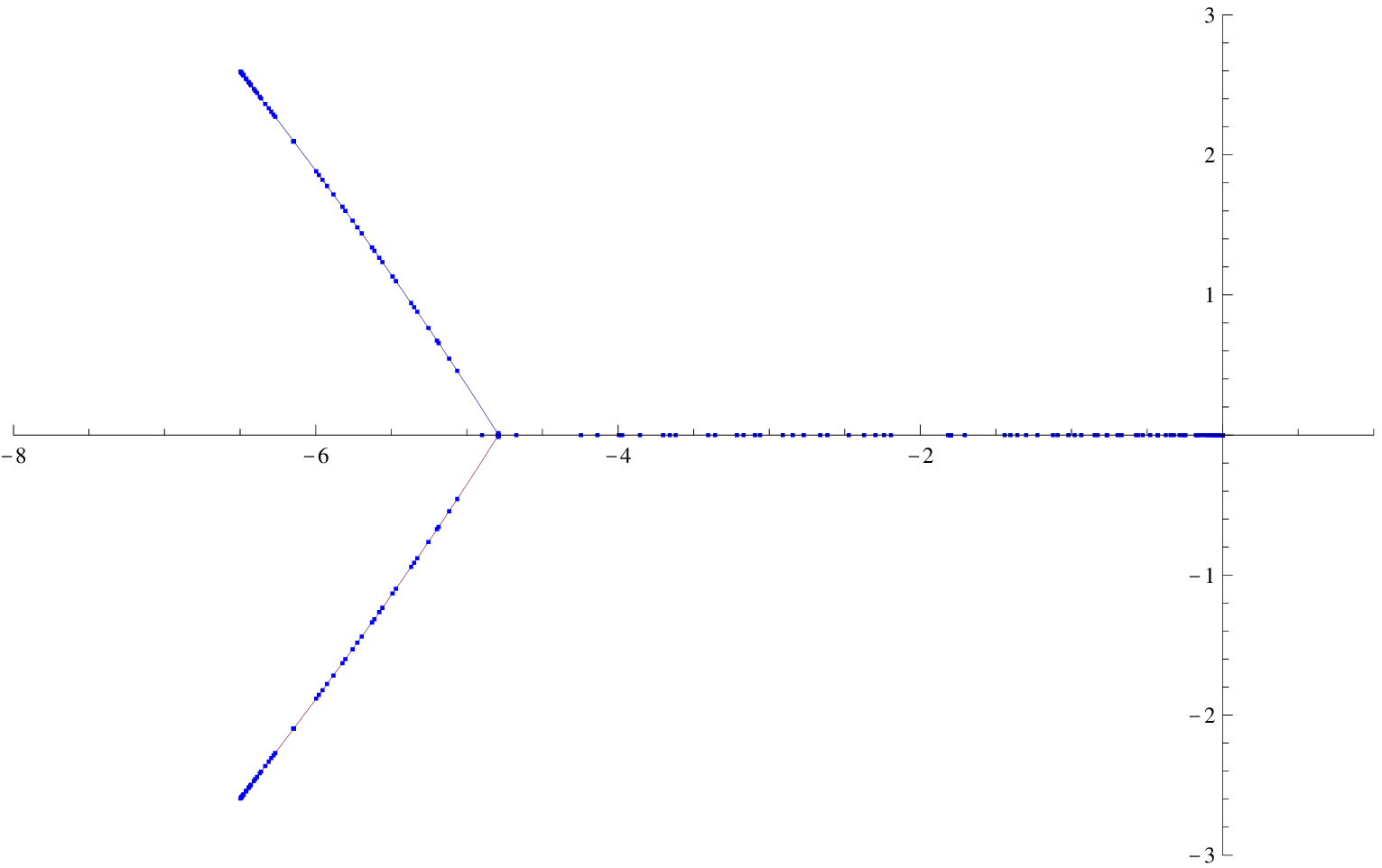}\\
The zeros of $E_{k, 7}^{\infty}$ for $4 \leqslant k \leqslant 40$\\
\includegraphics[width=6.3in]{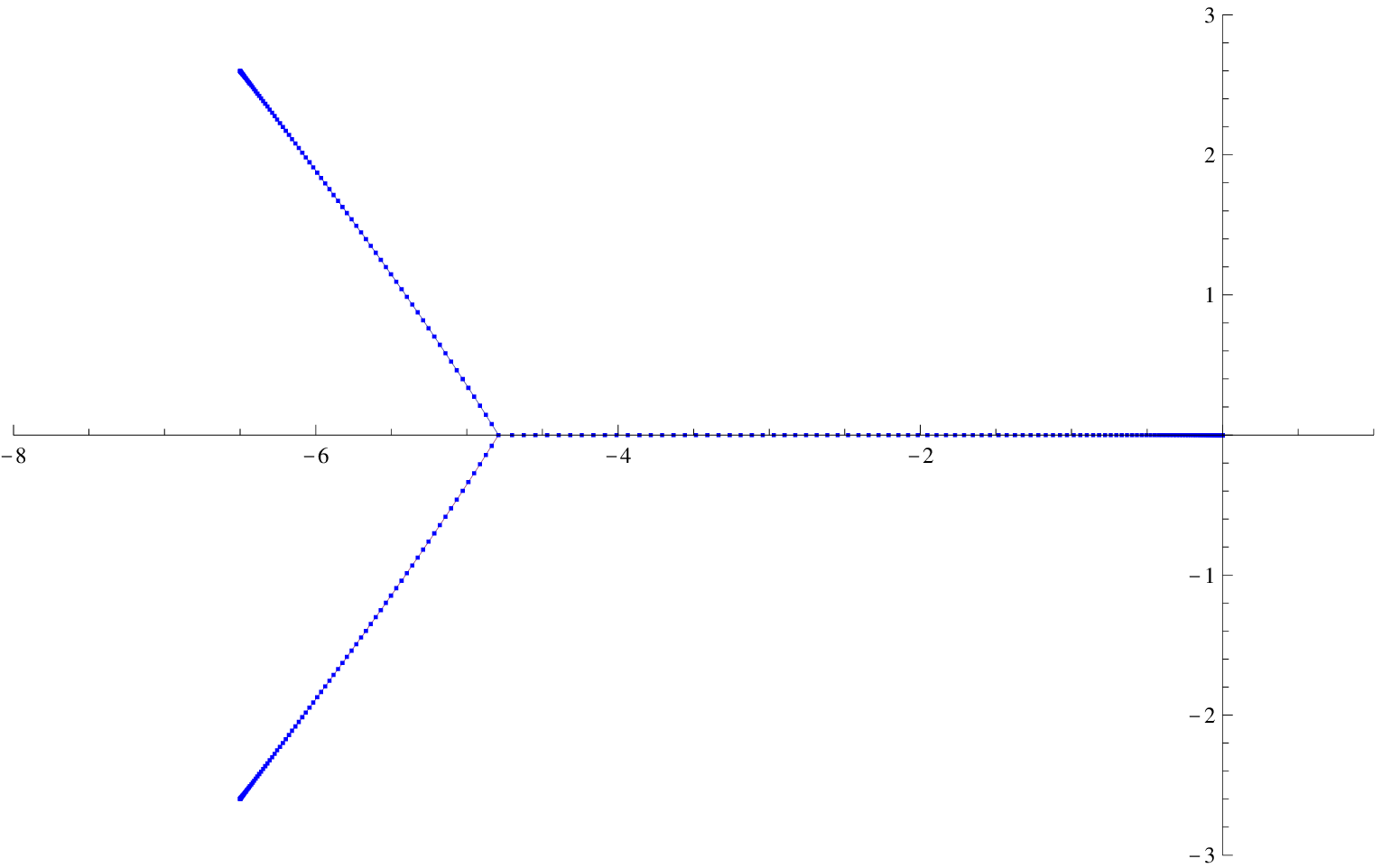}\\
The zeros of $E_{1000, 7}^{\infty}$
\end{center}
\caption{Image by $J_7$}\label{Im-J7Bz}
\end{figure}

Now, we can observe that some zeros of $E_{k, 7}^{\infty}$ do not lie on the lower arcs of $\partial \mathbb{F}_7$ for small weight $k$ by numerical calculation. However, when the weight $k$ increases, then the location of the zeros seems to approach to lower arcs of $\partial \mathbb{F}_7$. (see Figure \ref{Im-J7Bz})\\

\begin{figure}[hbtp]
\begin{center}
\includegraphics[width=6.3in]{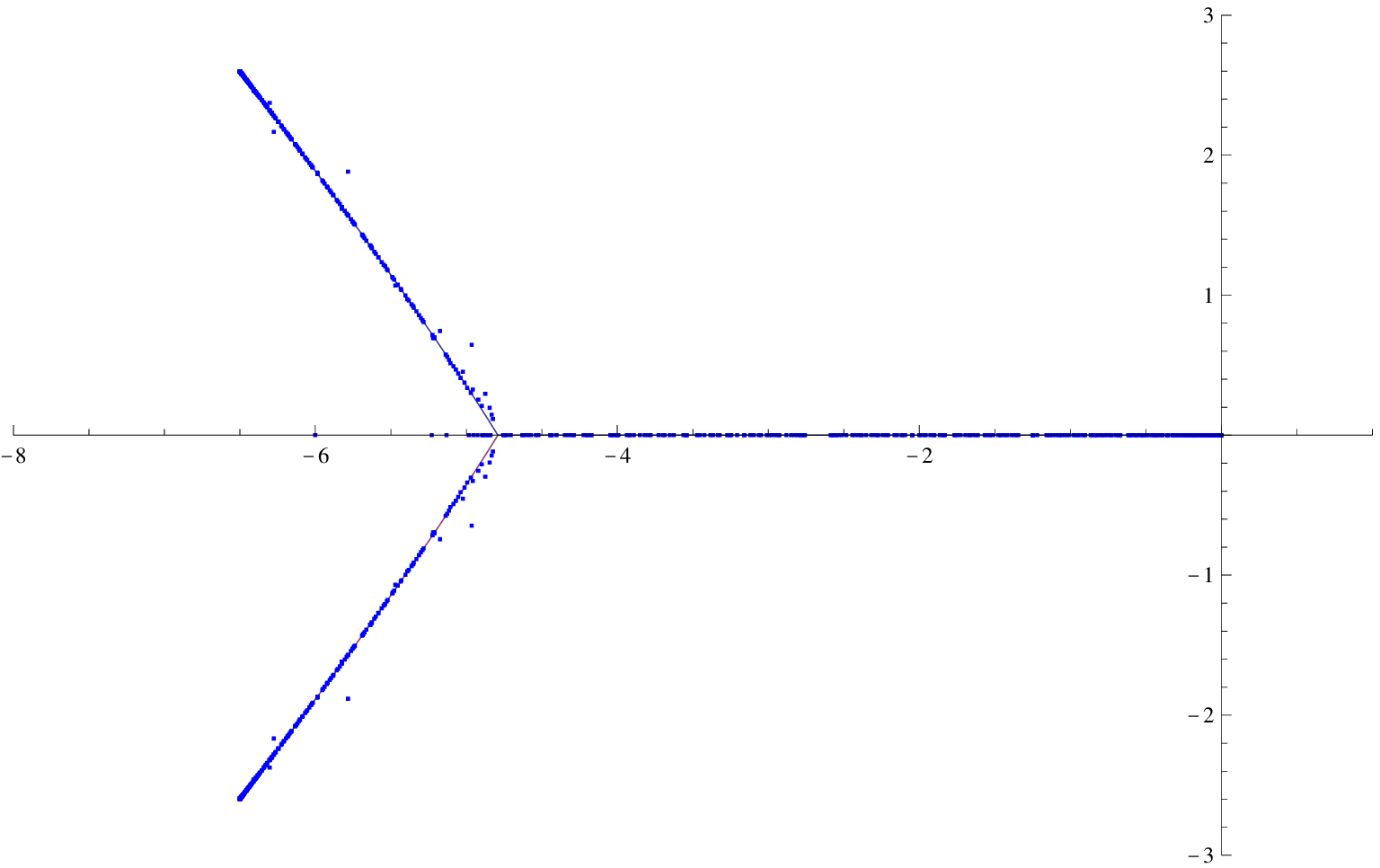}\\
The zeros of $F_{m, 7}$ for $m \leqslant 40$\\
\includegraphics[width=6.3in]{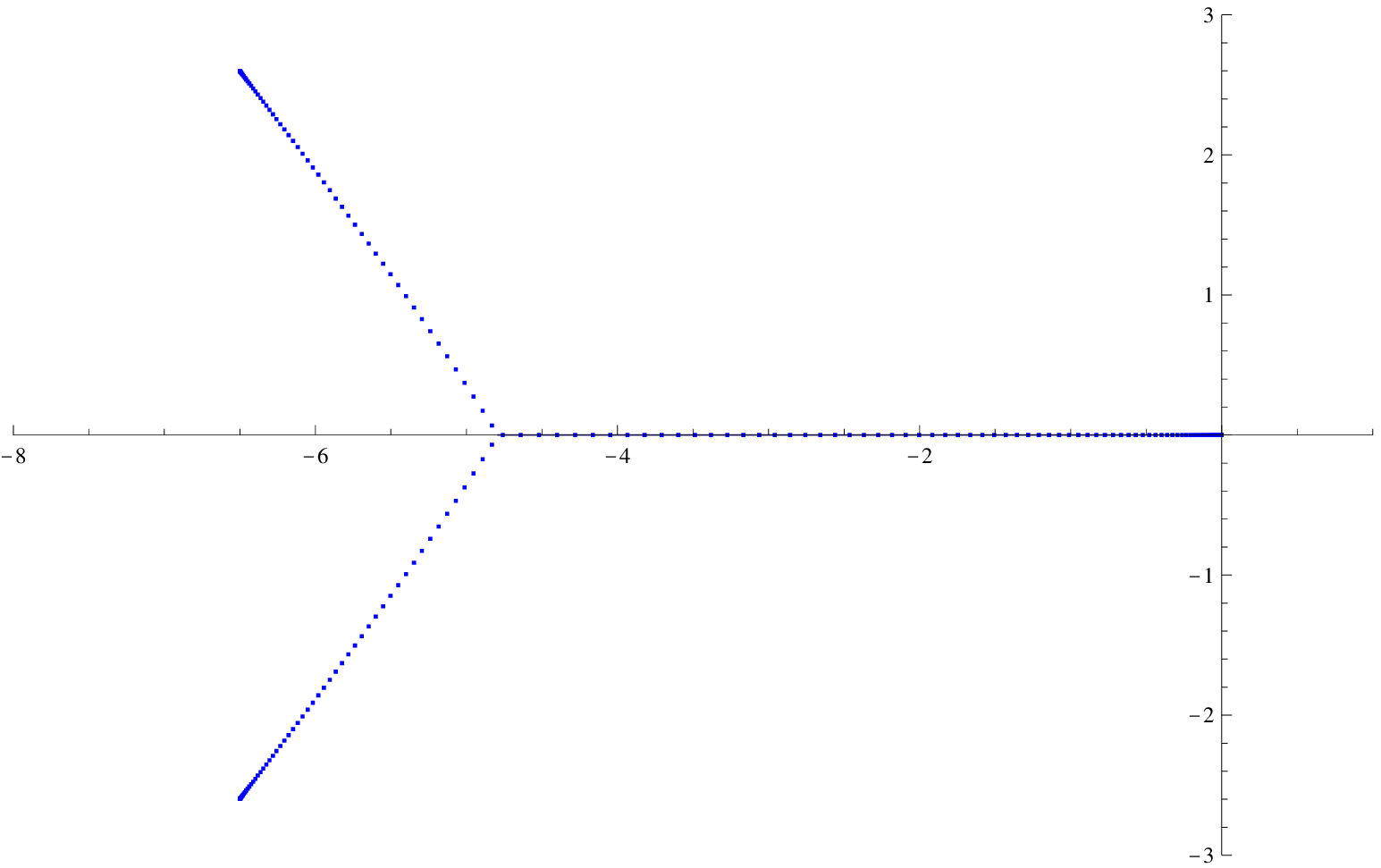}\\
The zeros of $F_{200, 7}$
\end{center}
\caption{Image by $J_7$}\label{Im-J7Bhz}
\end{figure}

\paragraph{\bf Location of the zeros of Hecke type Faber Polynomial}
Similarly to the Eisenstein series, we can observe that some zeros of $F_{m, 7}$ do not lie on the lower arcs of $\partial \mathbb{F}_7$ for small weight $m$ by numerical calculation. However, when the weight $m$ increases, then the location of the zeros seems to approach to lower arcs of $\partial \mathbb{F}_7$. (see Figure \ref{Im-J7Bhz})

\clearpage

\section{Level $8$}

We have $\Gamma_0(8)+ = \Gamma_0(8)+8 = \Gamma_0^{*}(8)$ and $\Gamma_0(8)-=\Gamma_0(8)$. We have $W_8 = \left(\begin{smallmatrix} 0 & -1 / (2 \sqrt{2}) \\ 2 \sqrt{2} & 0 \end{smallmatrix}\right)$, $W_{8-, 2} := \left(\begin{smallmatrix} -1 & -1/2 \\ 4 & 1 \end{smallmatrix}\right)$, and $W_{8-, 4} := \left(\begin{smallmatrix} -\sqrt{2} & -3/(2 \sqrt{2}) \\ 2 \sqrt{2} & \sqrt{2} \end{smallmatrix}\right)$.\\

\subsection{$\Gamma_0(8)+8 = \Gamma_0^{*}(8)$}

We have $\Gamma_0^{*}(8) = \langle \left( \begin{smallmatrix} 1 & 1 \\ 0 & 1 \end{smallmatrix} \right), \: W_8, \: \left( \begin{smallmatrix} 3 & 1 \\ 8 & 3 \end{smallmatrix} \right) \rangle$ and $\gamma_{-1/2} = W_{8-, 4}$.\\

\paragraph{\bf Location of the zeros of the Eisenstein series}
Since $W_{8-, 4}^{- 1} \Gamma_0^{*}(8) W_{8-, 4} = \Gamma_0^{*}(8)$, we have
\begin{equation}
E_{k, 8+8}^{-1/2}(W_{8-, 4} z) = (2 \sqrt{2} z + \sqrt{2})^k E_{k, 8+8}^{\infty}(z).
\end{equation}
Furthermore, we have
\begin{align*}
E_{k, 8+8}^{-1/2} (i \tan(\theta/2) / 2) &= ((e^{i \theta} + 1) / (2 \sqrt{2}))^k E_{k, 8+8}^{\infty}((e^{i \theta} - 5) / 12),\\
E_{k, 8+8}^{-1/2} (e^{i \theta'} / (2 \sqrt{2})) &= (\sqrt{2} e^{i \theta} + 1)^k E_{k, 8+8}^{\infty}(e^{i \theta} / (2 \sqrt{2})),
\end{align*}
where $e^{i \theta'} = (- 2 \sqrt{2} - 3 \cos\theta + i \sin\theta) / (3 + 2 \sqrt{2} \cos\theta)$.

For $k \leqslant 600$, we can prove that all of the zeros of $E_{k, 8+8}^{\infty}$ lie on the lower arcs of $\partial \mathbb{F}_{8+8}$ by numerical calculation.

\begin{figure}[htbp]
\begin{center}
{{$E_{k, 8+8}^{\infty}$}\includegraphics[width=1.5in]{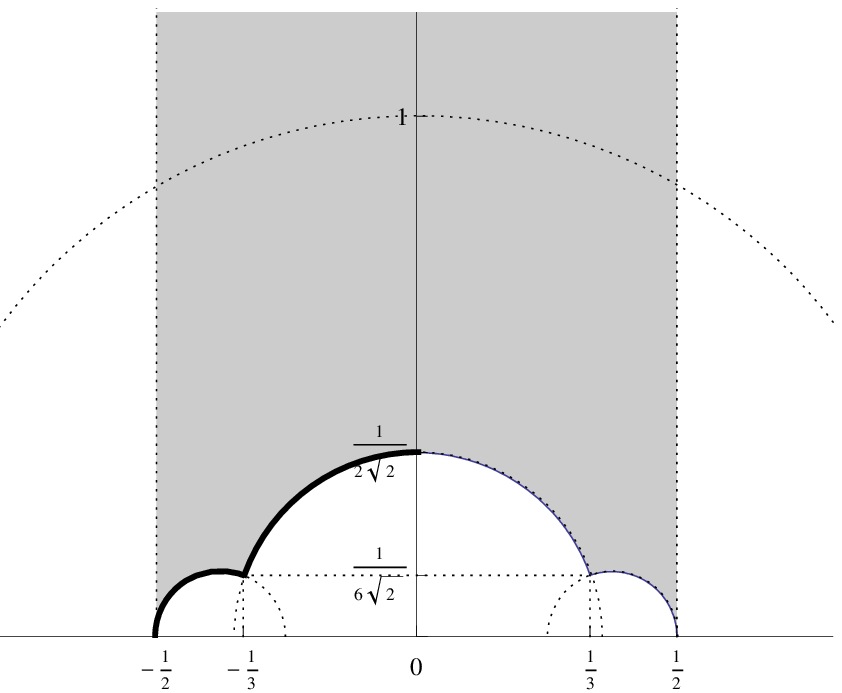}}
\quad \quad
{{$E_{k, 8+8}^{-1/2}$}\includegraphics[width=1.5in]{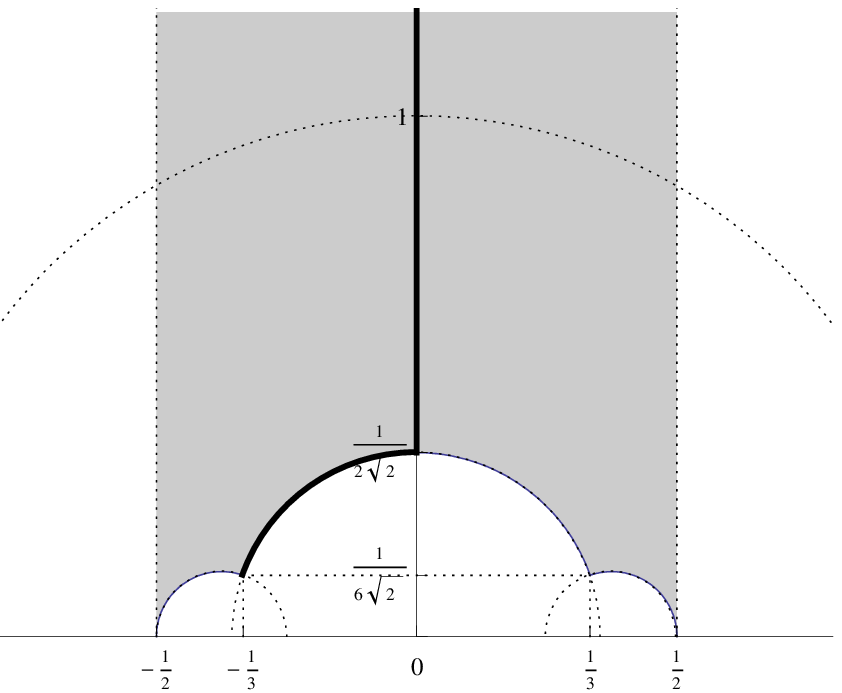}}
\end{center}
\caption{Location of the zeros of the Eisenstein series}
\end{figure}

\paragraph{\bf Location of the zeros of Hecke type Faber Polynomial}
For every integer $m \leqslant 200$ such that $m \not\equiv 0 \pmod{4}$, we can prove that all of the zeros of $F_{m, 8+8}$ lie on the lower arcs of $\partial \mathbb{F}_{8+8}$ by numerical calculation. On the other hand, by numerical calculation, for every integer $m \leqslant 200$ such that $m \equiv 0 \pmod{4}$, we can prove that all but one of the zeros of $F_{m, 8+8}$ lie on the lower arcs of $\partial \mathbb{F}_{8+8}$, and one of the zeros of $F_{m, 8+8}$ lies on $\partial \mathbb{F}_{8+8}$ but does not on the lower arcs.\\

\subsection{$\Gamma_0(8)$}

We have $\Gamma_0(8) = \langle - I, \: \left( \begin{smallmatrix} 1 & 1 \\ 0 & 1 \end{smallmatrix} \right), \: \left( \begin{smallmatrix} 1 & 0 \\ 8 & 1 \end{smallmatrix} \right), \: \left( \begin{smallmatrix} 3 & 1 \\ 8 & 3 \end{smallmatrix} \right) \rangle$, $\gamma_0 = W_8$, $\gamma_{-1/2} = W_{8-, 4}$, and $\gamma_{-1/4} = W_{8-, 2}$.\\

\paragraph{\bf Location of the zeros of the Eisenstein series}
Since $W_8^{- 1} \Gamma_0(8) W_8 = W_{8-, 4}^{- 1} \Gamma_0(8) W_{8-, 4} = W_{8-, 2}^{- 1} \Gamma_0(8) W_{8-, 2} = \Gamma_0(8)$, we have
\begin{equation*}
(\sqrt{8} z)^{-k} E_{k, 8}^0(W_8 z) = (2 \sqrt{2} z + \sqrt{2})^{-k} E_{k, 8}^{-1/2}(W_{8, 3} z)
 = (4 z + 1)^{-k} E_{k, 8}^{-1/4}(W_{8-, 2} z) = E_{k, 8}^{\infty}(z).
\end{equation*}
Furthermore, we have
\begin{align*}
E_{k, 8}^0 (-1/2 + i / (2 \tan(\theta/2))) &= ((e^{i \theta} - 1) / (2 \sqrt{2}))^k E_{k, 8}^{\infty}((e^{i \theta} - 1) / 8),\\
E_{k, 8}^0 ((e^{i \theta'} - 3) / 8) &= ((3 e^{i \theta} -1) / (2 \sqrt{2}))^k E_{k, 8}^{\infty}((e^{i \theta} - 3) / 8),\\
E_{k, 8}^{-1/2} ((e^{i \theta''} - 1) / 8) &= (- (3 e^{i \theta} + 1) / (2 \sqrt{2}))^k E_{k, 8}^{\infty}((e^{i \theta} - 1) / 8),\\
E_{k, 8}^{-1/2} (i \tan(\theta/2) / 2) &= ((e^{i \theta} + 1) / (2 \sqrt{2}))^k E_{k, 8}^{\infty}((e^{i \theta} - 3) / 8),\\
E_{k, 8}^{-1/4} (-1/2 + i \tan(\theta/2) / 2) &= ((e^{i \theta} + 1) / 2)^k E_{k, 8}^{\infty}((e^{i \theta} - 1) / 8),\\
E_{k, 8}^{-1/4} (i / (2 \tan(\theta/2))) &= ((e^{i \theta} - 1) / 4)^k E_{k, 8}^{\infty}((e^{i \theta} - 3) / 8),
\end{align*}
where $e^{i \theta'} = (3 - 5 \cos\theta + 4 i \sin\theta) / (5 - 3 \cos\theta)$ and $e^{i \theta''} = (-3 - 5 \cos\theta + 4 i \sin\theta) / (5 + 3 \cos\theta)$.

Now, recall that $E_{k, 8}^{\infty}(z) = E_{k, 2}^{\infty}(4 z)$. Similarly to $\Gamma_0(4)$, for $k \leqslant 1000$, we can prove that all of the zeros of $E_{k, 2}^{\infty}$ lie on the lower arcs of $\partial \mathbb{F}_2$ by numerical calculation, then we have all of the zeros of $E_{k, 8}^{\infty}$ in the lower arcs of $\partial \mathbb{F}_8$.

\begin{figure}[hbtp]
\begin{center}
{{$E_{k, 8}^{\infty}$}\includegraphics[width=1.5in]{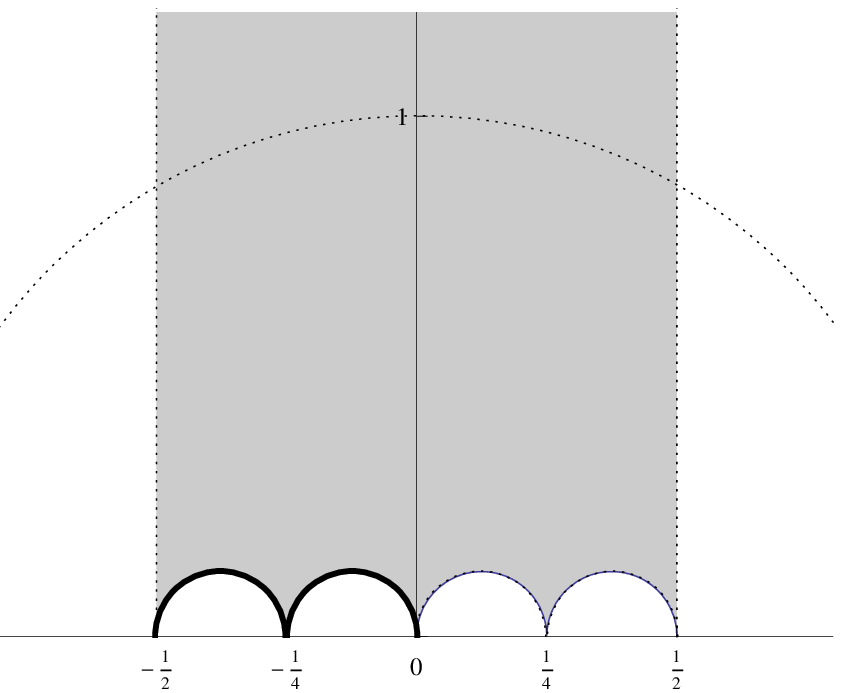}}
\quad \quad
{{$E_{k, 8}^0$}\includegraphics[width=1.5in]{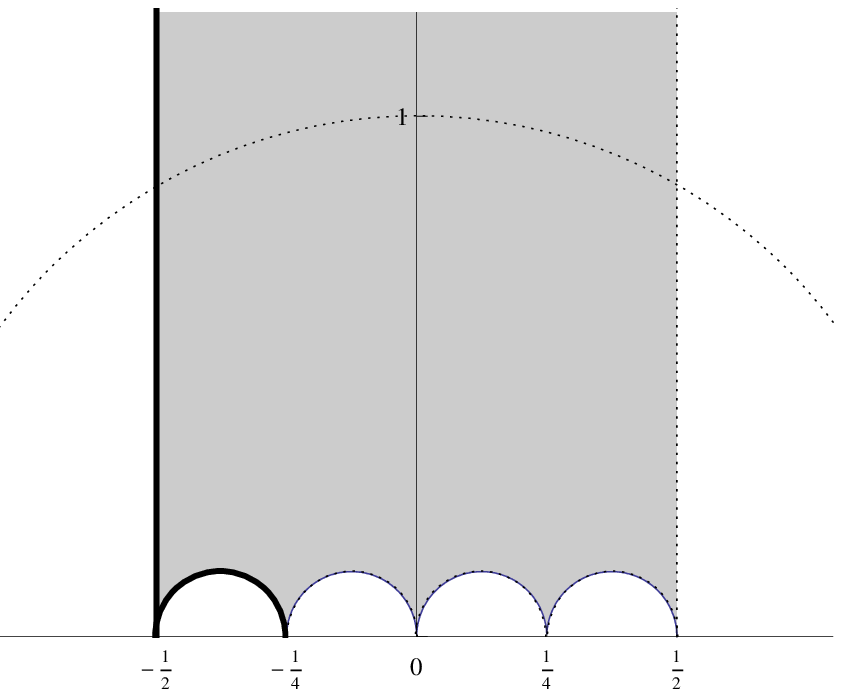}}\\

{{$E_{k, 8}^{-1/2}$}\includegraphics[width=1.5in]{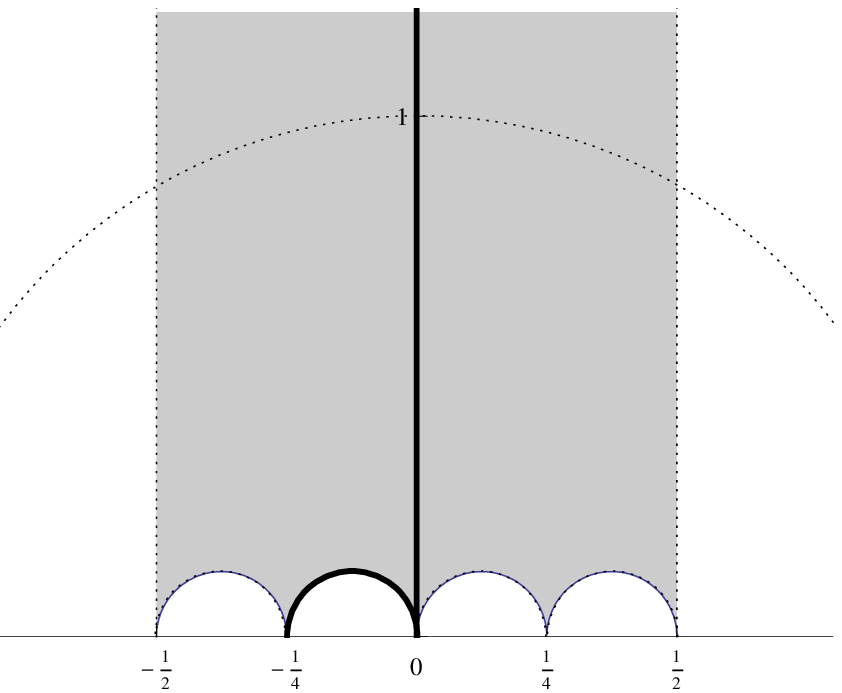}}
\quad \quad
{{$E_{k, 8}^{-1/4}$}\includegraphics[width=1.5in]{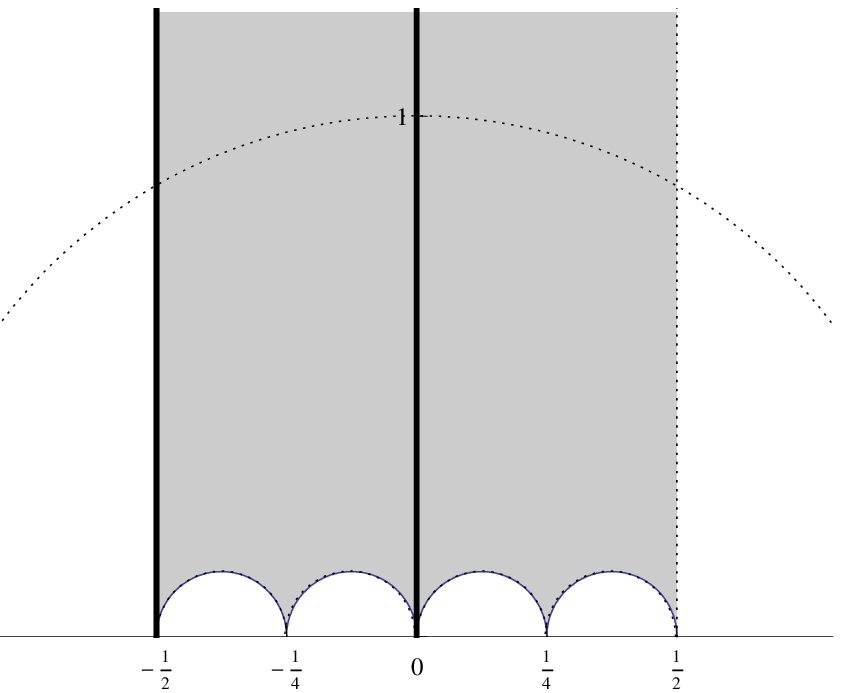}}
\end{center}
\caption{Location of the zeros of the Eisenstein series}
\end{figure}

\paragraph{\bf Location of the zeros of Hecke type Faber Polynomial}
For $m \leqslant 200$, we can prove that all of the zeros of $F_{m, 8}$ lie on the lower arcs of $\partial \mathbb{F}_8$ by numerical calculation.\\ \clearpage

\section{Level $9$}

We have $\Gamma_0(9)+ = \Gamma_0(9)+9 = \Gamma_0^{*}(9)$ and $\Gamma_0(9)-=\Gamma_0(9)$. We have $W_9 = \left(\begin{smallmatrix} 0 & -1/3 \\ 3 & 0 \end{smallmatrix}\right)$, $W_{9-, 3} := \left(\begin{smallmatrix} -1 & -2/3 \\ 3 & 1 \end{smallmatrix}\right)$, and $W_{9-, -3} := \left(\begin{smallmatrix} 1 & -2/3 \\ 3 & -1 \end{smallmatrix}\right)$.\\

\subsection{$\Gamma_0(9)+9 = \Gamma_0^{*}(9)$}

We define $\Gamma_0^{*}(9) = \langle \left( \begin{smallmatrix} 1 & 1 \\ 0 & 1 \end{smallmatrix} \right), \: W_9, \: \left( \begin{smallmatrix} 5 & 1 \\ 9 & 2 \end{smallmatrix} \right) \rangle$ and $\gamma_{-1/3} = W_{9-, 3}$.

\paragraph{\bf Location of the zeros of the Eisenstein series}
Since $W_{9-, 3}^{- 1} \Gamma_0^{*}(9) W_{9-, 3} = \Gamma_0^{*}(9)$, we have
\begin{equation}
E_{k, 9+9}^{-1/3}(W_{9-, 3} z) = (3 z)^k E_{k, 9+9}^{\infty}(z).
\end{equation}
Furthermore, we have
\begin{align*}
E_{k, 9+9}^{-1/3} (-1/2 + i / (6 \tan(\theta/2))) &= (e^{i \theta})^k E_{k, 9+9}^{\infty}(e^{i \theta} / 3),\\
E_{k, 9+9}^{-1/3} (i / (3 \tan(\theta/2))) &= ((e^{i \theta} - 3)/2)^k E_{k, 9+9}^{\infty}((e^{i \theta} - 3)/6).
\end{align*}

For $k \leqslant 600$, we can prove that all of the zeros of $E_{k, 9+9}^{\infty}$ lie on the lower arcs of $\partial \mathbb{F}_{9+9}$ by numerical calculation.

\begin{figure}[htbp]
\begin{center}
{{$E_{k, 9+9}^{\infty}$}\includegraphics[width=1.5in]{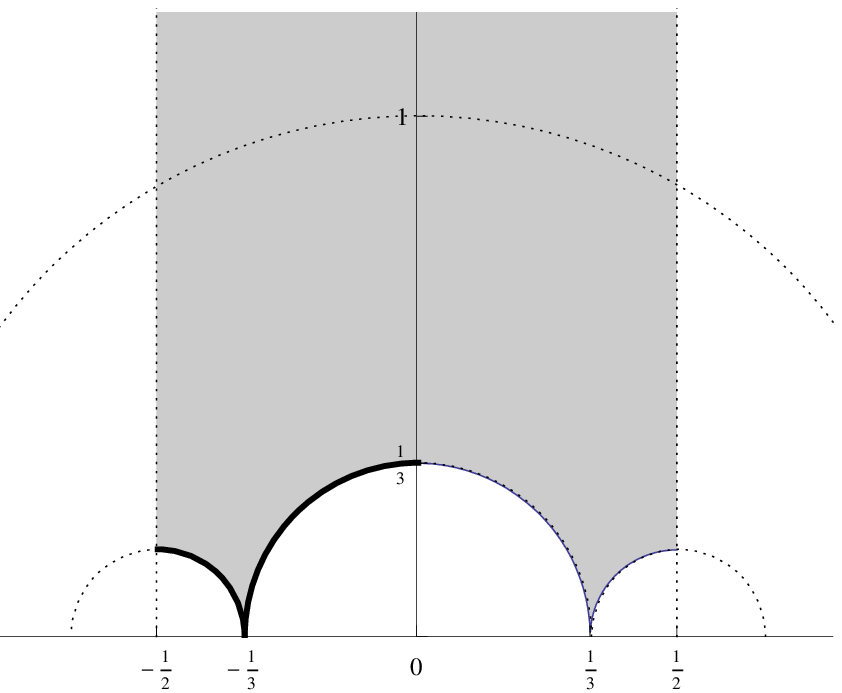}}
\quad \quad
{{$E_{k, 9+9}^{-1/3}$}\includegraphics[width=1.5in]{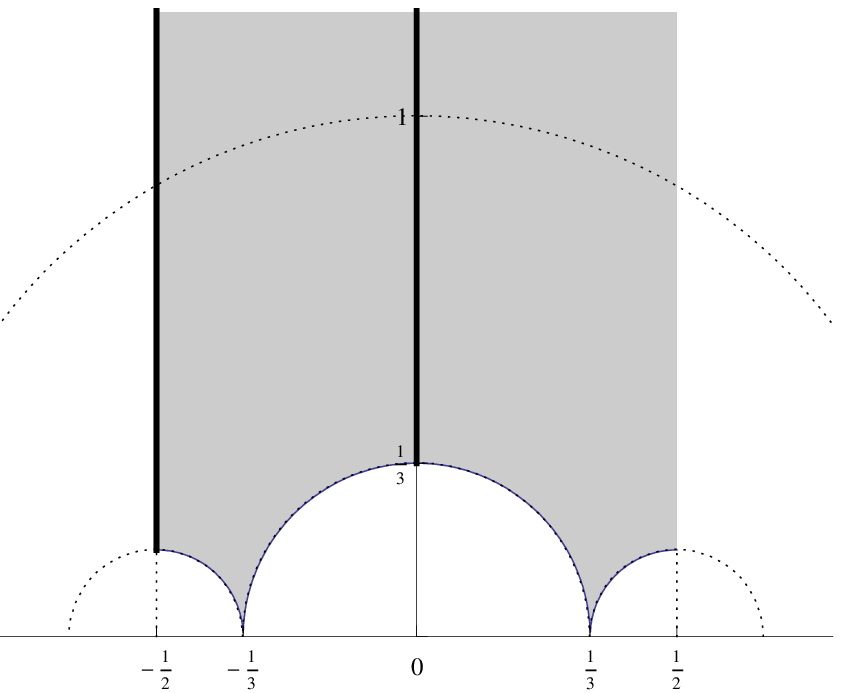}}
\end{center}
\caption{Location of the zeros of the Eisenstein series}
\end{figure}

\paragraph{\bf Location of the zeros of Hecke type Faber Polynomial}
For $m \leqslant 200$, we can prove that all of the zeros of $F_{m, 9+9}$ lie on the lower arcs of $\partial \mathbb{F}_{9+9}$ by numerical calculation.\\

\subsection{$\Gamma_0(9)$}

We have $\Gamma_0(9) = \langle - I, \: \left( \begin{smallmatrix} 1 & 1 \\ 0 & 1 \end{smallmatrix} \right), \: \left( \begin{smallmatrix} 1 & 0 \\ 9 & 1 \end{smallmatrix} \right), \: \left( \begin{smallmatrix} 5 & 1 \\ 9 & 2 \end{smallmatrix} \right) \rangle$, $\gamma_0 = W_9$, $\gamma_{-1/3} = W_{9-, 3}$, and $\gamma_{1/3} = W_{9-, -3}$.\\

\paragraph{\bf Location of the zeros of the Eisenstein series}
Since $W_9^{- 1} \Gamma_0(9) W_9 = W_{9-, 3}^{- 1} \Gamma_0(9) W_{9-, 3} = W_{9-, -3}^{- 1} \Gamma_0(9) W_{9-, -3} = \Gamma_0(9)$, we have
\begin{equation*}
(\sqrt{9} z)^{-k} E_{k, 9}^0(W_9 z) = (3 z + 1)^{-k} E_{k, 9}^{-1/3}(W_{9-, 3} z)
 = (3 z - 1)^{-k} E_{k, 9}^{1/3}(W_{9-, -3} z) = E_{k, 9}^{\infty}(z).
\end{equation*}
Furthermore, we have
\begin{allowdisplaybreaks}
\begin{align*}
E_{k, 9}^0 (-1/2 + i / (2 \tan(\theta/2))) &= ((e^{i \theta} - 1) / 3)^k E_{k, 9}^{\infty}((e^{i \theta} - 1) / 9),\\
E_{k, 9}^0 ((e^{i \theta'} + 2) / 3) &= ((- 2 e^{i \theta} + 1) / 3)^k E_{k, 9}^{\infty}((e^{i \theta} - 4) / 9),\\
E_{k, 9}^0 ((e^{i (\pi - \theta')} - 2) / 3) &= ((e^{i \theta} + 2) / 3)^k E_{k, 9}^{\infty}((e^{i \theta} + 2) / 9),\\
E_{k, 9}^{-1/3} (e^{i (\pi - \theta')} / 3) &= ((e^{i \theta} + 2) / 3)^k E_{k, 9}^{\infty}((e^{i \theta} - 1) / 9),\\
E_{k, 9}^{-1/3} (1/6 + i / (2 \tan(\theta/2))) &= ((e^{i \theta} - 1) / 3)^k E_{k, 9}^{\infty}((e^{i \theta} - 4) / 9),\\
E_{k, 9}^{-1/3} ((e^{i \theta'} + 1) / 3) &= ((2 e^{i \theta} + 1) / 3)^k E_{k, 9}^{\infty}((e^{i \theta} + 2) / 9),\\
E_{k, 9}^{1/3} (e^{i \theta'} / 3) &= ((- 2 e^{i \theta} - 1) / 3)^k E_{k, 9}^{\infty}((e^{i \theta} - 1) / 9),\\
E_{k, 9}^{1/3} ((e^{i (\pi - \theta')} - 1) / 3) &= ((- e^{i \theta} - 2) / 3)^k E_{k, 9}^{\infty}((e^{i \theta} - 4) / 9),\\
E_{k, 9}^{1/3} (-1/6 + i / (2 \tan(\theta/2))) &= ((e^{i \theta} - 1) / 3)^k E_{k, 9}^{\infty}((e^{i \theta} + 2) / 9),
\end{align*}
\end{allowdisplaybreaks}
where $e^{i \theta'} = (- 4 - 5 \cos\theta + 3 i \sin\theta) / (5 + 4 \cos\theta)$.

Now, recall that $E_{k, 9}^{\infty}(z) = E_{k, 3}^{\infty}(3 z)$. Moreover, by the transformation with $\left( \begin{smallmatrix} 1 & \pm 1 \\ 0 & 1 \end{smallmatrix} \right)$ for $E_{k, 3}^{\infty}$, we have
\begin{equation*}
E_{k, 9}^{\infty}((e^{i \theta} - 1) / 9) = E_{k, 3}^{\infty}((e^{i \theta} - 1) / 3)
 = E_{k, 2}^{\infty}((e^{i (\pi - \theta)} - 1 \pm 3) / 3) = E_{k, 9}^{\infty}((e^{i (\pi - \theta)} - 1 \pm 3) / 9).
\end{equation*}
\begin{figure}[hbtp]
\begin{center}
{{$E_{k, 9}^{\infty}$}\includegraphics[width=1.5in]{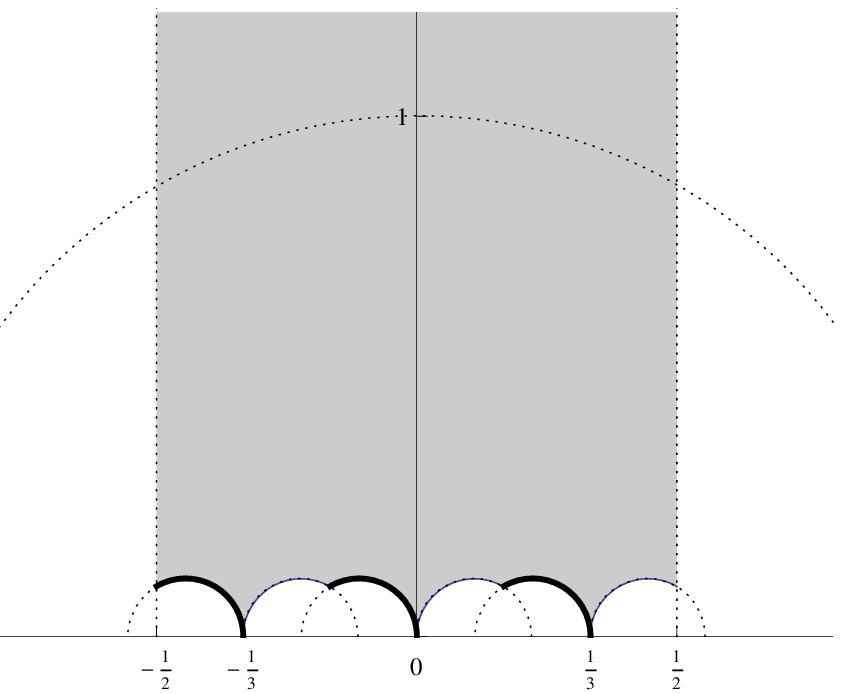}}
\quad \quad
{{$E_{k, 9}^0$}\includegraphics[width=1.5in]{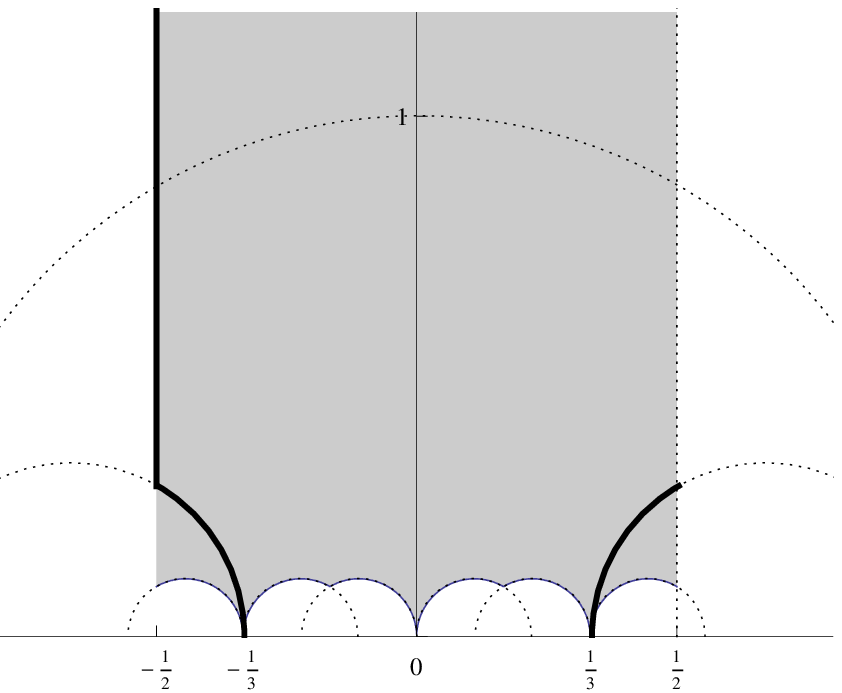}}\\

{{$E_{k, 9}^{-1/3}$}\includegraphics[width=1.5in]{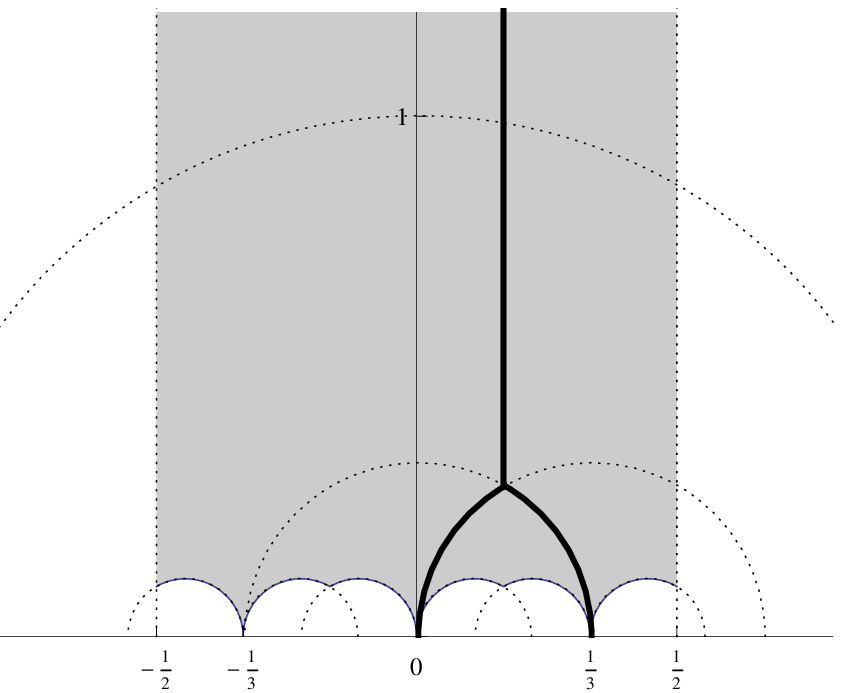}}
\quad \quad
{{$E_{k, 9}^{1/3}$}\includegraphics[width=1.5in]{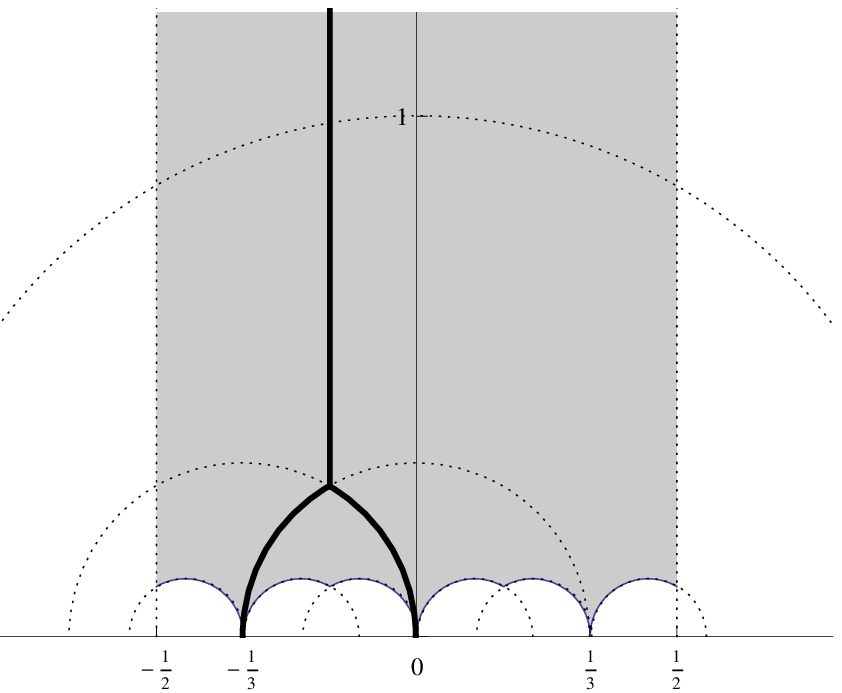}}
\end{center}
\caption{Location of the zeros of the Eisenstein series}
\end{figure}

\begin{figure}[hbtp]
\begin{center}
{{\small Lower arcs of $\partial \mathbb{F}_9$}\includegraphics[width=2.5in]{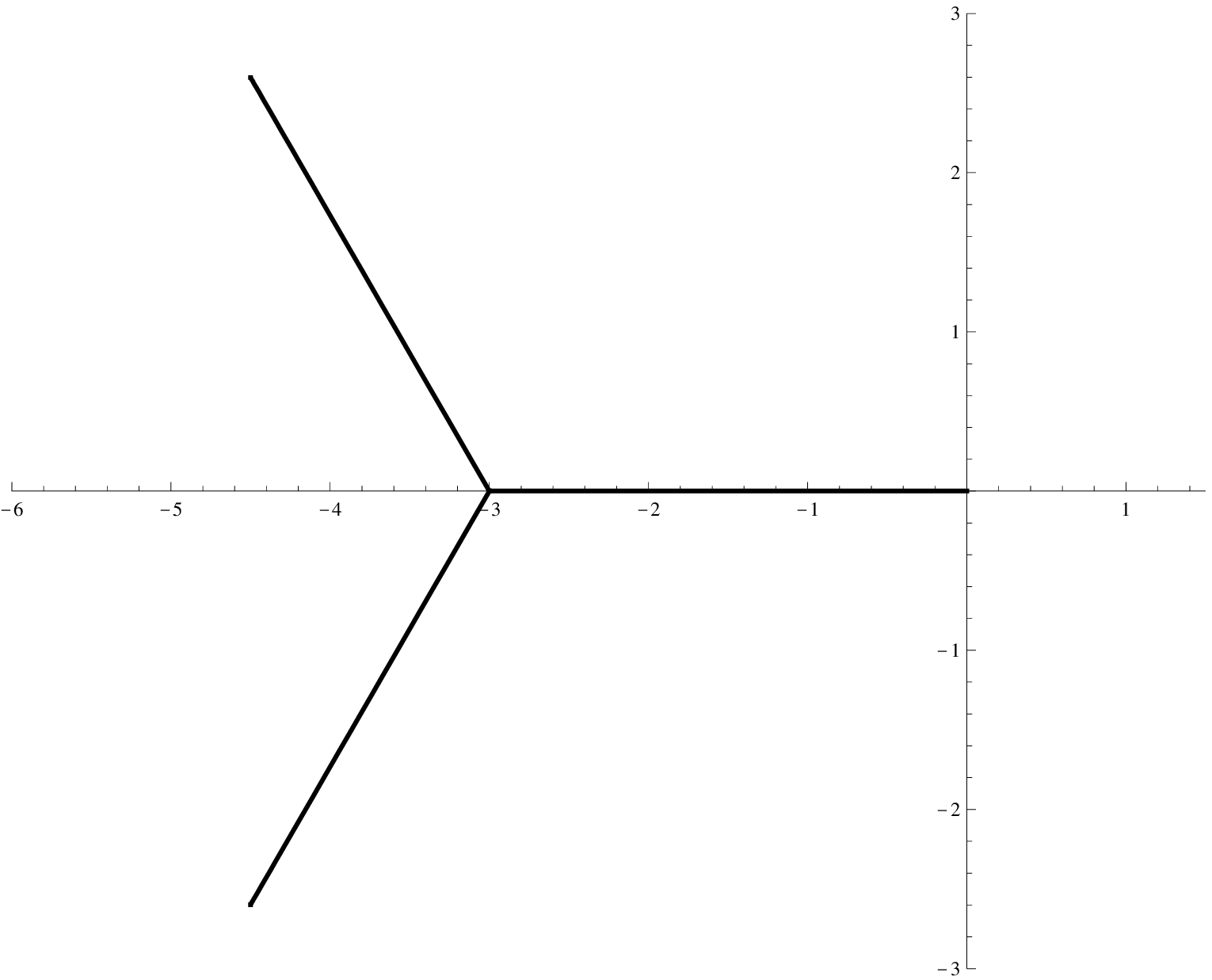}}
\end{center}
\caption{Image by $J_9$}\label{Im-J9B}
\end{figure}

For $k \leqslant 1000$, we can prove that all of the zeros of $E_{k, 3}^{\infty}$ lie on the lower arcs of $\partial \mathbb{F}_3$ by numerical calculation, then we have all of the zeros of $E_{k, 9}^{\infty}$ in the lower arcs of $\partial \mathbb{F}_9$. Thus, this case is very interesting. Though $J_9$ does not take real value on the some arcs of $\partial \mathbb{F}_9$, all of the zeros of $E_{k, 9}^{\infty}$ seems to lie on the lower arcs.\\

\paragraph{\bf Location of the zeros of Hecke type Faber Polynomial}
For $m \leqslant 200$, we can prove that all of the zeros of $F_{m, 9}$ lie on the lower arcs of $\partial \mathbb{F}_9$ by numerical calculation.\\

\clearpage

\section{Level $10$}

We have $\Gamma_0(10)+$, $\Gamma_0(10)+10=\Gamma_0^{*}(10)$, $\Gamma_0(10)+5$, $\Gamma_0(10)+2$, and $\Gamma_0(10)-=\Gamma_0(10)$. We have $W_{10} = \left(\begin{smallmatrix} 0 & - 1 / \sqrt{10} \\ \sqrt{10} & 0 \end{smallmatrix}\right)$, $W_{10, 2} := \left(\begin{smallmatrix} -\sqrt{2} & -1/\sqrt{2} \\ 5\sqrt{2} & 2 \sqrt{2} \end{smallmatrix}\right)$, and $W_{10, 5} := \left(\begin{smallmatrix} -\sqrt{5} & -3/\sqrt{5}\\ 2\sqrt{5} & \sqrt{5}\end{smallmatrix}\right)$.\\

\subsection{$\Gamma_0(10)+$}

We have $\Gamma_0(10)+ = \langle \left( \begin{smallmatrix} 1 & 1 \\ 0 & 1 \end{smallmatrix} \right), \: W_{10}, \: W_{10, 5} \rangle$.\\

\paragraph{\bf Location of the zeros of the Eisenstein series}
For $k \leqslant 800$, we can prove that all of the zeros of $E_{k, 10+}$ lie on the lower arcs of $\partial \mathbb{F}_{10+}$ by numerical calculation.
\begin{figure}[hbtp]
\begin{center}
\includegraphics[width=1.5in]{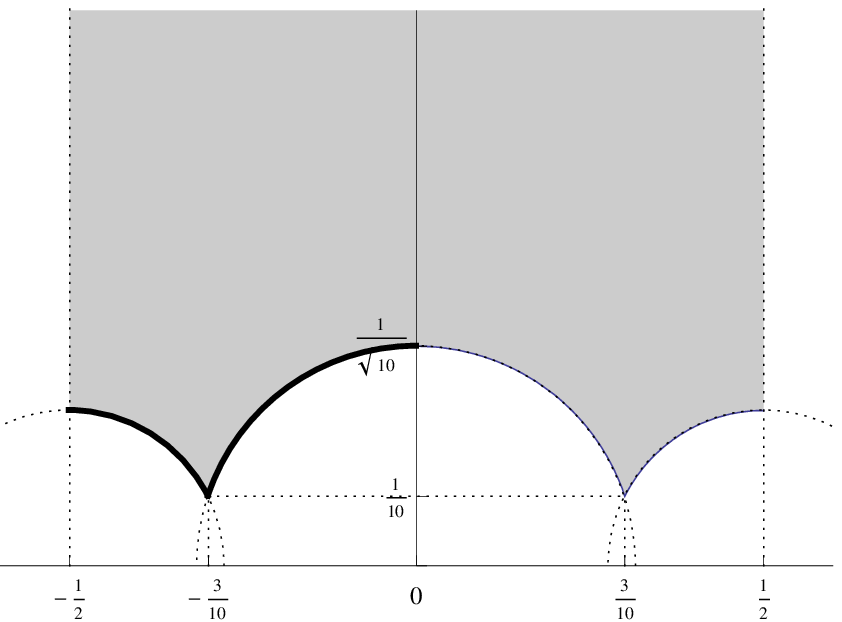}
\end{center}
\caption{Location of the zeros of the Eisenstein series}
\end{figure}

\paragraph{\bf Location of the zeros of Hecke type Faber Polynomial}
For $m \leqslant 200$, we can prove that all of the zeros of $F_{m, 10+}$ lie on the lower arcs of $\partial \mathbb{F}_{10+}$ by numerical calculation.\\

\subsection{$\Gamma_0(10)+10 = \Gamma_0^{*}(10)$}

We have $\Gamma_0^{*}(10) = \langle \left( \begin{smallmatrix} 1 & 1 \\ 0 & 1 \end{smallmatrix} \right), \: W_{10}, \: \left(\begin{smallmatrix} -3 & -1 \\ 10 & 3 \end{smallmatrix} \right), \: \left( \begin{smallmatrix} 9 & 4 \\ 20 & 9 \end{smallmatrix} \right) \rangle$ and $\gamma_{-1/2} = W_{10, 5}$.\\

\paragraph{\bf Location of the zeros of the Eisenstein series}
Since $W_{10, 5}^{- 1} \Gamma_0^{*}(10) W_{10, 5} = \Gamma_0^{*}(10)$, we have
\begin{equation}
E_{k, 10+10}^{-1/2}(W_{10, 5} z) = (2 \sqrt{5} z + \sqrt{5})^k E_{k, 10+10}^{\infty}(z).
\end{equation}
Furthermore, we have
\begin{align*}
E_{k, 10+10}^{-1/2} (e^{i \theta'} / (3 \sqrt{10}) - 1/3) &= (- \sqrt{5} e^{i \theta} - \sqrt{2})^k E_{k, 10+10}^{\infty}(e^{i \theta} / \sqrt{10}),\\
E_{k, 10+10}^{-1/2} (e^{i \theta'} / \sqrt{10}) &= (\sqrt{2} (e^{i \theta} + 1) / 3)^k E_{k, 10+10}^{\infty}(e^{i \theta} / (3 \sqrt{10}) - 1/3),\\
E_{k, 10+10}^{-1/2} (i \tan(\theta/2) / 2) &= ((e^{i \theta} + 1) / (2 \sqrt{5}))^k E_{k, 10+10}^{\infty}((e^{i \theta} - 9) / 20),
\end{align*}
where $e^{i \theta'} = (- 2 \sqrt{10} - 7 \cos\theta + 3 i \sin\theta) / (7 + 2 \sqrt{10} \cos\theta)$.

For $k \leqslant 500$, we can prove that all of the zeros of $E_{k, 10+10}^{\infty}$ lie on the lower arcs of $\partial \mathbb{F}_{10+10}$ by numerical calculation.

\begin{figure}[htbp]
\begin{center}
{{$E_{k, 10+10}^{\infty}$}\includegraphics[width=1.5in]{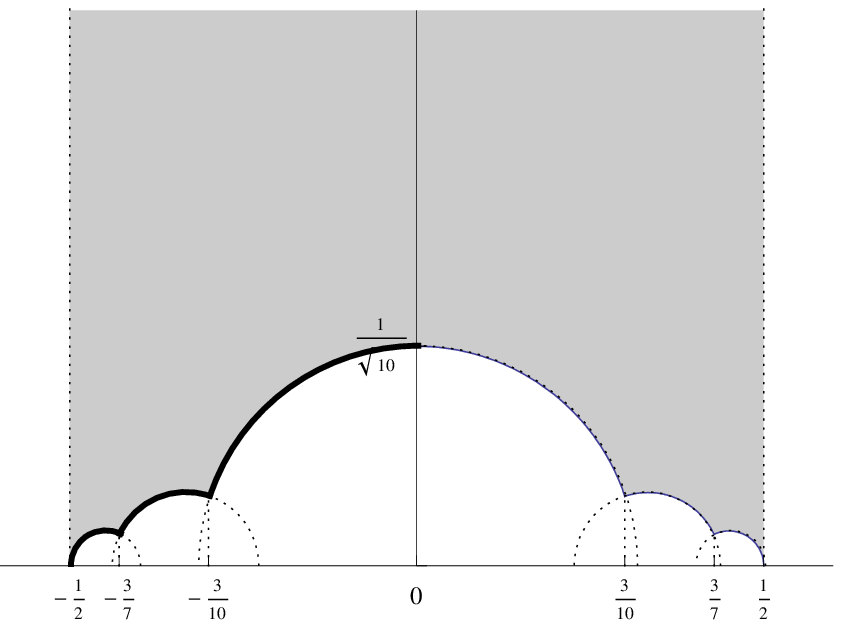}}
\quad \quad
{{$E_{k, 10+10}^{-1/2}$}\includegraphics[width=1.5in]{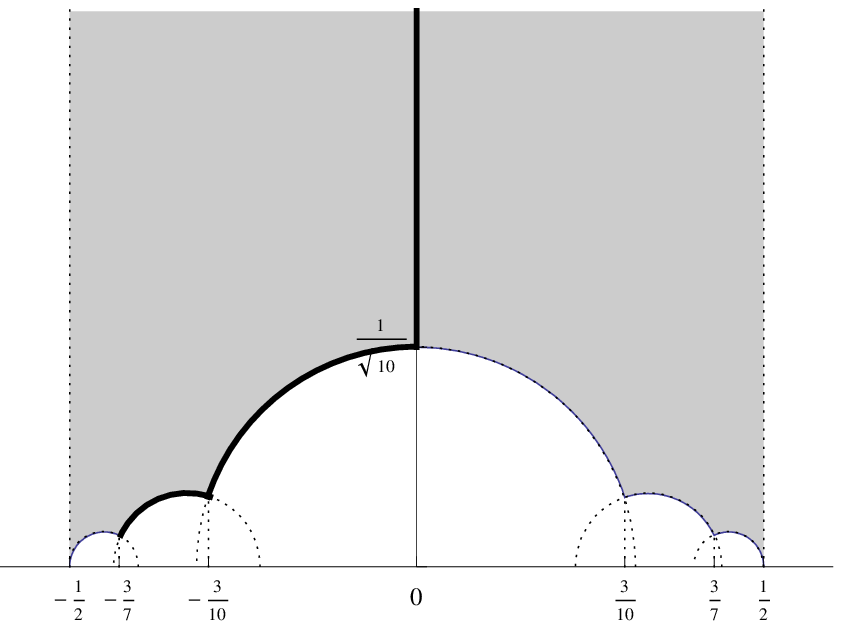}}
\end{center}
\caption{Location of the zeros of the Eisenstein series}
\end{figure}

\paragraph{\bf Location of the zeros of Hecke type Faber Polynomial}
For every odd integer $m \leqslant 200$ but $m = 7, 9, 11$, we can prove that all of the zeros of $F_{m, 10+10}$ lie on the lower arcs of $\partial \mathbb{F}_{10+10}$ by numerical calculation. On the other hand, by numerical calculation, for $m = 7, 9, 11$, we can prove that all but two of the zeros of $F_{m, 10+10}$ lie on the lower arcs of $\partial \mathbb{F}_{10+10}$, and two of the zeros of $F_{m, 10+10}$ do not lie on $\partial \mathbb{F}_{10+10}$. For the other cases where $m$ is even and $m \leqslant 200$, by numerical calculation, we can prove that all but one of the zeros of $F_{m, 10+10}$ lie on the lower arcs of $\partial \mathbb{F}_{10+10}$, and one of the zeros of $F_{m, 10+10}$ lies on $\partial \mathbb{F}_{10+10}$ but does not on the lower arcs.\\

\subsection{$\Gamma_0(10)+5$}

We have $\Gamma_0(10)+5 = \langle \left( \begin{smallmatrix} 1 & 1 \\ 0 & 1 \end{smallmatrix} \right), \: W_{10, 5}, \: \left( \begin{smallmatrix} 1 & 0 \\ 10 & 1 \end{smallmatrix} \right), \: \left(\begin{smallmatrix} -3 & -1 \\ 10 & 3 \end{smallmatrix} \right) \rangle$ and $\gamma_0 = W_{10}$.\\

\paragraph{\bf Location of the zeros of the Eisenstein series}
Since $W_{10}^{- 1} (\Gamma_0(10)+5) W_{10} = \Gamma_0(10)+5$, we have
\begin{equation}
E_{k, 10+5}^0(W_{10} z) = (\sqrt{10} z)^k E_{k, 10+5}^{\infty}(z).
\end{equation}
Furthermore, we have
\begin{align*}
E_{k, 10+5}^0 (-1/2 + i / (2 \tan(\theta/2))) &= ((e^{i \theta} - 1) / \sqrt{10})^k E_{k, 10+5}^{\infty}((e^{i \theta} - 1) / 10),\\
E_{k, 10+5}^0 (e^{i \theta'} / (2 \sqrt{5}) - 1/2) &= (- (\sqrt{5} e^{i \theta} -1) / (2 \sqrt{2}))^k E_{k, 10+5}^{\infty}(e^{i \theta} / (4 \sqrt{5}) - 1/4),\\
E_{k, 10+5}^0 (e^{i \theta'} / (4 \sqrt{5}) - 1/4) &= (- (\sqrt{5} e^{i \theta} -1) / \sqrt{2})^k E_{k, 10+5}^{\infty}(e^{i \theta} / (2 \sqrt{5}) - 1/2),
\end{align*}
where $e^{i \theta'} = (2 \sqrt{5} - 6 \cos\theta + 4 i \sin\theta) / (6 - 2 \sqrt{5} \cos\theta)$.

For $k \leqslant 450$, we can prove that all of the zeros of $E_{k, 10+5}^{\infty}$ lie on the lower arcs of $\partial \mathbb{F}_{10+5}$ by numerical calculation.

\begin{figure}[hbtp]
\begin{center}
{{$E_{k, 10+5}^{\infty}$}\includegraphics[width=1.5in]{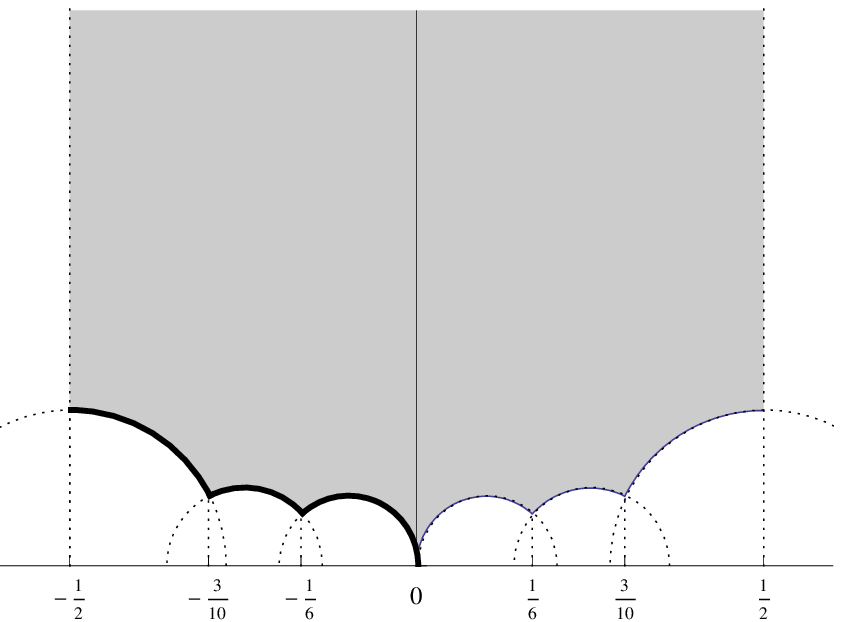}}
\quad \quad
{{$E_{k, 10+5}^0$}\includegraphics[width=1.5in]{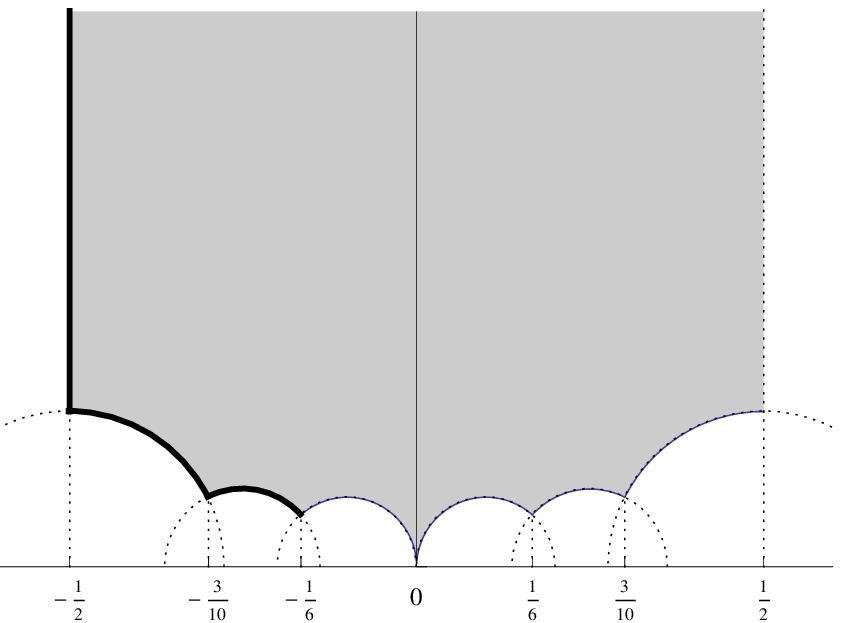}}
\end{center}
\caption{Location of the zeros of the Eisenstein series}
\end{figure}

\paragraph{\bf Location of the zeros of Hecke type Faber Polynomial}
For $m \leqslant 200$, we can prove that all of the zeros of $F_{m, 10+5}$ lie on the lower arcs of $\partial \mathbb{F}_{10+5}$ by numerical calculation.\\

\subsection{$\Gamma_0(10)+2$}

We have $\Gamma_0(10)+2 = \langle \left( \begin{smallmatrix} 1 & 1 \\ 0 & 1 \end{smallmatrix} \right), \: W_{10, 2}, \: \left( \begin{smallmatrix} 1 & 0 \\ 10 & 1 \end{smallmatrix} \right), \: \left(\begin{smallmatrix} 9 & 2 \\ 40 & 9 \end{smallmatrix}\right) \rangle$ and $\gamma_0 = W_{10}$.\\

\paragraph{\bf Location of the zeros of Eisenstein series}

Since $W_{10}^{- 1} (\Gamma_0(10)+2) W_{10} = \Gamma_0(10)+2$, we have
\begin{equation}
E_{k, 10+2}^0(W_{10} z) = (\sqrt{10} z)^k E_{k, 10+2}^{\infty}(z).
\end{equation}
Furthermore, we have
\begin{align*}
E_{k, 10+2}^0 (- 1 / 2 + i / (2 \tan\theta/2)) &= ((e^{i \theta} - 1) / \sqrt{10})^k E_{k, 10+2}^{\infty}((e^{i \theta} - 1) / 10),\\
E_{k, 10+2}^0 (e^{i \theta'} / \sqrt{2} + 1) &= (- (\sqrt{2} e^{i \theta} + 1) / \sqrt{5})^k E_{k, 10+2}^{\infty}(e^{i \theta} / (5 \sqrt{2}) - 2/5),\\
E_{k, 10+2}^0 (e^{i \theta''} / \sqrt{2} - 1) &= ((e^{i \theta} + \sqrt{2}) / \sqrt{5})^k E_{k, 10+2}^{\infty}(e^{i \theta} / (5 \sqrt{2}) + 1/5),
\end{align*}
where $e^{i \theta'} = (- 2 \sqrt{2} - 3 \cos\theta + i \sin\theta) / (3 + 2 \sqrt{2} \cos\theta)$ and $e^{i \theta''} = ( - 2 \sqrt{2} + 3 \cos\theta + i \sin\theta) / (3 - 2 \sqrt{2} \cos\theta)$.
\begin{figure}[hbtp]
\begin{center}
{{$E_{k, 10+2}^{\infty}$}\includegraphics[width=1.5in]{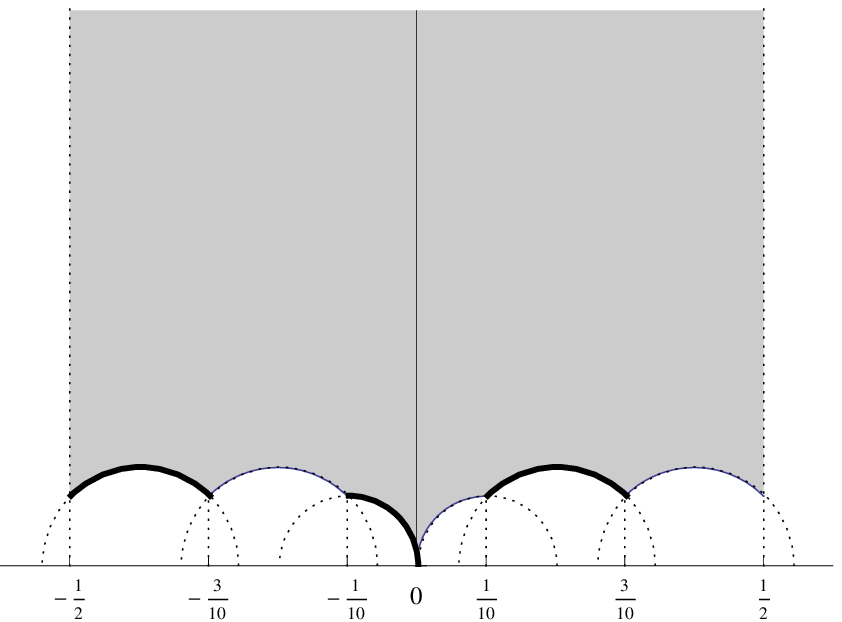}}
\quad \quad
{{$E_{k, 10+2}^{0}$}\includegraphics[width=1.5in]{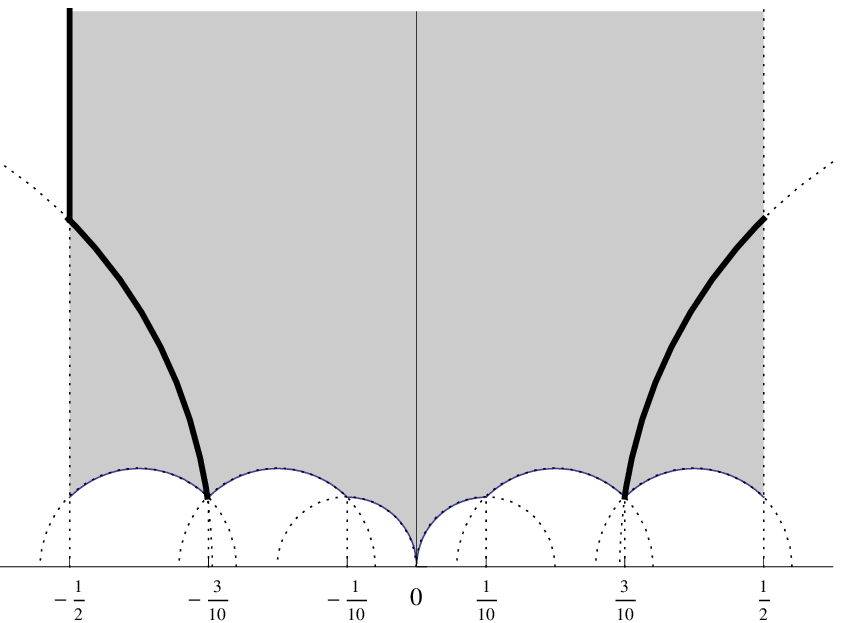}}
\end{center}
\caption{Neighborhood of location of the zeros of the Eisenstein series}
\end{figure}

\begin{figure}[hbtp]
\begin{center}
{{\small Lower arcs of $\partial \mathbb{F}_{10+2}$}\includegraphics[width=2.5in]{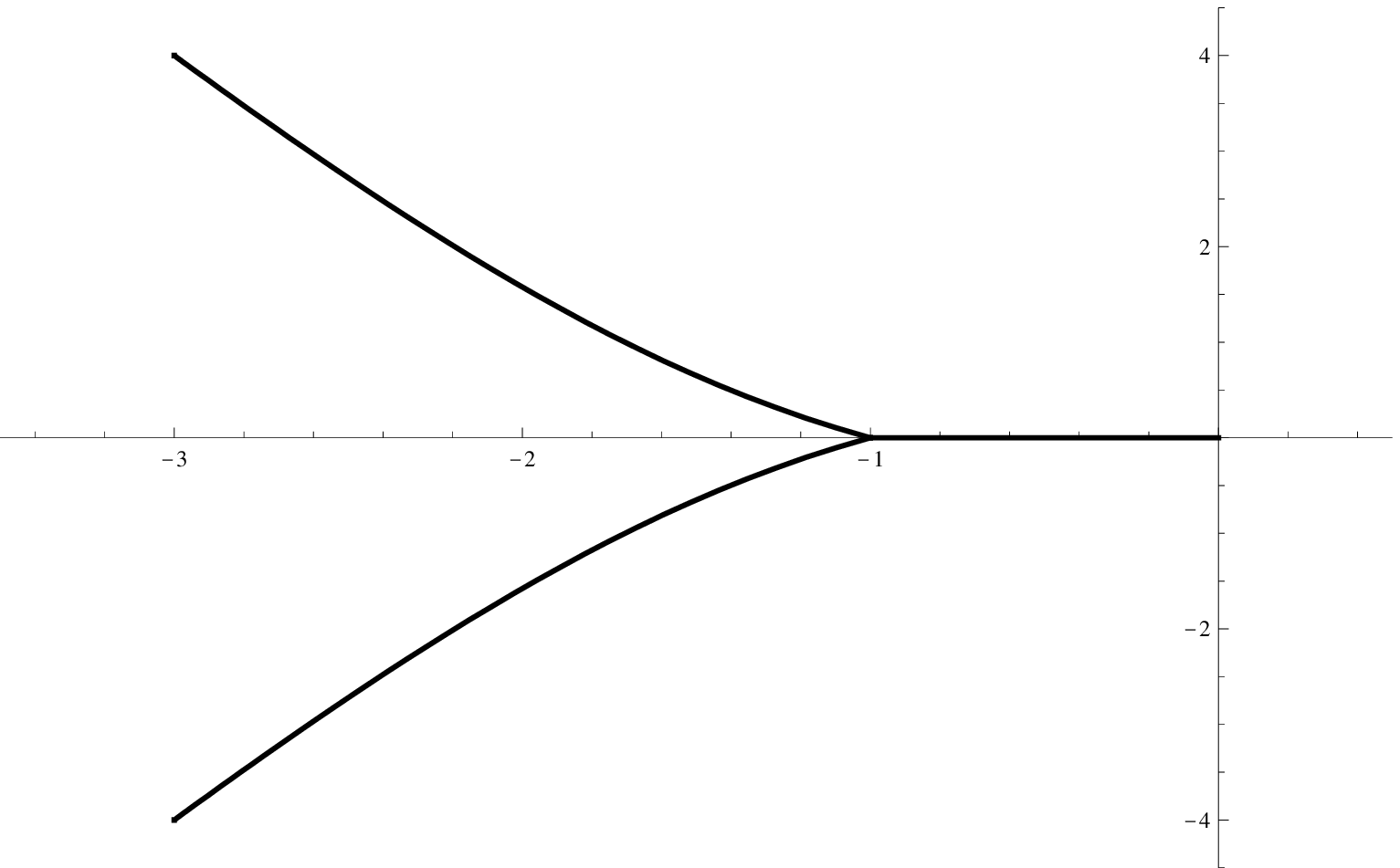}}
\end{center}
\caption{Image by $J_{10+2}$}\label{Im-J10C}
\end{figure}

\begin{figure}[hbtp]
\begin{center}
\includegraphics[width=6.3in]{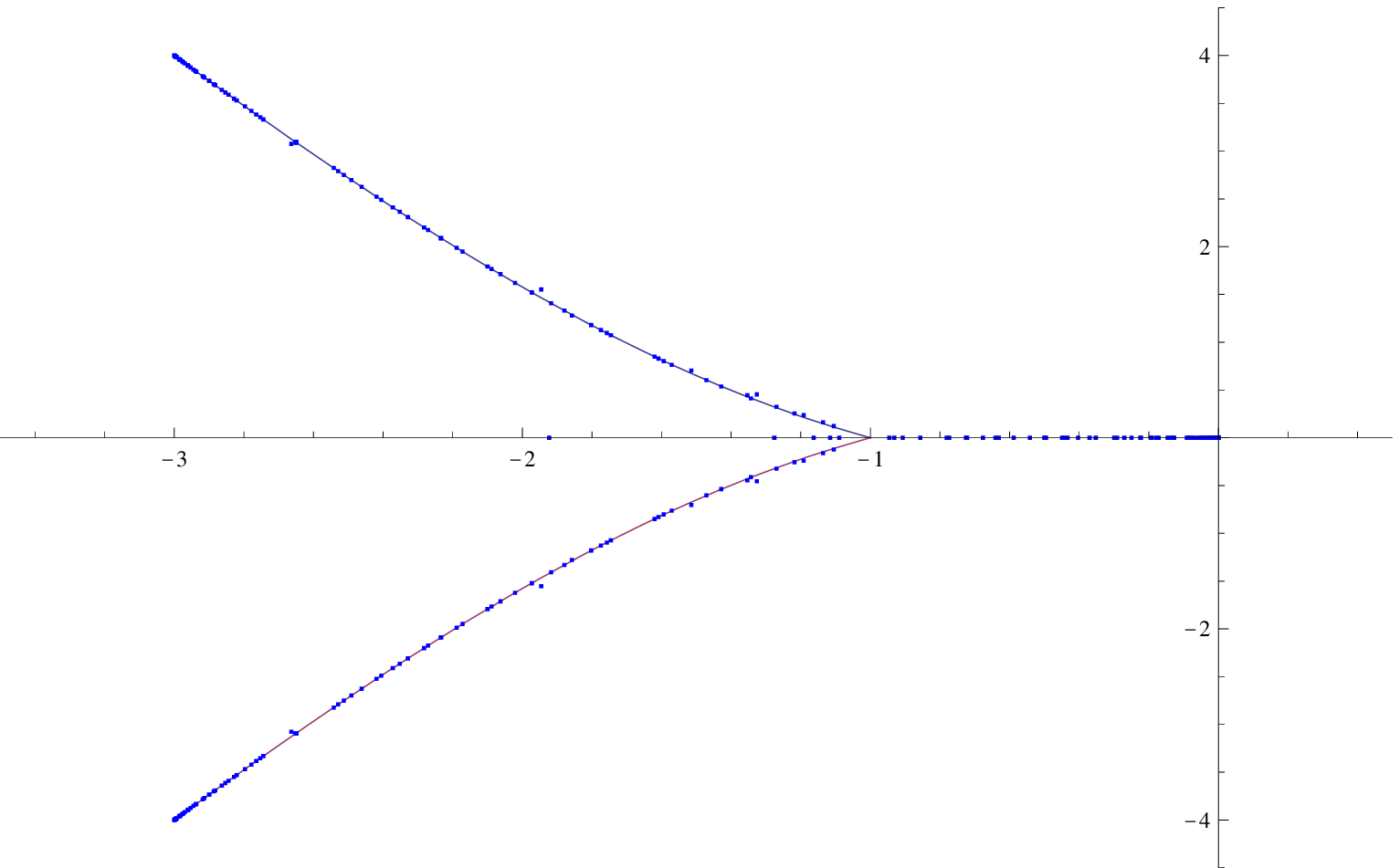}\\
The zeros of $E_{k, 10+2}^{\infty}$ for $4 \leqslant k \leqslant 40$\\
\includegraphics[width=6.3in]{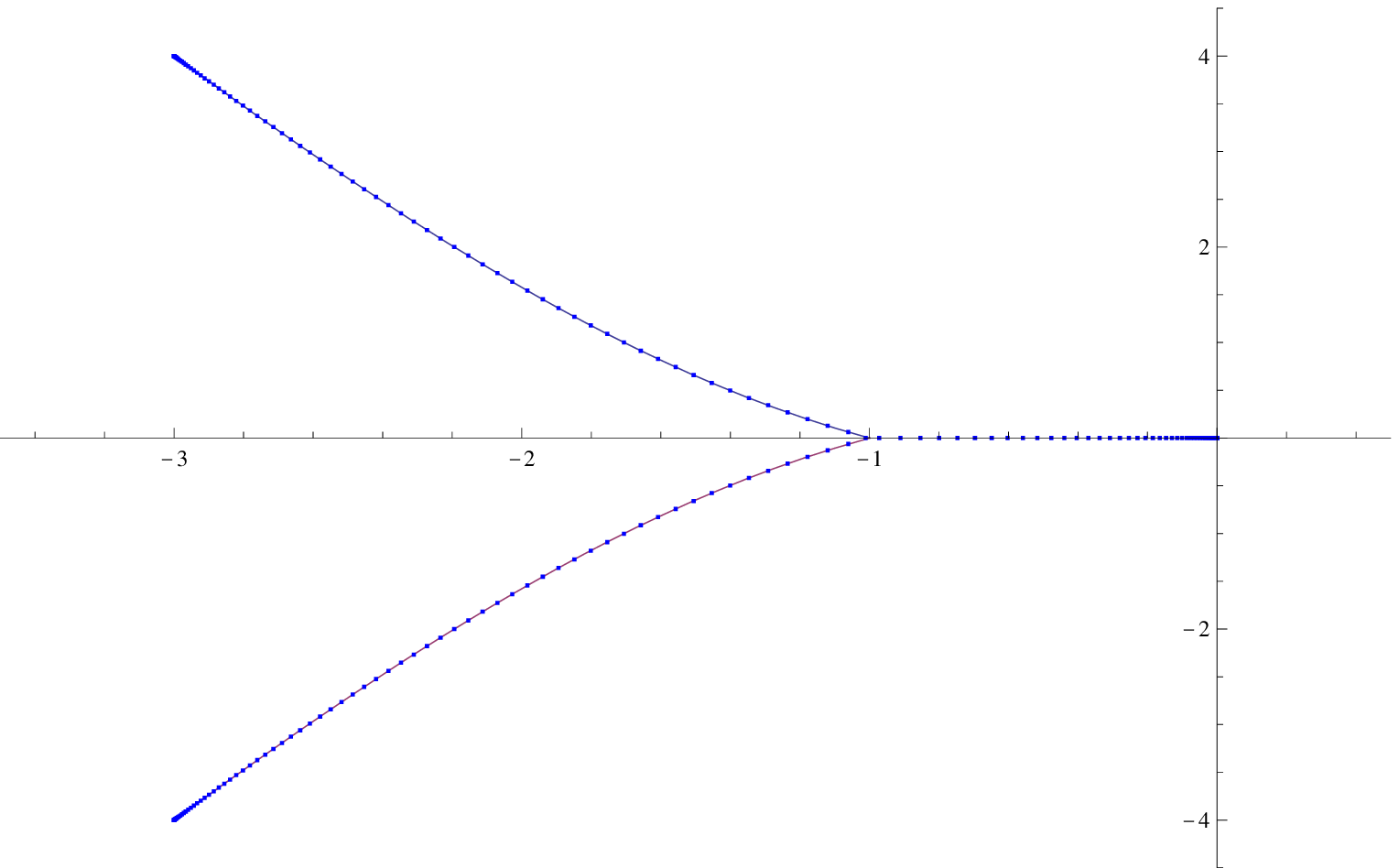}\\
The zeros of $E_{1000, 10+2}^{\infty}$
\end{center}
\caption{Image by $J_{10+2}$}\label{Im-J10Cz}
\end{figure}

Now, we can observe that the zeros of $E_{k, 10+2}^{\infty}$ do not lie on the arcs of $\partial \mathbb{F}_{10+2}$ for small weight $k$ by numerical calculation. However, when the weight $k$ increases, then the location of the zeros seems to approach to lower arcs of $\partial \mathbb{F}_{10+2}$. (See Figure \ref{Im-J10Cz})

\begin{figure}[hbtp]
\begin{center}
\includegraphics[width=6.3in]{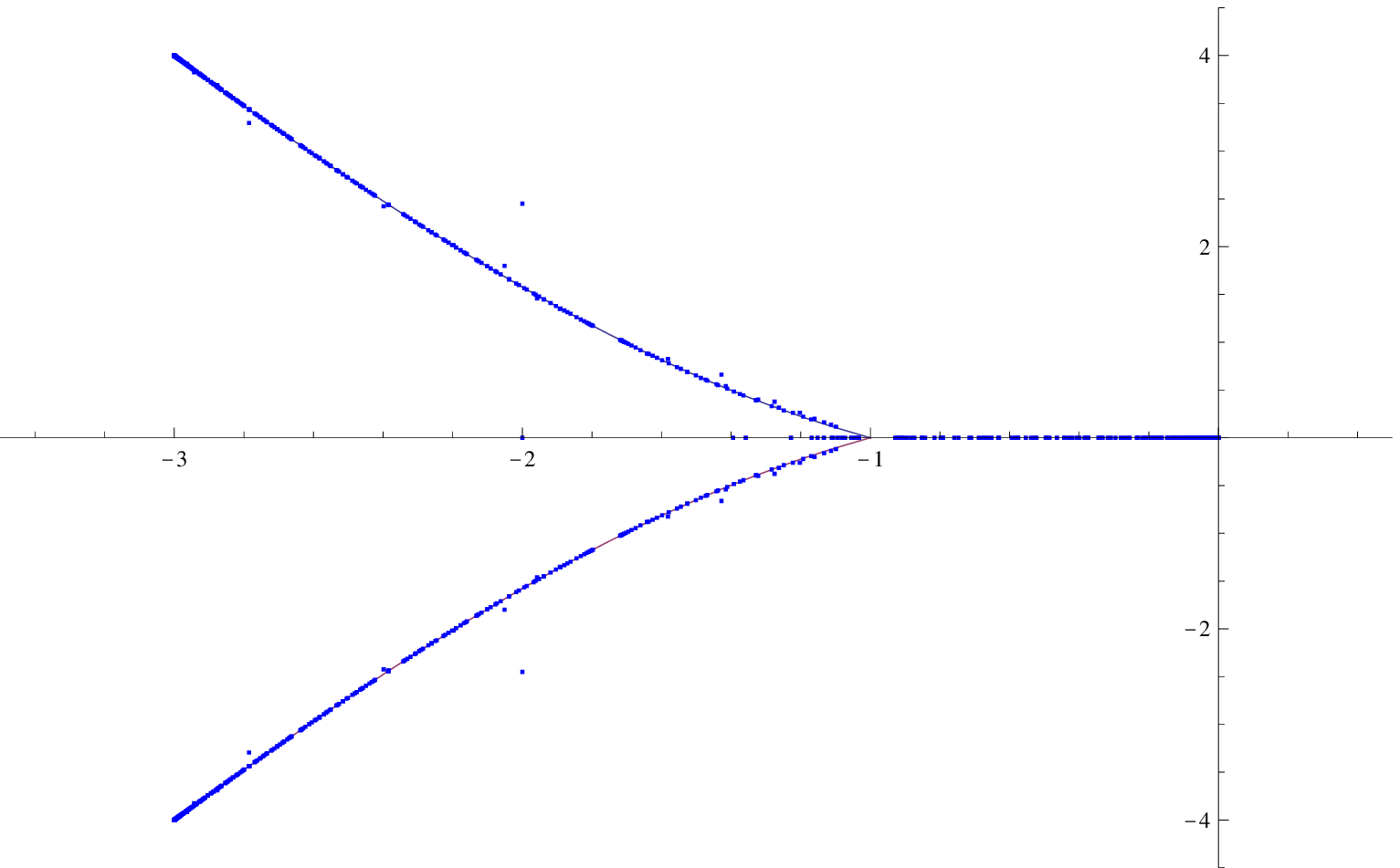}\\
The zeros of $F_{m, 10+2}$ for $m \leqslant 40$\\
\includegraphics[width=6.3in]{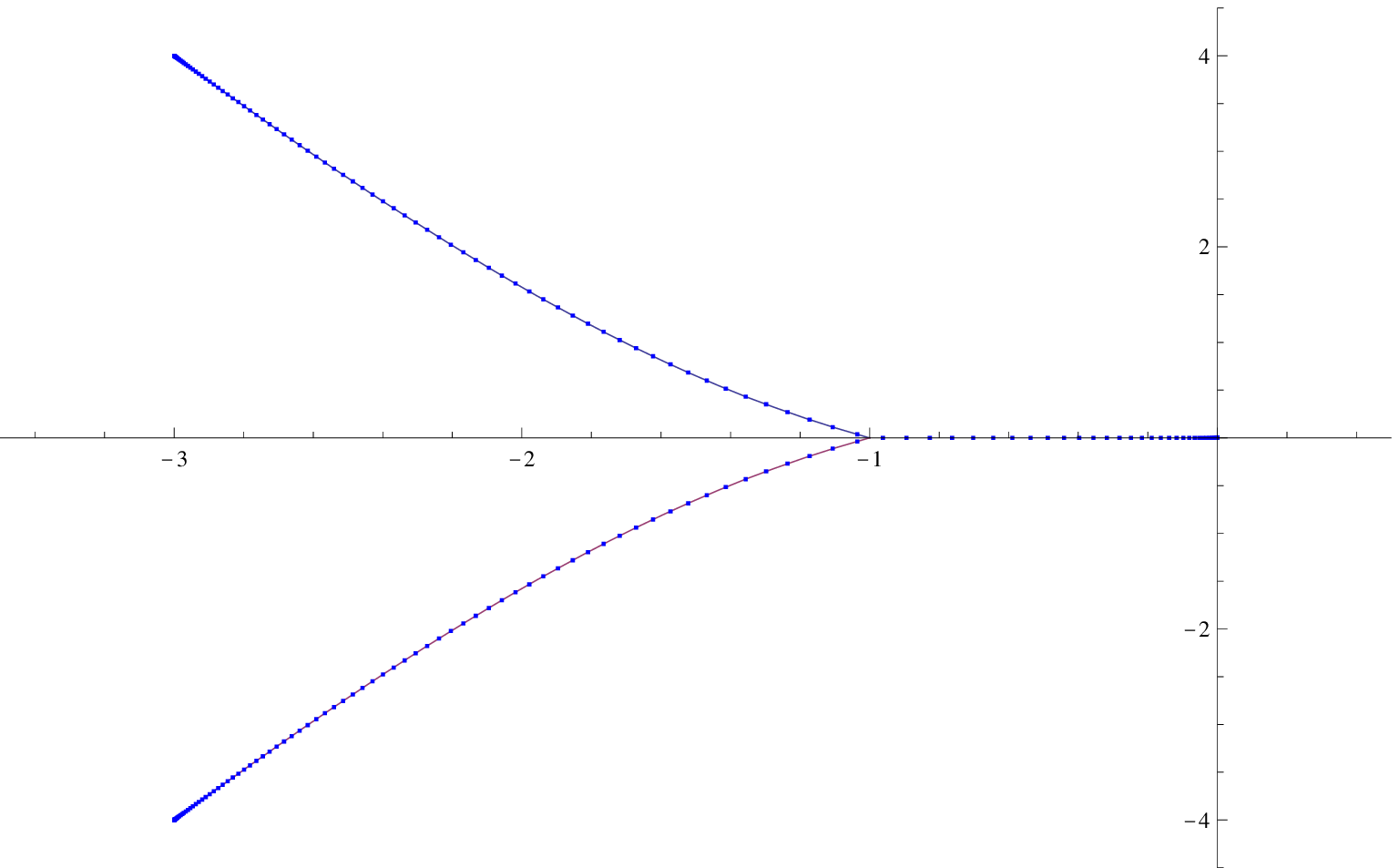}\\
The zeros of $F_{200, 10+2}$
\end{center}
\caption{Image by $J_{10+2}$}\label{Im-J10Chz}
\end{figure}

\paragraph{\bf Location of the zeros of Hecke type Faber Polynomial}
We can observe that some zeros of $F_{m, 10+2}$ do not lie on the lower arcs of $\partial \mathbb{F}_{10+2}$ for small weight $m$ by numerical calculation. However, when the weight $m$ increases, then the location of the zeros seems to approach to lower arcs of $\partial \mathbb{F}_{10+2}$. (see Figure \ref{Im-J10Chz})\\

\newpage

\subsection{$\Gamma_0(10)$}

We have $\Gamma_0(10) = \langle \left( \begin{smallmatrix} 1 & 1 \\ 0 & 1 \end{smallmatrix} \right), \: \left( \begin{smallmatrix} 1 & 0 \\ 10 & 1 \end{smallmatrix} \right), \: \left( \begin{smallmatrix} -3 & -1 \\ 10 & 3 \end{smallmatrix} \right), \: \left( \begin{smallmatrix} 3 & -1 \\ 10 & -3 \end{smallmatrix} \right), \: \left( \begin{smallmatrix} 9 & 4 \\ 20 & 9 \end{smallmatrix} \right) \rangle$, $\gamma_0 = W_{10}$, $\gamma_{-1/2} = W_{10, 5}$, and $\gamma_{-1/5} = W_{10, 2}$.\\

\paragraph{\bf Location of the zeros of the Eisenstein series}
Since $W_{10}^{- 1} \Gamma_0(10) W_{10} = W_{10, 5}^{- 1} \Gamma_0(10) W_{10, 5} = W_{10, 2}^{- 1} \Gamma_0(10) W_{10, 2} = \Gamma_0(10)$, we have
\begin{equation*}
(\sqrt{10} z)^{-k} E_{k, 10}^0(W_{10} z) = (2 \sqrt{5} z + \sqrt{5})^{-k} E_{k, 10}^{-1/2}(W_{10, 5} z)
 = (5 \sqrt{2} z + 2 \sqrt{2})^{-k} E_{k, 10}^{-1/5}(W_{10, 2} z) = E_{k, 10}^{\infty}(z).
\end{equation*}
Furthermore, we have
\begin{align*}
E_{k, 10}^0 (-1/2 + i / (2 \tan(\theta/2))) &= ((e^{i \theta} - 1) / \sqrt{10})^k E_{k, 10}^{\infty}((e^{i \theta} - 1) / 10),\\
E_{k, 10}^0 ((e^{i \theta_1} + 3) / 8) &= ((e^{i \theta} - 3) / \sqrt{10})^k E_{k, 10}^{\infty}((e^{i \theta} - 3) / 10),\\
E_{k, 10}^0 ((e^{i (\pi - \theta_1)} - 3) / 8) &= ((3 e^{i \theta} - 1) / \sqrt{10})^k E_{k, 10}^{\infty}((e^{i \theta} + 3) / 10),\\
E_{k, 10}^0 ((e^{i \theta_2} - 9) / 20) &= ((5 e^{i \theta} -1) / (2 \sqrt{10}))^k E_{k, 10}^{\infty}((e^{i \theta} - 9) / 20),\\
E_{k, 10}^{-1/2} ((e^{i \theta_3} - 1) / 10) &= ((3 e^{i \theta} + 2) / \sqrt{5})^k E_{k, 10}^{\infty}((e^{i \theta} - 1) / 10),\\
E_{k, 10}^{-1/2} ((e^{i \theta_4} - 1) / 6) &= ((e^{i \theta} + 2) / \sqrt{5})^k E_{k, 10}^{\infty}((e^{i \theta} - 3) / 10),\\
E_{k, 10}^{-1/2} ((e^{i (\pi - \theta_4)} - 1) / 6) &= ((e^{i \theta} + 1) /  \sqrt{5})^k E_{k, 10}^{\infty}((e^{i \theta} + 3) / 10),\\
E_{k, 10}^{-1/2} (i \tan(\theta/2) / 2) &= ((e^{i \theta} + 1) / (2 \sqrt{5}))^k E_{k, 10}^{\infty}((e^{i \theta} - 5) / 20),\\
E_{k, 10}^{-1/5} (-1/2 + i \tan(\theta/2) / 10) &= ((e^{i \theta} + 1) / \sqrt{2})^k E_{k, 10}^{\infty}((e^{i \theta} - 1) / 10),\\
E_{k, 10}^{-1/5} (-3/10 + i \tan(\theta/2) / 10) &= ((e^{i \theta} + 1) / \sqrt{2})^k E_{k, 10}^{\infty}((e^{i \theta} - 3) / 10),\\
E_{k, 10}^{-1/5} (3/10 + i \tan(\theta/2) / 10) &= ((e^{i \theta} + 1) / \sqrt{2})^k E_{k, 10}^{\infty}((e^{i \theta} + 3) / 10),\\
E_{k, 10}^{-1/5} (i \tan(\theta/2) / 5) &= ((e^{i \theta} - 1) / (2 \sqrt{2}))^k E_{k, 10}^{\infty}((e^{i \theta} - 9) / 20),
\end{align*}
where $e^{i \theta_1} = (-3 + 5 \cos\theta + 4 i \sin\theta) / (5 - 3 \cos\theta)$, $e^{i \theta_2} = (21 - 29 \cos\theta + 20 i \sin\theta) / (29 - 21 \cos\theta)$, $e^{i \theta_3} = (-12 -13 \cos\theta + 5 i \sin\theta) / (13 + 12 \cos\theta)$, and $e^{i \theta_4} = (4 + 5 \cos\theta + 3 i \sin\theta) / (5 + 4 \cos\theta)$.
\begin{figure}[hbtp]
\begin{center}
{{$E_{k, 10}^{\infty}$}\includegraphics[width=1.5in]{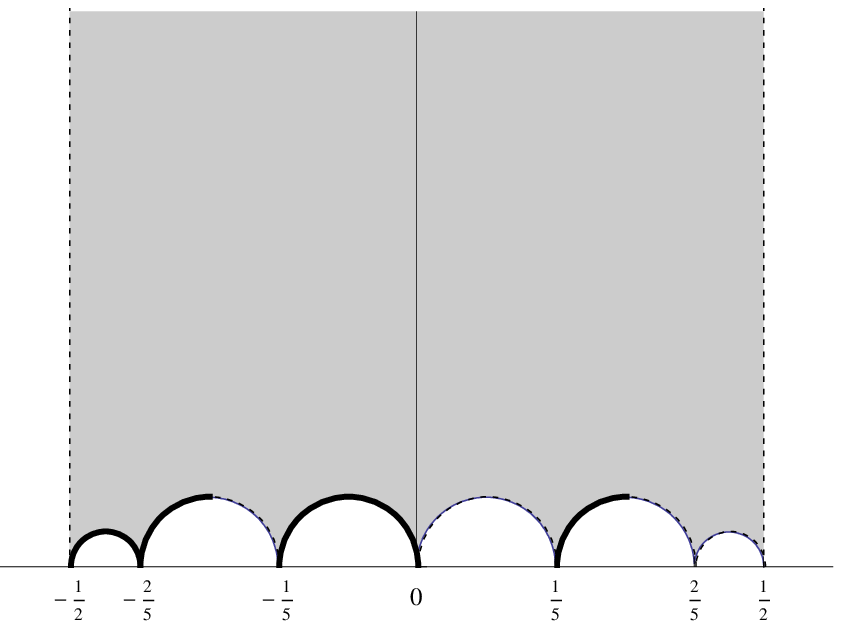}}
\quad \quad
{{$E_{k, 10}^0$}\includegraphics[width=1.5in]{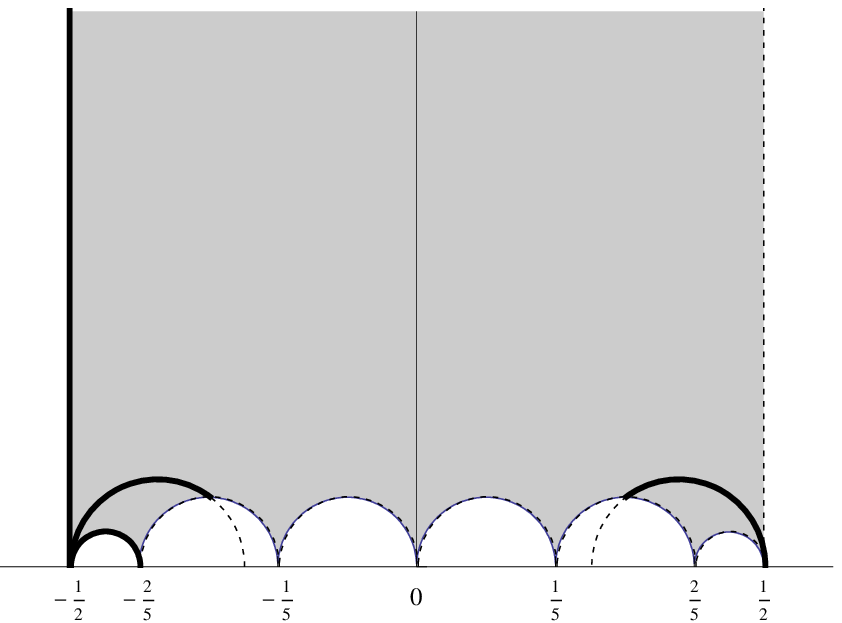}}\\

{{$E_{k, 10}^{-1/2}$}\includegraphics[width=1.5in]{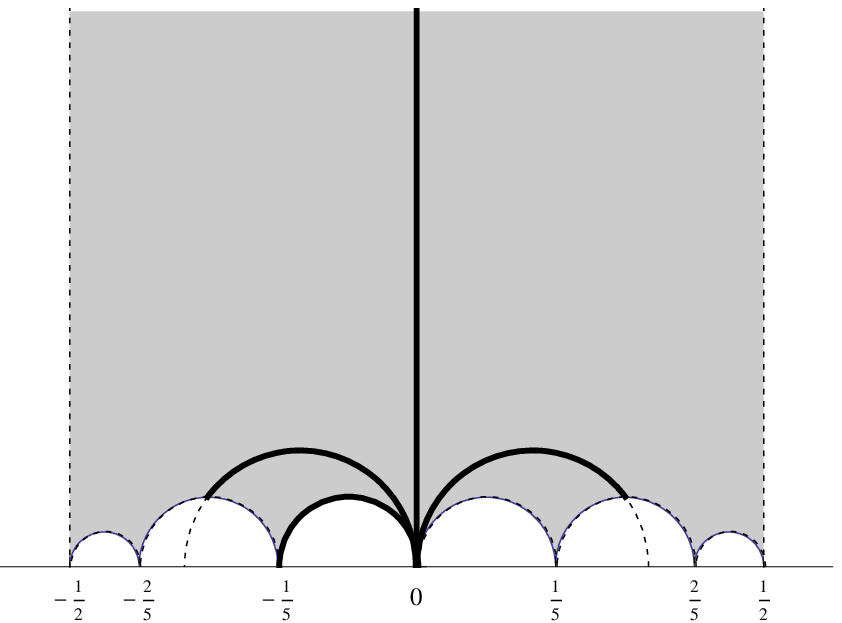}}
\quad \quad
{{$E_{k, 10}^{-1/5}$}\includegraphics[width=1.5in]{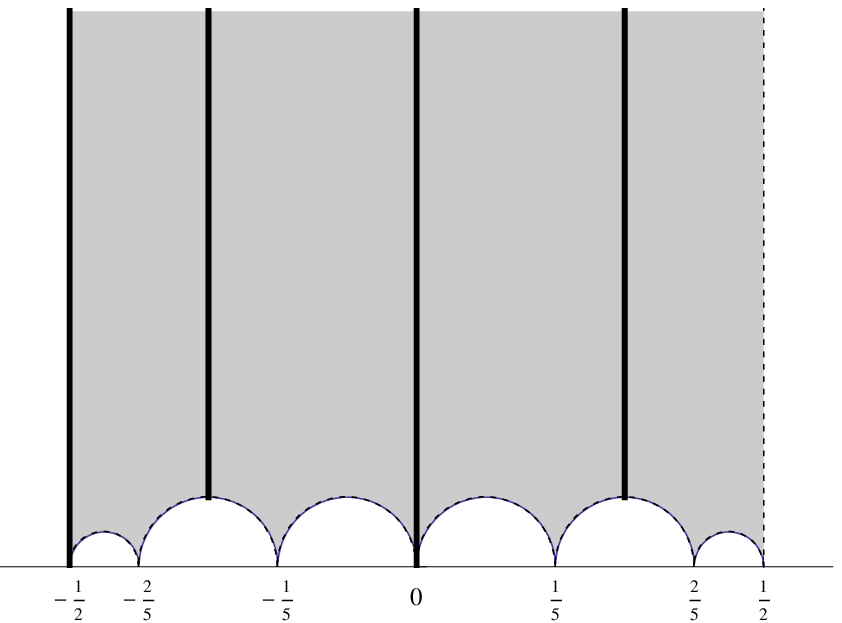}}
\end{center}
\caption{Location of the zeros of the Eisenstein series}
\end{figure}

\begin{figure}[hbtp]
\begin{center}
\includegraphics[width=6.3in]{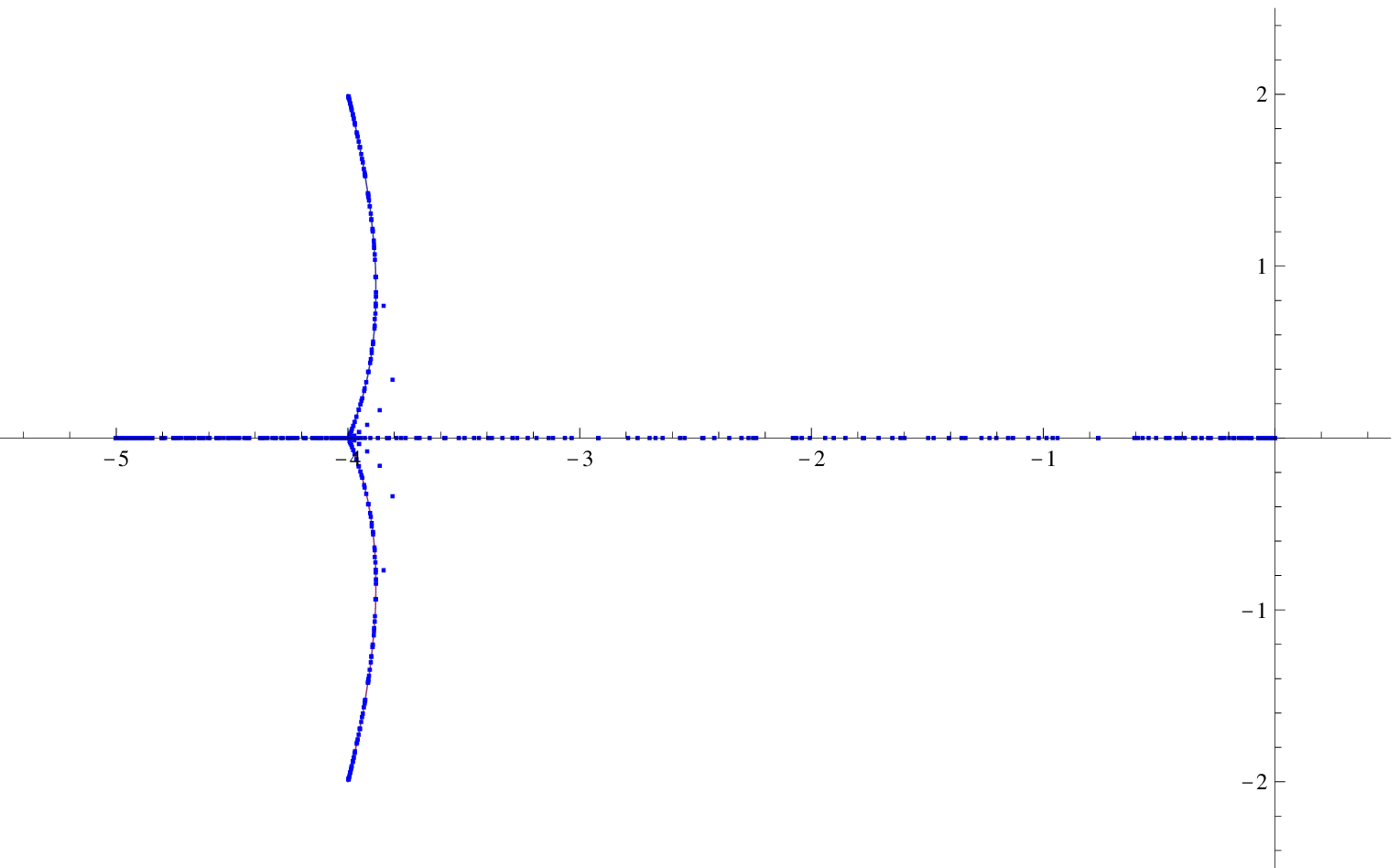}\\
The zeros of $E_{k, 10}^{\infty}$ for $4 \leqslant k \leqslant 40$\\
\includegraphics[width=6.3in]{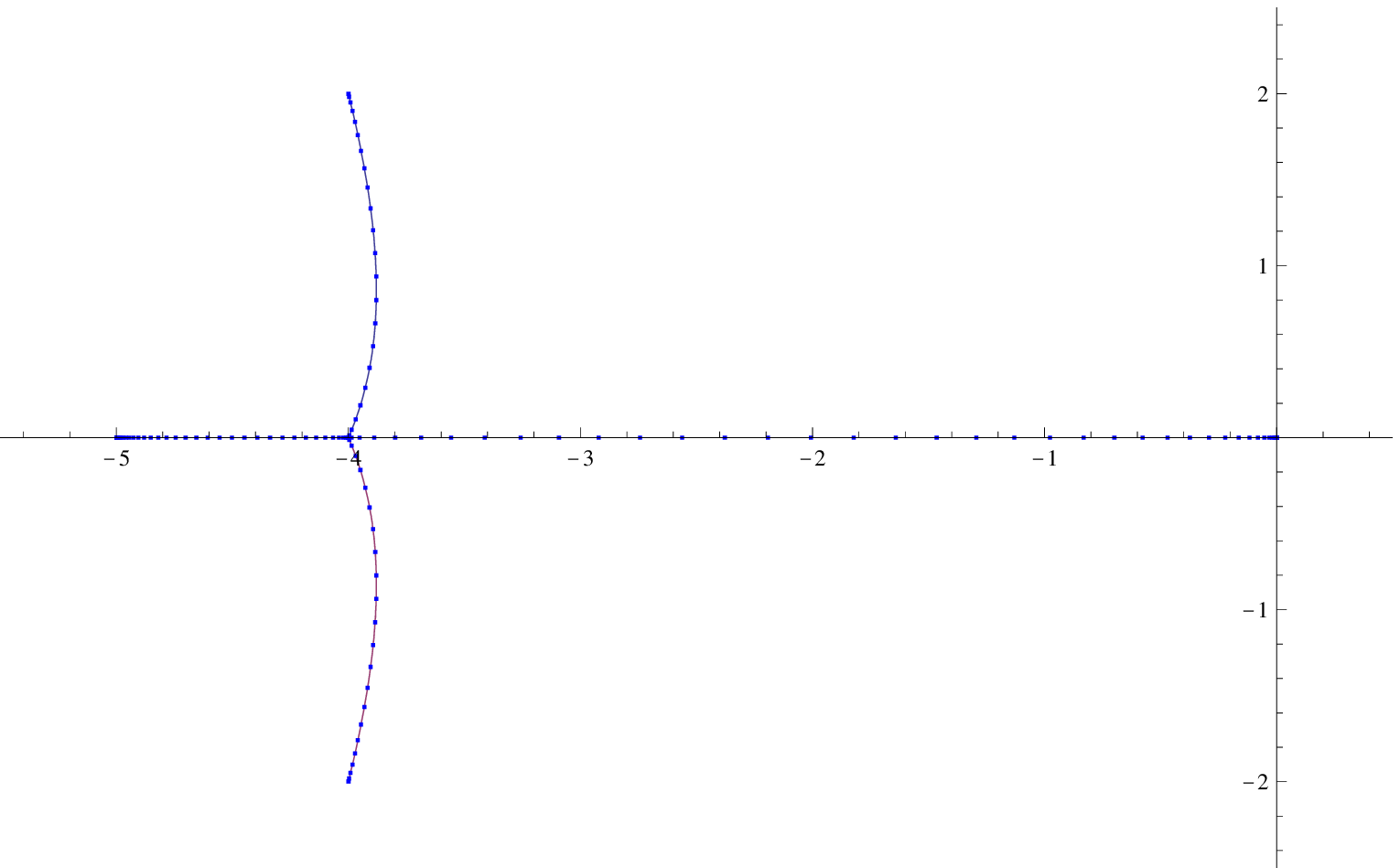}\\
The zeros of $E_{1000, 10}^{\infty}$
\end{center}
\caption{Image by $J_{10}$}\label{Im-J10Ez}
\end{figure}

Now, we can observe that the zeros of $E_{k, 10}^{\infty}$ do not lie on the arcs of $\partial \mathbb{F}_{10}$ for small weight $k$ by numerical calculation. However, when the weight $k$ increases, then the location of the zeros seems to approach to lower arcs of $\partial \mathbb{F}_{10}$. (See Figure \ref{Im-J10Ez})\\

\begin{figure}[hbtp]
\begin{center}
\includegraphics[width=6.3in]{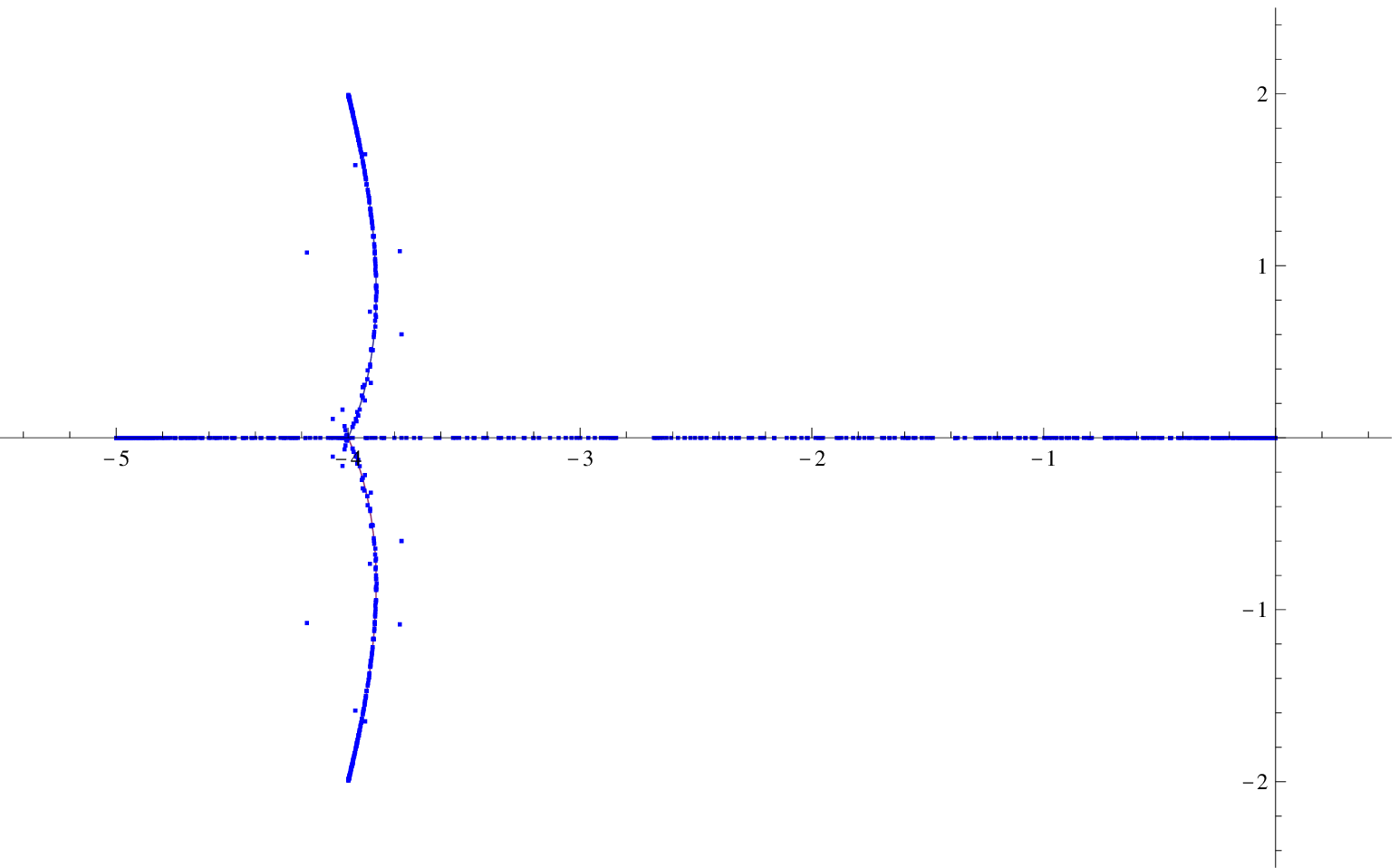}\\
The zeros of $F_{m, 10}$ for $m \leqslant 40$\\
\includegraphics[width=6.3in]{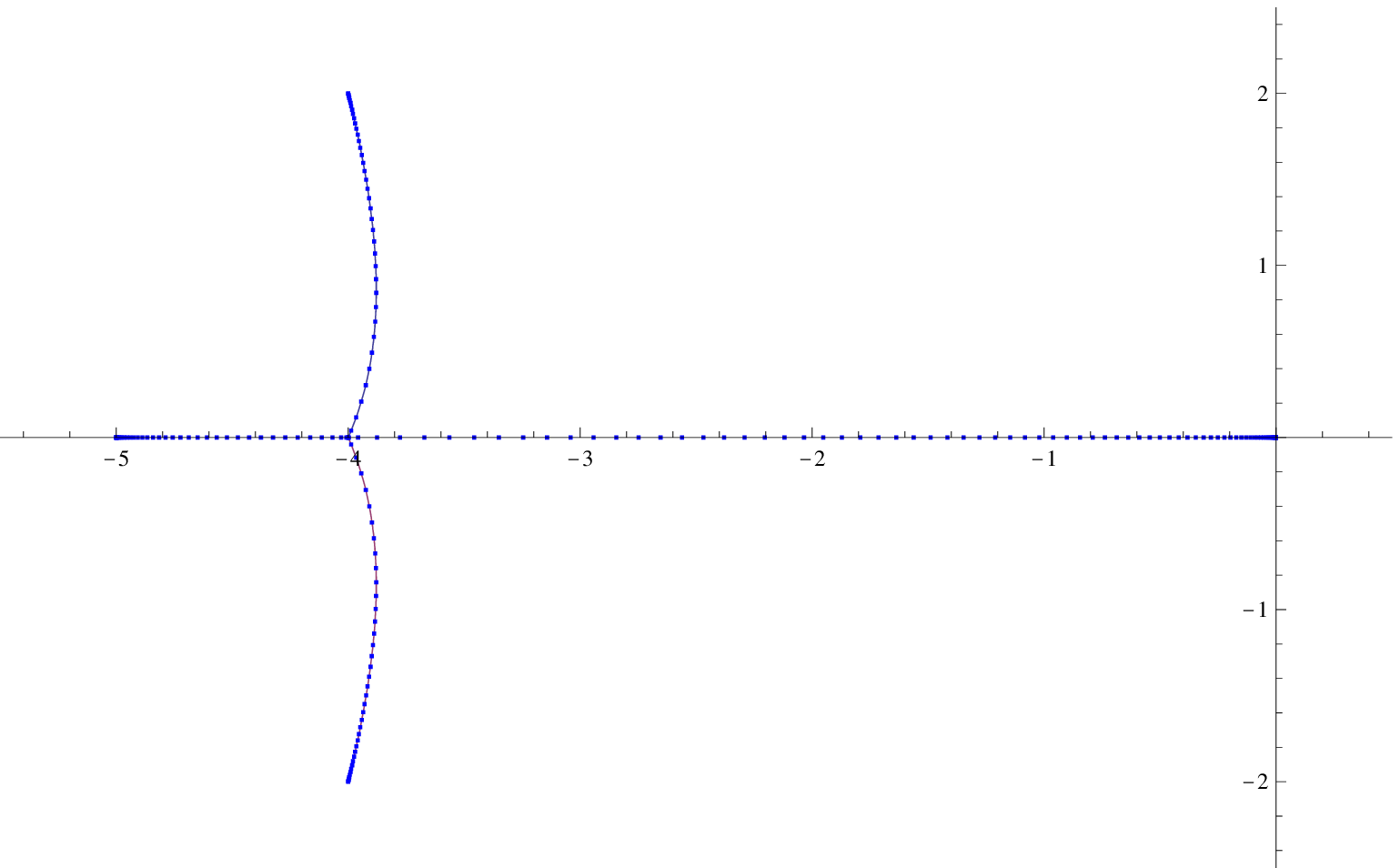}\\
The zeros of $F_{200, 10}$
\end{center}
\caption{Image by $J_{10}$}\label{Im-J10Ehz}
\end{figure}

\paragraph{\bf Location of the zeros of Hecke type Faber Polynomial}
We can observe that some zeros of $F_{m, 10}$ do not lie on the lower arcs of $\partial \mathbb{F}_{10}$ for small weight $m$ by numerical calculation. However, when the weight $m$ increases, then the location of the zeros seems to approach to lower arcs of $\partial \mathbb{F}_{10}$. (see Figure \ref{Im-J10Ehz})\\

\clearpage

\section{Level $11$}

We have $\Gamma_0(11)+=\Gamma_0^{*}(11)$ and $\Gamma_0(11)-=\Gamma_0(11)$, but $\Gamma_0(11)$ is of genus $1$. We have $W_{11} = \left(\begin{smallmatrix}0&-1 / \sqrt{11}\\ \sqrt{11}&0\end{smallmatrix}\right)$.\\

\subsection{$\Gamma_0^{*}(11)$}

We have $\Gamma_0^{*}(11) = \langle \left( \begin{smallmatrix} 1 & 1 \\ 0 & 1 \end{smallmatrix} \right), \: W_{11}, \: \left( \begin{smallmatrix} 4 & 1 \\ 11 & 3 \end{smallmatrix} \right), \: \left( \begin{smallmatrix} 3 & 1 \\ 11 & 4 \end{smallmatrix} \right) \rangle$.\\

\paragraph{\bf Location of the zeros of the Eisenstein series}
\begin{figure}[hbtp]
\begin{center}
\includegraphics[width=1.5in]{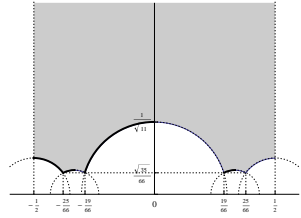}
\end{center}
\caption{Neighborhood of location of the zeros of $E_{k, 11+}^{\infty}$}
\end{figure}

\begin{figure}[hbtp]
\begin{center}
{{\small Lower arcs of $\partial \mathbb{F}_{11+}$}\includegraphics[width=2.5in]{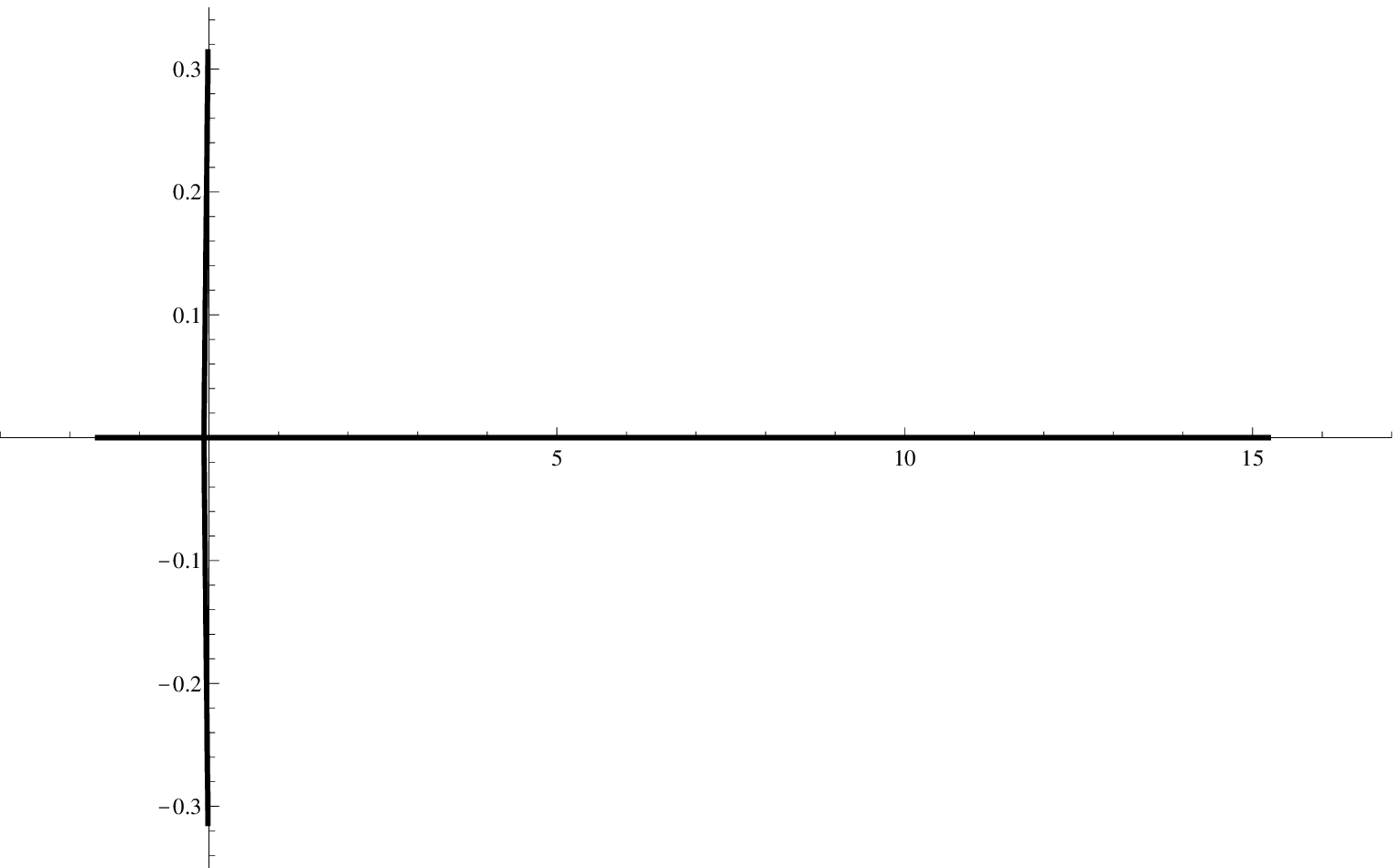}}
\end{center}
\caption{Image by $J_{11+}$}\label{Im-J11A}
\end{figure}

\begin{figure}[hbtp]
\begin{center}
\includegraphics[width=6.3in]{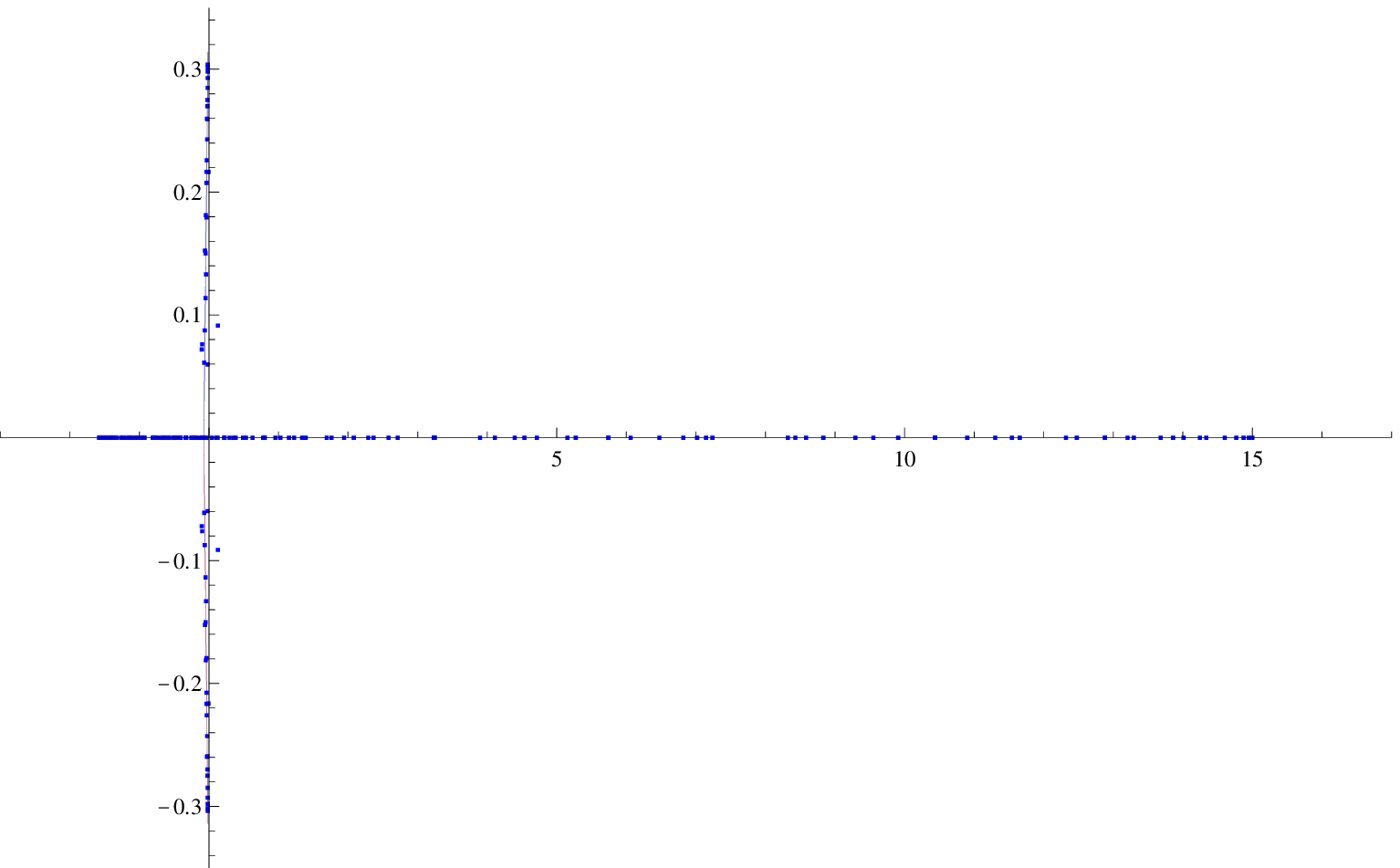}\\
The zeros of $E_{k, 11+}^{\infty}$ for $4 \leqslant k \leqslant 40$\\
\includegraphics[width=6.3in]{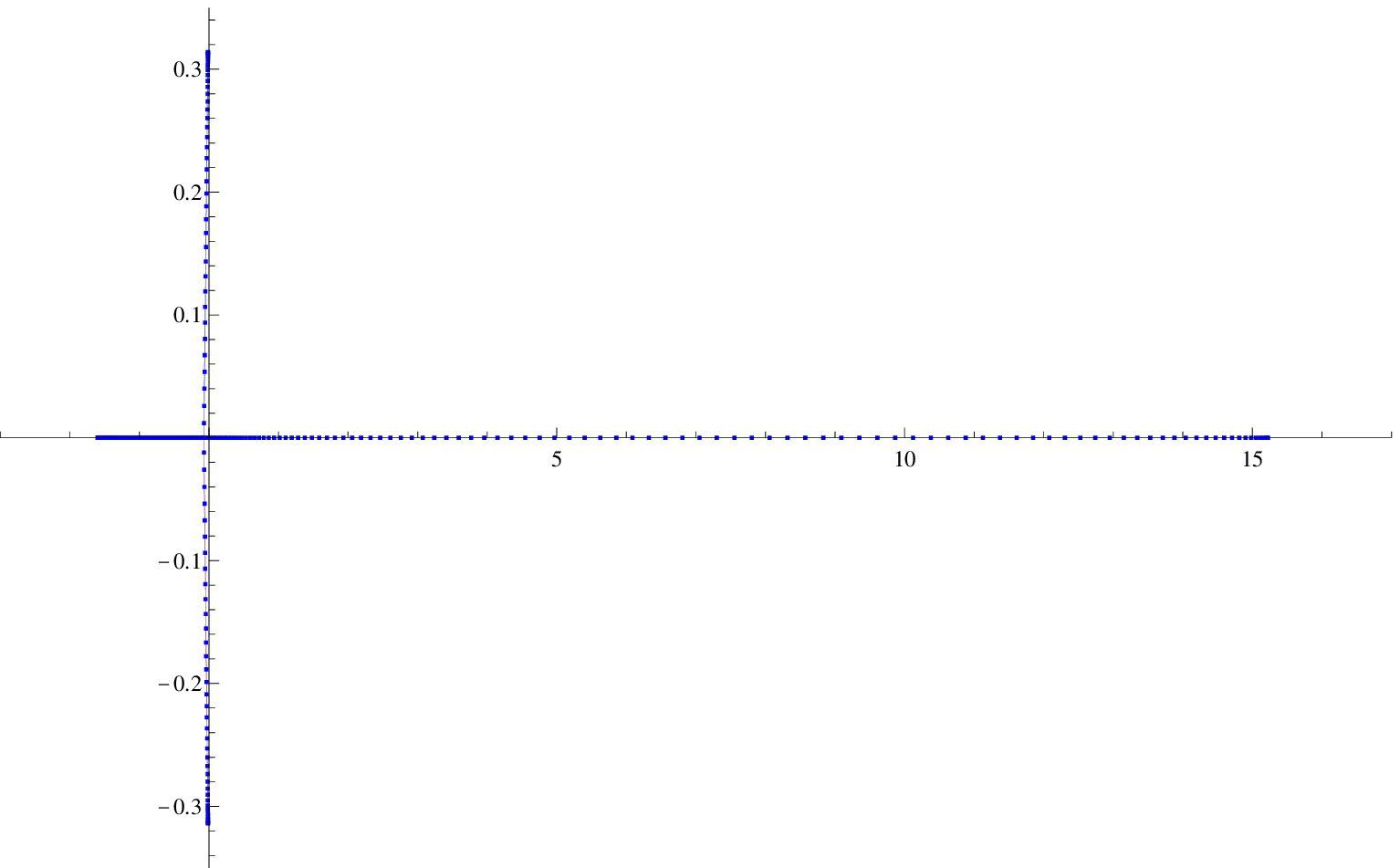}\\
The zeros of $E_{1000, 11+}^{\infty}$
\end{center}
\caption{Image by $J_{11+}$}\label{Im-J11Az}
\end{figure}

We can observe that the zeros of $E_{k, 11+}^{\infty}$ do not lie on the arcs of $\partial \mathbb{F}_{11+}$ for small weight $k$ by numerical calculation. However, when the weight $k$ increases, then the location of the zeros seems to approach to lower arcs of $\partial \mathbb{F}_{11+}$. (See Figure \ref{Im-J11Az})\\

\begin{figure}[hbtp]
\begin{center}
\includegraphics[width=6.3in]{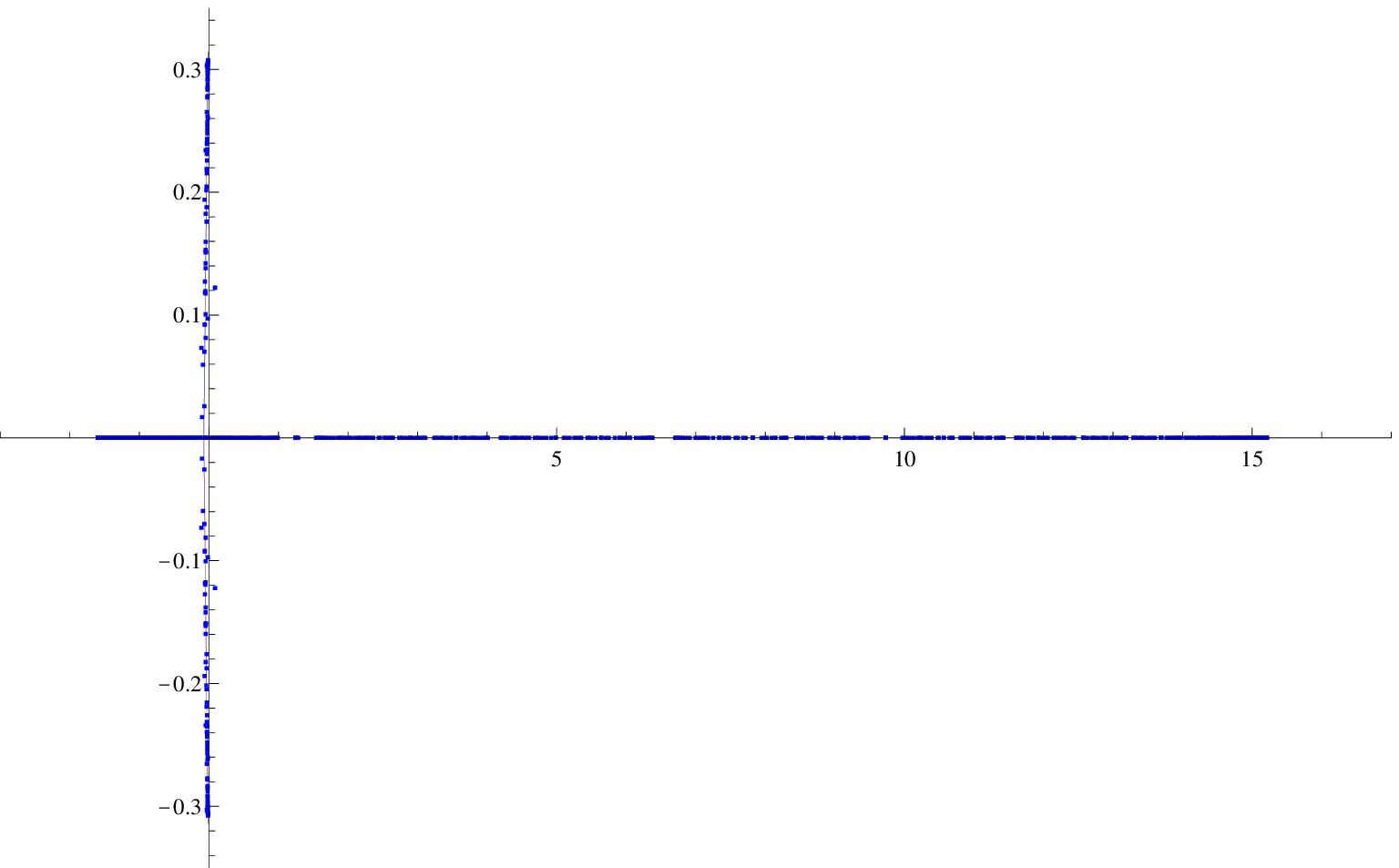}\\
The zeros of $F_{m, 11+}$ for $m \leqslant 40$\\
\includegraphics[width=6.3in]{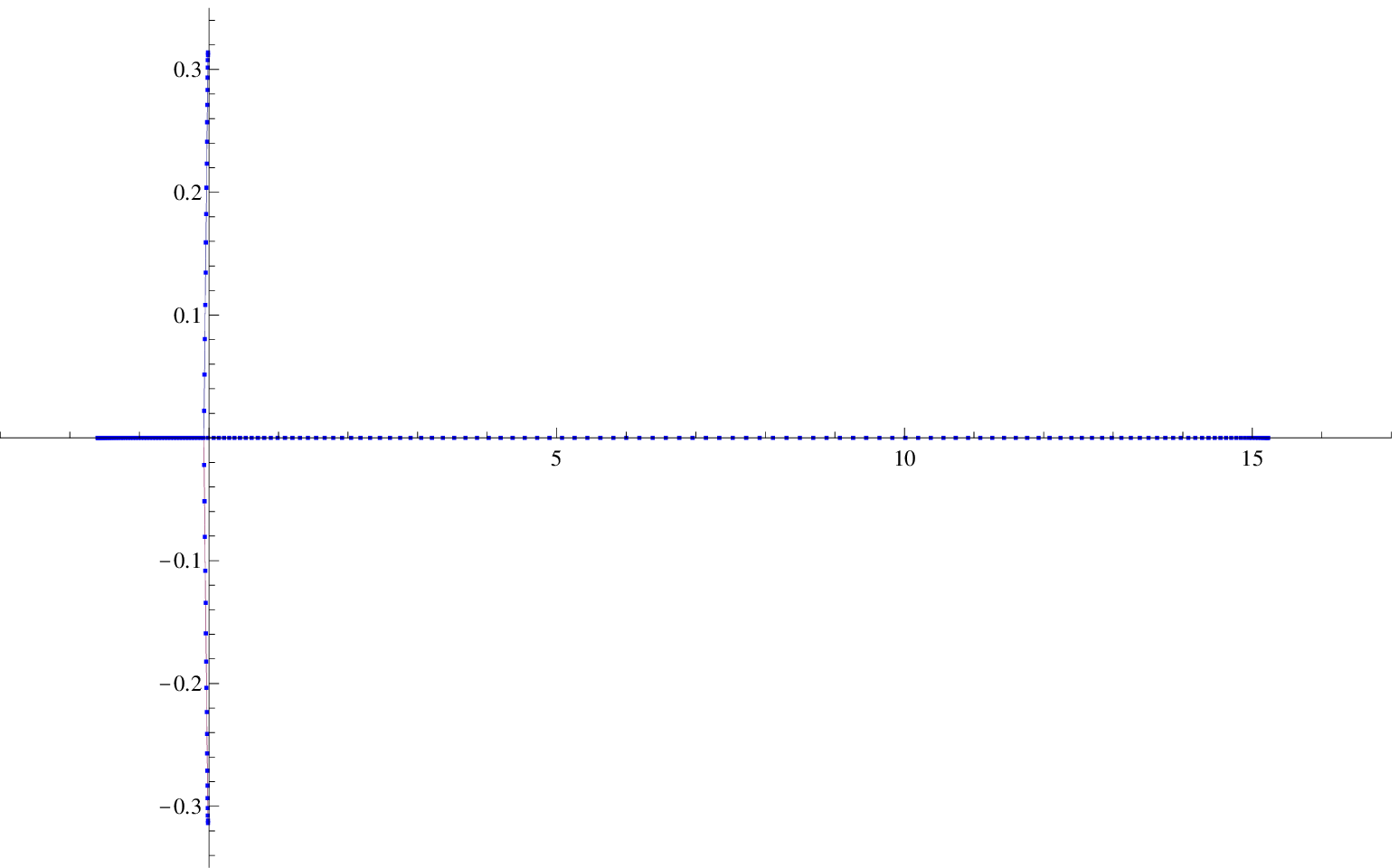}\\
The zeros of $F_{200, 11+}$
\end{center}
\caption{Image by $J_{11+}$}\label{Im-J11Ahz}
\end{figure}

\paragraph{\bf Location of the zeros of Hecke type Faber Polynomial}
We can observe that some zeros of $F_{m, 11+}$ do not lie on the lower arcs of $\partial \mathbb{F}_{11+}$ for small weight $m$ by numerical calculation. However, when the weight $m$ increases, then the location of the zeros seems to approach to lower arcs of $\partial \mathbb{F}_{11+}$. (see Figure \ref{Im-J11Ahz})\\
 \clearpage

\section{Level $12$}

We have $\Gamma_0(12)+$, $\Gamma_0(12)+12=\Gamma_0^{*}(12)$, $\Gamma_0(12)+4$, $\Gamma_0(12)+3$, and $\Gamma_0(12)-=\Gamma_0(12)$. We have $W_{12} = \left(\begin{smallmatrix}0 & - 1 / (2 \sqrt{3})\\ 2 \sqrt{3} & 0\end{smallmatrix}\right)$, $W_{12, 3} := \left(\begin{smallmatrix} - \sqrt{3} & -1 / \sqrt{3} \\ 4 \sqrt{3} & \sqrt{3} \end{smallmatrix}\right)$, $W_{12, 4} := \left(\begin{smallmatrix} -2 & 1/2 \\ 6 & -2 \end{smallmatrix}\right)$, $W_{12-, 2} := \left(\begin{smallmatrix} -1 & 0 \\ 6 & -1 \end{smallmatrix}\right)$, $W_{12+, 2} := \left(\begin{smallmatrix} -1 / \sqrt{2} & -1 / (2 \sqrt{2}) \\ 3 \sqrt{2} & 1 / \sqrt{2} \end{smallmatrix}\right)$, $W_{12-, 6} := \left(\begin{smallmatrix} - \sqrt{3} & -2 / \sqrt{3} \\ 2 \sqrt{3} & \sqrt{3} \end{smallmatrix}\right)$, and $W_{12+, 6} := \left(\begin{smallmatrix} - \sqrt{6} / 2 & -5 / (2 \sqrt{6}) \\ \sqrt{6} & \sqrt{6} / 2 \end{smallmatrix}\right)$.\\

\subsection{$\Gamma_0(12)+$}

We have $\Gamma_0(12)+ = T_{1/2}^{-1} (\Gamma_0(6)+3) T_{1/2}$ and $\Gamma_0(12)+ = \langle \left( \begin{smallmatrix} 1 & 1 \\ 0 & 1 \end{smallmatrix} \right), \: W_{12}, \: W_{12, 4} \rangle$. Furthermore, we have $\gamma_{-1/2} = W_{12+, 6}$.\\

\paragraph{\bf Location of the zeros of the Eisenstein series}
Since $W_{12+, 6}^{- 1} (\Gamma_0(12)+) W_{12+, 6} = \Gamma_0(12)+$, we have
\begin{equation}
E_{k, 12+}^{-1/2}(W_{12+, 6} z) = (\sqrt{6} z + \sqrt{6} / 2)^k E_{k, 12+}^{\infty}(z).
\end{equation}
Furthermore, we have
\begin{align*}
E_{k, 12+}^{-1/2} (e^{i \theta'} / (2 \sqrt{3})) &= ((\sqrt{3} e^{i \theta} -1) / \sqrt{2})^k E_{k, 12+}^{\infty}(e^{i \theta} / (2 \sqrt{3})),\\
E_{k, 12+}^{-1/2} (i / (2 \tan(\theta/2))) &= ((e^{i \theta} + 1) / \sqrt{6})^k E_{k, 12+}^{\infty}(e^{i \theta} / 6 - 1/3),
\end{align*}
where $e^{i \theta'} = (- \sqrt{3} - 2 \cos\theta + i \sin\theta) / (2 + \sqrt{3} \cos\theta)$.

Now, recall that $\Gamma_0(12)+ = T_{1/2}^{-1} \: (\Gamma_0(6)+3) \: T_{1/2}$. Then, for $k \leqslant 600$, since we can prove that all of the zeros of $E_{k, 6+3}^{\infty}$ lie on the lower arcs of $\partial \mathbb{F}_{6+3}$ by numerical calculation, we have all of the zeros of $E_{k, 12+}^{\infty}$ in the lower arcs of $\partial \mathbb{F}_{12+}$.

\begin{figure}[hbtp]
\begin{center}
{{$E_{k, 12+}^{\infty}$}\includegraphics[width=1.5in]{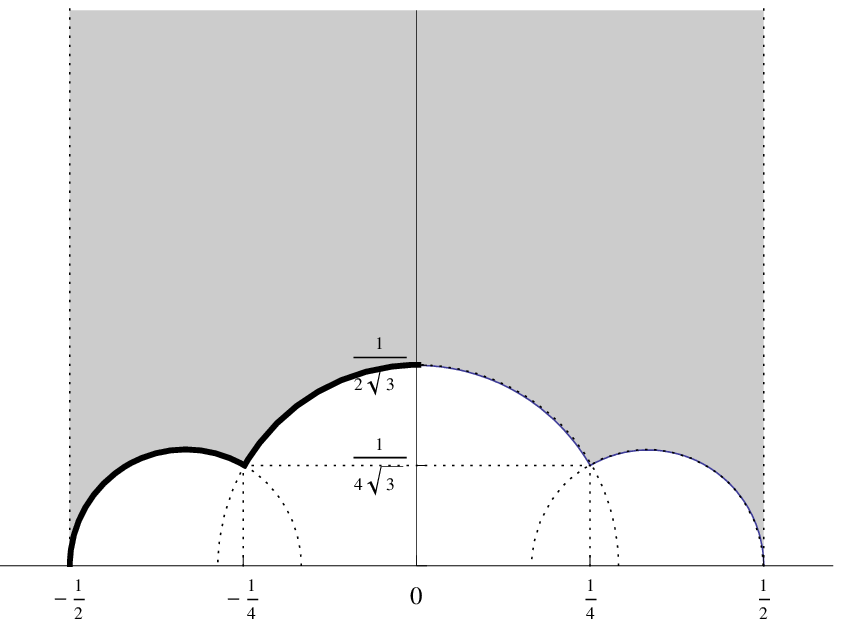}}
\quad \quad
{{$E_{k, 12+}^{-1/2}$}\includegraphics[width=1.5in]{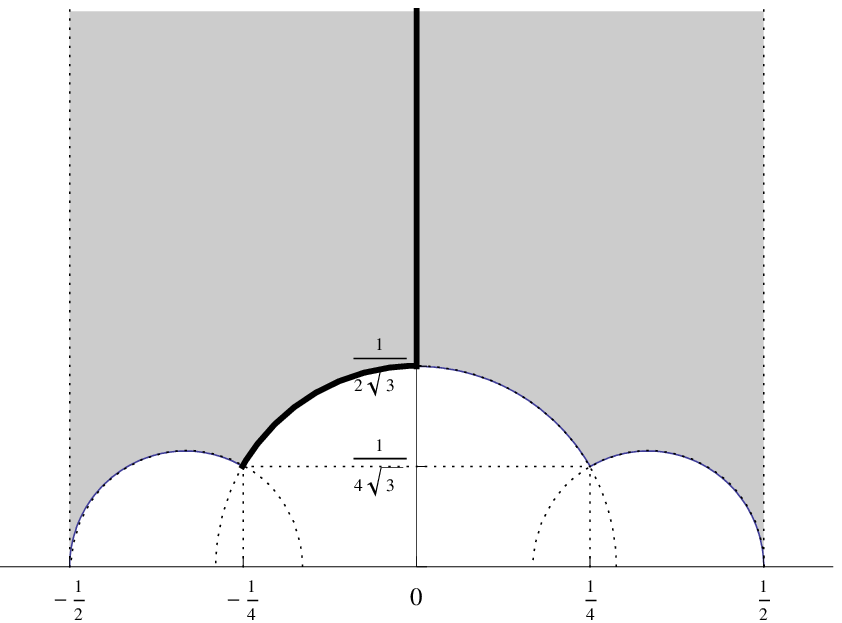}}
\end{center}
\caption{Location of the zeros of the Eisenstein series}
\end{figure}

\paragraph{\bf Location of the zeros of Hecke type Faber Polynomial}
Similarly to the Eisenstein series, for $m \leqslant 200$, since we can prove that all of the zeros of $F_{m, 6+3}$ lie on the lower arcs of $\partial \mathbb{F}_{6+3}$ by numerical calculation, we have all of the zeros of $F_{m, 12+}$ in the lower arcs of $\partial \mathbb{F}_{12+}$.\\

\subsection{$\Gamma_0(12)+12 = \Gamma_0^{*}(12)$}

We have $\Gamma_0^{*}(12) = \langle \left( \begin{smallmatrix} 1 & 1 \\ 0 & 1 \end{smallmatrix} \right), \: W_{12}, \: \left( \begin{smallmatrix} 5 & 2 \\ 12 & 5 \end{smallmatrix} \right), \: \left( \begin{smallmatrix} 7 & 2 \\ 24 & 7 \end{smallmatrix} \right) \rangle$, $\gamma_{-1/3} = W_{12, 4}$, and $\gamma_{-1/2} = W_{12-, 6}$.

\paragraph{\bf Location of the zeros of the Eisenstein series}
Since $W_{12, 4}^{- 1} \Gamma_0^{*}(12) W_{12, 4} = \Gamma_0^{*}(12)$, we have
\begin{equation}
E_{k, 12+12}^{-1/3}(W_{12, 4} z) = (6 z - 2)^k E_{k, 12+12}^{\infty}(z).
\end{equation}
Furthermore, we have
\begin{align*}
E_{k, 12+12}^{-1/3} (-1/2 + i / (6 \tan(\theta/2))) &= (- (e^{i \theta} - 1) / 2)^k E_{k, 12+12}^{\infty}((e^{i \theta} - 5)/12),\\
E_{k, 12+12}^{-1/3} (i \tan(\theta/2) / 3) &= (- (e^{i \theta} + 1) / 4)^k E_{k, 12+12}^{\infty}((e^{i \theta} - 7)/24),\\
E_{k, 12+12}^{-1/3} (e^{i \theta'} / (2 \sqrt{3})) &= (- (\sqrt{3} e^{i \theta} + 2))^k E_{k, 12+12}^{\infty}(e^{i \theta} / (2 \sqrt{3})),
\end{align*}
where $e^{i \theta'} = (4 \sqrt{3} + 7 \cos\theta + i \sin\theta) / (7 + 4 \sqrt{3} \cos\theta)$. On the other hand, recall that $E_{k, 12+12}^{-1/2}(z) = 2^{-1} 2^{-k/2} E_{k, 12+}^{\infty}(z)$. Moreover, by the transformation with $W_{12, 3}$ for $E_{k, 12+}^{-1/2}$, we have
\begin{align*}
E_{k, 12+12}^{-1/2} (e^{i (\pi - \theta')} / (2 \sqrt{3})) &= (2 e^{i \theta} + \sqrt{3})^k E_{k, 12+12}^{-1/2}(e^{i \theta} / (2 \sqrt{3})),\\
E_{k, 12+12}^{-1/2} ((e^{i \theta} - 7)/24) &= (\sqrt{3} (1 + i / \tan(\theta/2)))^k E_{k, 12+12}^{-1/2}(i / (4 \tan(\theta/2))).
\end{align*}

For $k \leqslant 500$, we can prove that all of the zeros of $E_{k, 12+12}^{\infty}$ lie on the lower arcs of $\partial \mathbb{F}_{12+12}$ by numerical calculation.

\begin{figure}[htbp]
\begin{center}
{{$E_{k, 12+12}^{\infty}$}\includegraphics[width=1.5in]{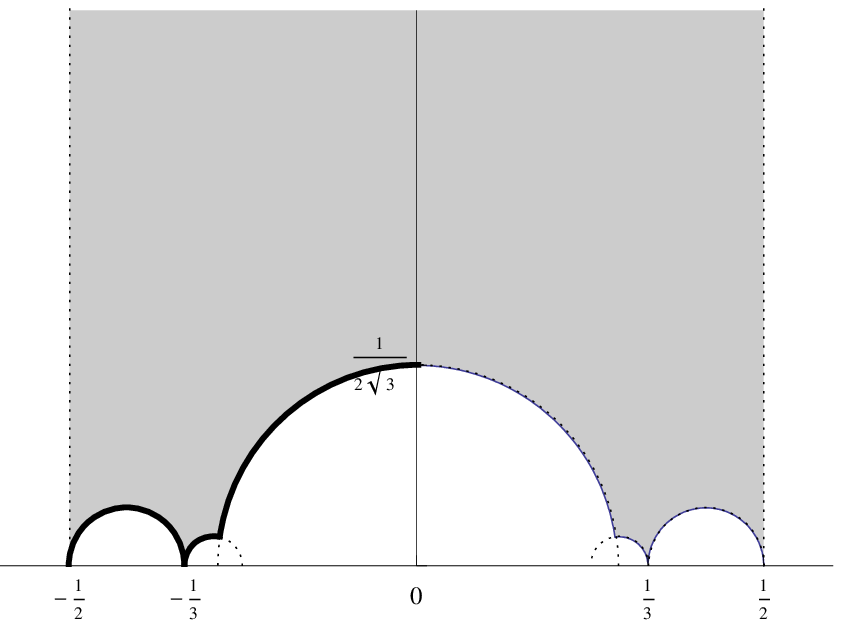}}
\;
{{$E_{k, 12+12}^{-1/3}$}\includegraphics[width=1.5in]{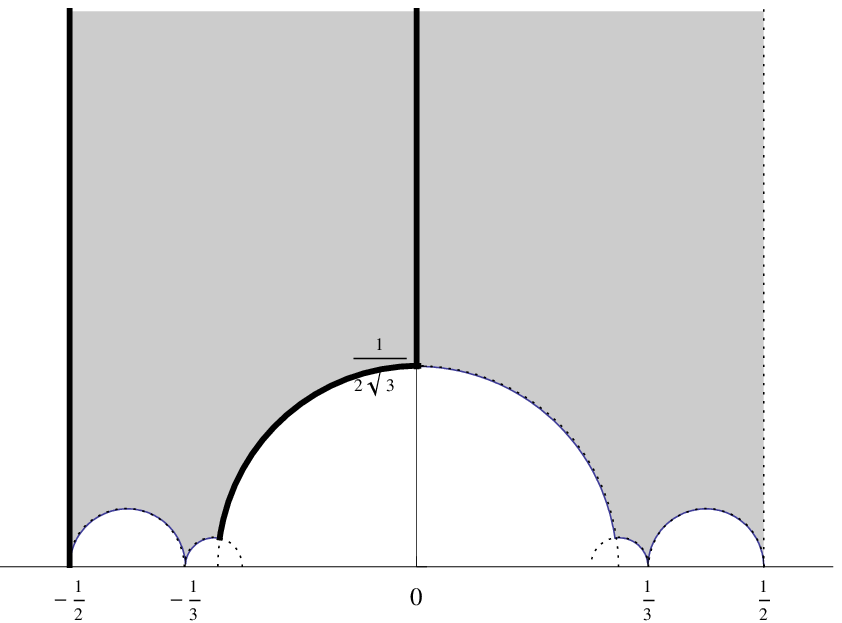}}
\;
{{$E_{k, 12+12}^{-1/2}$}\includegraphics[width=1.5in]{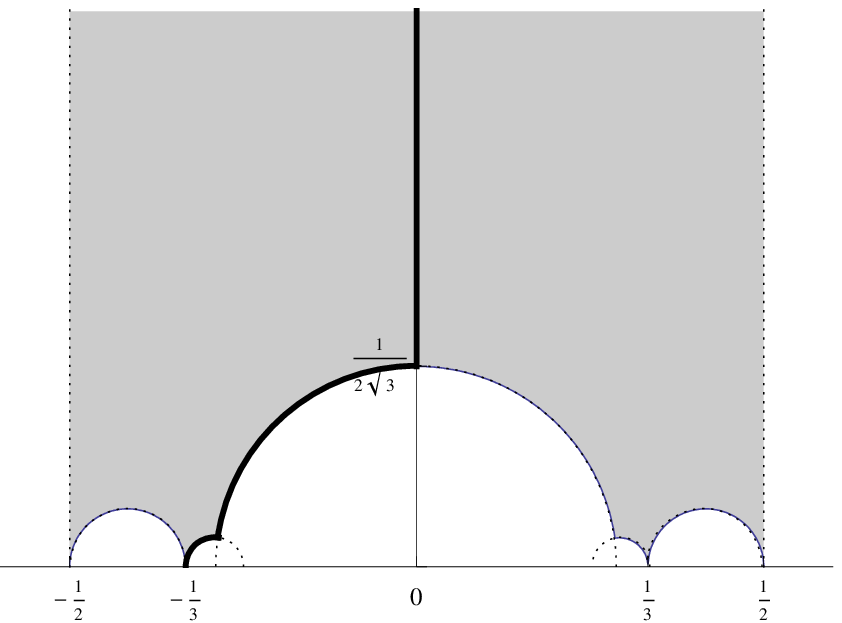}}
\end{center}
\caption{Location of the zeros of the Eisenstein series}
\end{figure}

\paragraph{\bf Location of the zeros of Hecke type Faber Polynomial}
For every integer $m \leqslant 200$ such that $m \not\equiv 2, 4 \pmod{6}$ but $m = 3, 6, 12, 13, 15$, we can prove that all of the zeros of $F_{m, 12+12}$ lie on the lower arcs of $\partial \mathbb{F}_{12+12}$ by numerical calculation. On the other hand, by numerical calculation, for $m = 3, 6, 12, 13, 15$, we can prove that all but two of the zeros of $F_{m, 12+12}$ lie on the lower arcs of $\partial \mathbb{F}_{12+12}$, and two of the zeros of $F_{m, 12+12}$ do not lie on $\partial \mathbb{F}_{12+12}$. For the other cases where $m$ is $m \leqslant 200$ such that $m \equiv 2, 4 \pmod{6}$, by numerical calculation, we can prove that all but one of the zeros of $F_{m, 12+12}$ lie on the lower arcs of $\partial \mathbb{F}_{12+12}$, and one of the zeros of $F_{m, 12+12}$ lies on $\partial \mathbb{F}_{12+12}$ but does not on the lower arcs.\\

\subsection{$\Gamma_0(12)+4$}

We have $\Gamma_0(12)+4 = T_{1/2}^{-1} \Gamma_0(6) T_{1/2}$ and $\Gamma_0(12)+4 = \langle - I, \: \left( \begin{smallmatrix} 1 & 1 \\ 0 & 1 \end{smallmatrix} \right), \: W_{12, 4}, \: \left( \begin{smallmatrix} 1 & 0 \\ 12 & 1 \end{smallmatrix} \right) \rangle$. Furthermore, we have $\gamma_0 = W_{12}$, $\gamma_{-1/2} = W_{12+, 6}$, and $\gamma_{-1/6} = W_{12+, 2}$.\\

\paragraph{\bf Location of the zeros of the Eisenstein series}
Since $W_{12}^{- 1} (\Gamma_0(12)+4) W_{12} = W_{12+, 6}^{- 1} (\Gamma_0(12)+4) W_{12+, 6} = W_{12+, 2}^{- 1} (\Gamma_0(12)+4) W_{12+, 2} = \Gamma_0(12)+4$, we have
\begin{align*}
(2 \sqrt{3} z)^{-k} E_{k, 12+4}^0(W_{12} z) = (\sqrt{6} z + \sqrt{6} / 2)^{-k} E_{k, 12+4}^{-1/2}(W_{12+, 6} z)
 = (3 \sqrt{2} z + 1 / \sqrt{2})^{-k} E_{k, 12+4}^{-1/6}(W_{12+, 2} z) = E_{k, 12+4}^{\infty}(z).
\end{align*}
Furthermore, we have
\begin{align*}
E_{k, 12+4}^0 (-1/2 + i / (2 \tan(\theta/2))) &= ((e^{i \theta} -1) / (2 \sqrt{3}))^k E_{k, 12+4}^{\infty}((e^{i \theta} - 1) / 12),\\
E_{k, 12+4}^0 (-1/3 + e^{i \theta'} / 6) &= (- (2 e^{i \theta} - 1) / \sqrt{3})^k E_{k, 12+4}^{\infty}(-1/3 + e^{i \theta} / 6),\\
E_{k, 12+4}^{-1/2} ((e^{i \theta''} - 1) / 12) &= ((5 e^{i \theta} + 1) / (2 \sqrt{6}))^k E_{k, 12+4}^{\infty}((e^{i \theta} - 1) / 12),\\
E_{k, 12+4}^{-1/2} (i \tan(\theta/2) / 2) &= ((e^{i \theta} + 1) / \sqrt{6})^k E_{k, 12+4}^{\infty}(-1/3 + e^{i \theta} / 6),\\
E_{k, 12+4}^{-1/6} (-1/2 + i \tan(\theta/2) / 3) &= ((e^{i \theta} + 1) / (2 \sqrt{2}))^k E_{k, 12+4}^{\infty}((e^{i \theta} - 1) / 12),\\
E_{k, 12+4}^{-1/6} (i / (6 \tan(\theta/2))) &= ((e^{i \theta} - 1) / \sqrt{2})^k E_{k, 12+4}^{\infty}(-1/3 + e^{i \theta} / 6),
\end{align*}
where $e^{i \theta'} = (4 - 5 \cos\theta + 3 i \sin\theta) / (5 - 4 \cos\theta)$ and $e^{i \theta''} = (- 5 - 13 \cos\theta + 12 i \sin\theta) / (13 + 5 \cos\theta)$.

Now, recall that $\Gamma_0(12)+4 = T_{1/2}^{-1} \: \Gamma_0(6) \: T_{1/2}$. Then, for $k \leqslant 500$, since we can prove that all of the zeros of $E_{k, 6}^{\infty}$ lie on the lower arcs of $\partial \mathbb{F}_{6}$ by numerical calculation, we have all of the zeros of $E_{k, 12+4}^{\infty}$ in the lower arcs of $\partial \mathbb{F}_{12+4}$.\\

\begin{figure}[hbtp]
\begin{center}
{{$E_{k, 12+4}^{\infty}$}\includegraphics[width=1.5in]{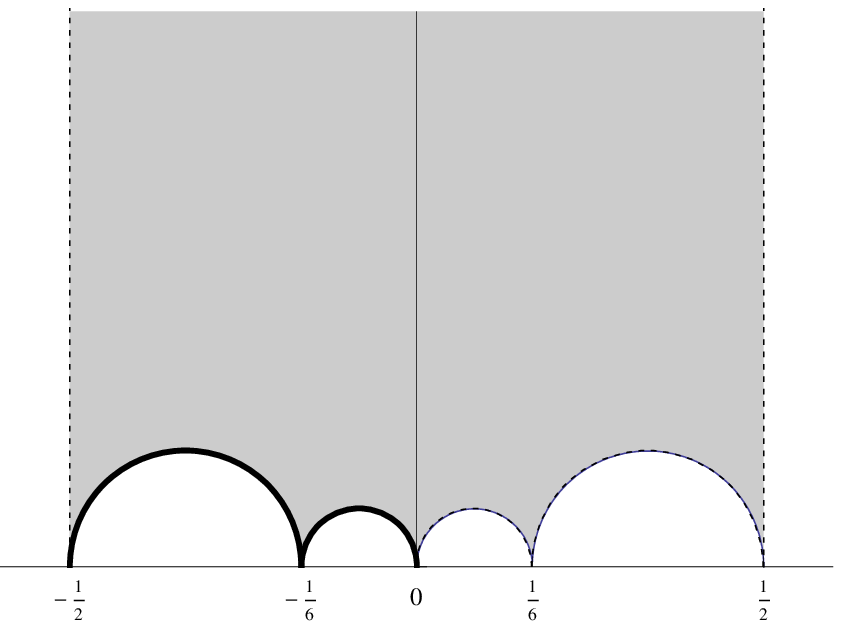}}
\quad \quad
{{$E_{k, 12+4}^0$}\includegraphics[width=1.5in]{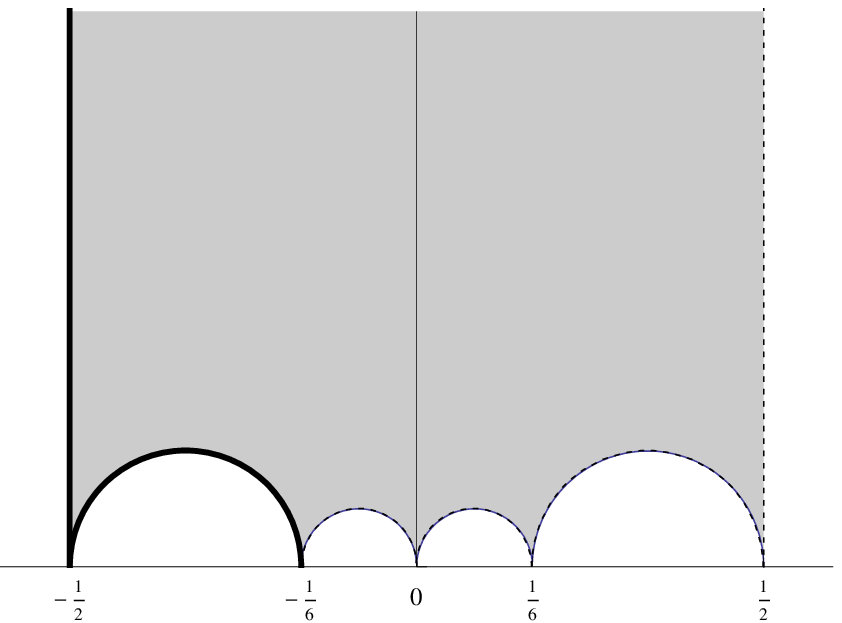}}\\

{{$E_{k, 12+4}^{-1/2}$}\includegraphics[width=1.5in]{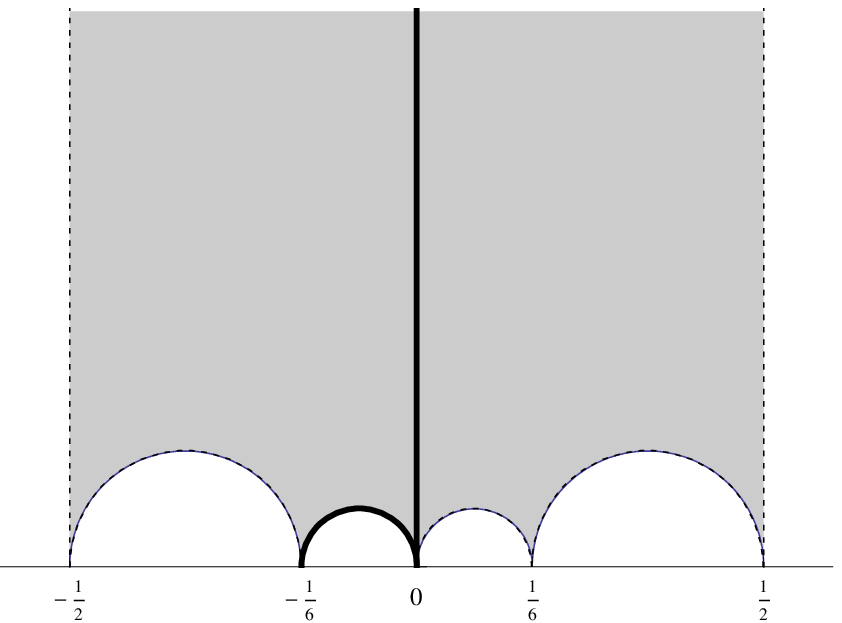}}
\quad \quad
{{$E_{k, 12+4}^{-1/6}$}\includegraphics[width=1.5in]{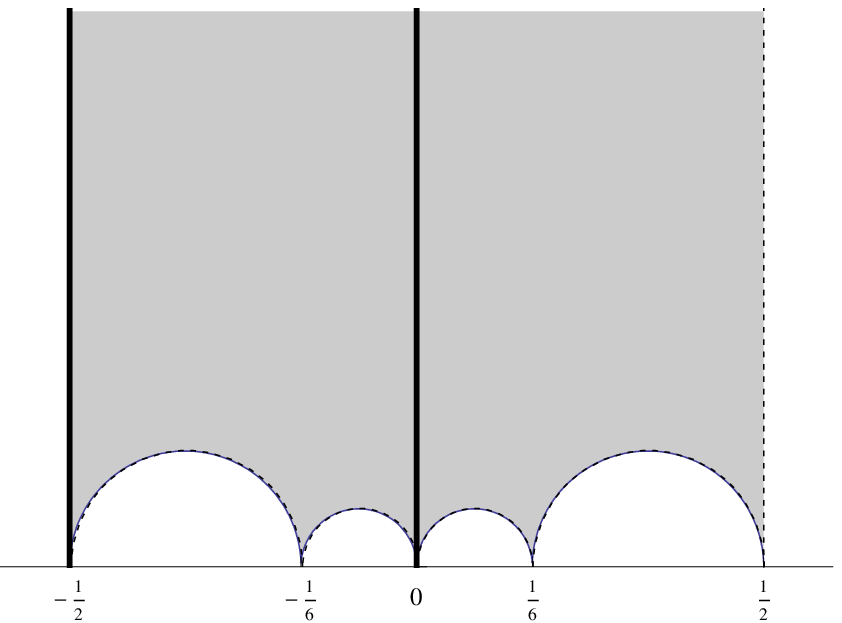}}
\end{center}
\caption{Location of the zeros of the Eisenstein series}
\end{figure}

\paragraph{\bf Location of the zeros of Hecke type Faber Polynomial}
Similarly to the Eisenstein series, for $m \leqslant 200$, since we can prove that all of the zeros of $F_{m, 6}$ lie on the lower arcs of $\partial \mathbb{F}_{6}$ by numerical calculation, we have all of the zeros of $F_{m, 12+4}$ in the lower arcs of $\partial \mathbb{F}_{12+4}$.

\newpage

\subsection{$\Gamma_0(12)+3$}

We have $\Gamma_0(12)+3 = \langle \left( \begin{smallmatrix} 1 & 1 \\ 0 & 1 \end{smallmatrix} \right), \: W_{12, 3}, \: \left( \begin{smallmatrix} 1 & 0 \\ 12 & 1 \end{smallmatrix} \right), \: \left( \begin{smallmatrix} 5 & 2 \\ 12 & 5 \end{smallmatrix} \right) \rangle$, .$\gamma_0 = W_{12}$, and $\gamma_{-1/2} = W_{12-, 6}$.\\

\paragraph{\bf Location of the zeros of the Eisenstein series}
Since $W_{12}^{- 1} (\Gamma_0(12)+3) W_{12} = W_{12-, 6}^{- 1} (\Gamma_0(12)+3) W_{12-, 6} = \Gamma_0(12)+3$, we have
\begin{equation}
(2 \sqrt{3} z)^{-k} E_{k, 12+3}^0(W_{12} z) = (2 \sqrt{3} z + \sqrt{3})^{-k} E_{k, 12+3}^{-1/2}(W_{12-, 6} z) = E_{k, 12+3}^{\infty}(z).
\end{equation}
Furthermore, we have
\begin{align*}
E_{k, 12+3}^0 (-1/2 + i / (2 \tan(\theta/2))) &= ((e^{i \theta} - 1) / (2 \sqrt{3}))^k E_{k, 12+3}^{\infty}((e^{i \theta} - 1)/12),\\
E_{k, 12+3}^0 (-1/2 + i / (6 \tan(\theta/2))) &= (- (e^{i \theta} - 1) / 2)^k E_{k, 12+3}^{\infty}((e^{i \theta} - 5)/12),\\
E_{k, 12+3}^0 (1/2 + e^{i \theta'} / (2 \sqrt{3})) &= (- (\sqrt{3} e^{i \theta} + 1) / 2)^k E_{k, 12+3}^{\infty}(-1/4 + e^{i \theta} / (4 \sqrt{3})),\\
E_{k, 12+3}^0 (-1/2 + e^{i (\pi - \theta')} / (2 \sqrt{3})) &= ((e^{i \theta} + \sqrt{3}) / 2)^k E_{k, 12+3}^{\infty}(1/4 + e^{i \theta} / (4 \sqrt{3})),\\
E_{k, 12+3}^{-1/2} (i \tan(\theta/2) / 6) &= ((e^{i \theta} + 1) / 2)^k E_{k, 12+3}^{\infty}((e^{i \theta} - 1)/12),\\
E_{k, 12+3}^{-1/2} (i \tan(\theta/2) / 2) &= ((e^{i \theta} + 1) / (2 \sqrt{3}))^k E_{k, 12+3}^{\infty}((e^{i \theta} - 5)/12),\\
E_{k, 12+3}^{-1/2} (e^{i (\pi - \theta')} / (2 \sqrt{3})) &= ((e^{i \theta} + \sqrt{3}) / 2)^k E_{k, 12+3}^{\infty}(-1/4 + e^{i \theta} / (4 \sqrt{3})),\\
E_{k, 12+3}^{-1/2} (e^{i \theta'} / (2 \sqrt{3})) &= (- (\sqrt{3} e^{i \theta} + 1) / 2)^k E_{k, 12+3}^{\infty}(1/4 + e^{i \theta} / (4 \sqrt{3})),
\end{align*}
where $e^{i \theta'} = (- \sqrt{3} - 2 \cos\theta + i \sin\theta) / (2 + \sqrt{3} \cos\theta)$.

\begin{figure}[htbp]
\begin{center}
{{$E_{k, 12+3}^{\infty}$}\includegraphics[width=1.5in]{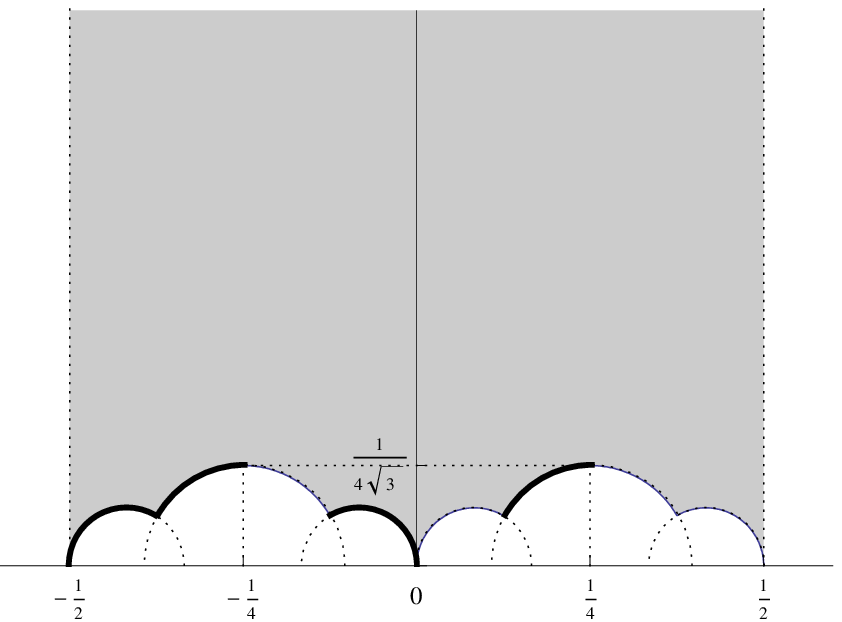}}
\;
{{$E_{k, 12+3}^0$}\includegraphics[width=1.5in]{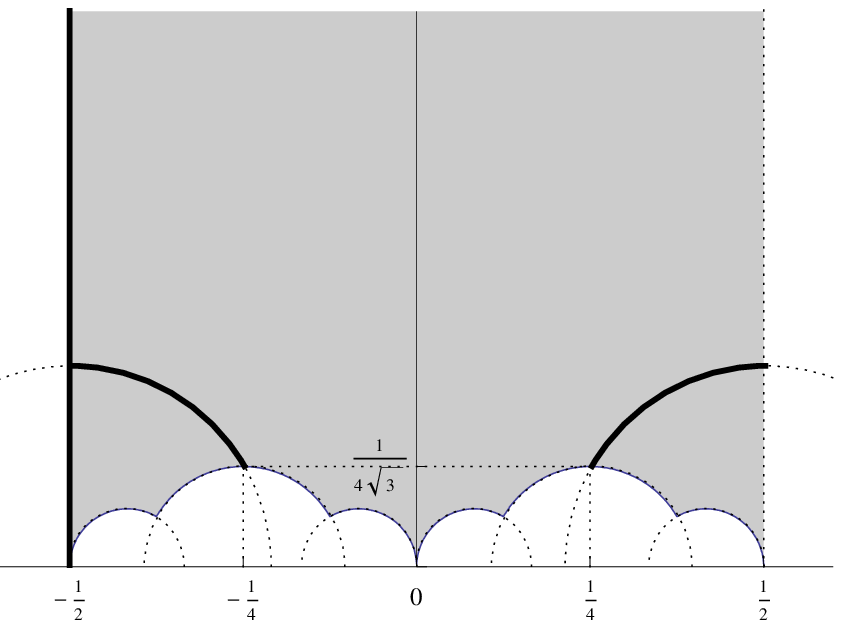}}
\;
{{$E_{k, 12+3}^{-1/2}$}\includegraphics[width=1.5in]{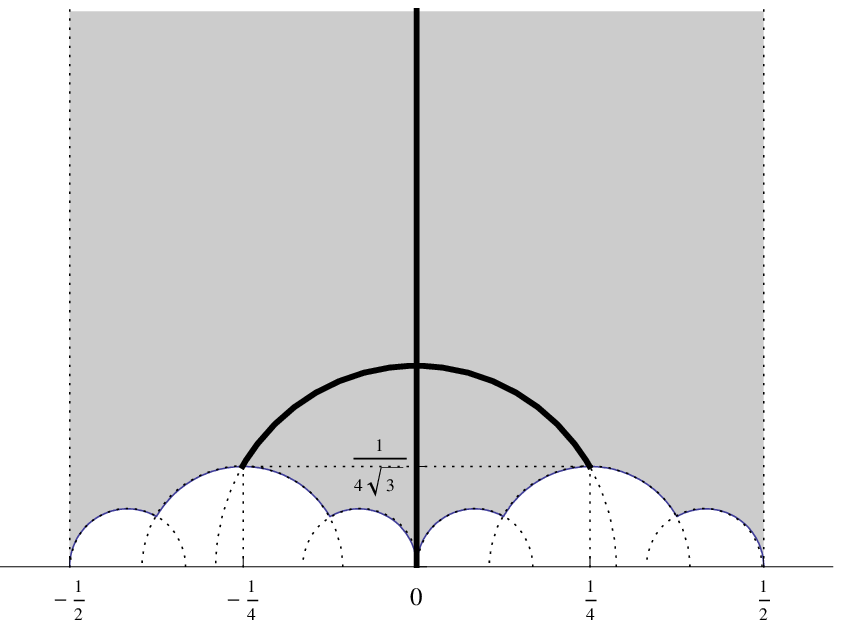}}
\end{center}
\caption{Location of the zeros of the Eisenstein series}
\end{figure}

\begin{figure}[hbtp]
\begin{center}
{{\small Lower arcs of $\partial \mathbb{F}_{12+3}$}\includegraphics[width=2.5in]{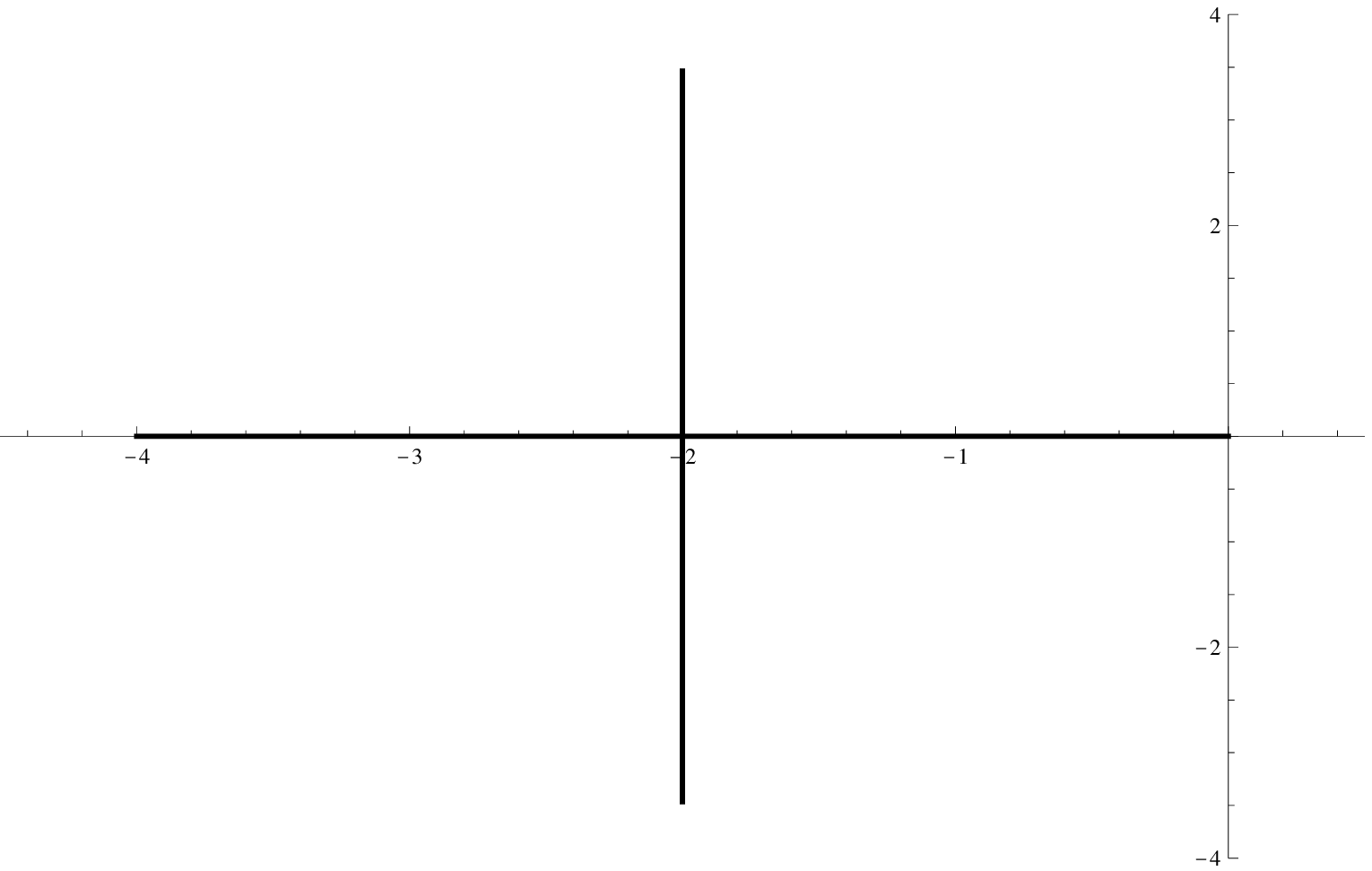}}
\end{center}
\caption{Image by $J_{12+3}$}\label{Im-J12E}
\end{figure}

Now, recall that $E_{k, 12+3}^{\infty}(z) = E_{k, 6+3}^{\infty}(2 z)$. For $k \leqslant 600$, we can prove that all of the zeros of $E_{k, 6+3}^{\infty}$ lie on the lower arcs of $\partial \mathbb{F}_{6+3}$ by numerical calculation, then we have all of the zeros of $E_{k, 12+3}^{\infty}$ in the lower arcs of $\partial \mathbb{F}_{12+3}$. Similarly to $\Gamma_0(9)$, this case is interesting.\\

\paragraph{\bf Location of the zeros of Hecke type Faber Polynomial}
For $m \leqslant 200$, we can prove that all of the zeros of $F_{m, 12+3}$ lie on the lower arcs of $\partial \mathbb{F}_{12+3}$ by numerical calculation.

\newpage

\subsection{$\Gamma_0(12)$}

We have $\Gamma_0(12) = \langle - I, \: \left( \begin{smallmatrix} 1 & 1 \\ 0 & 1 \end{smallmatrix} \right), \: \left( \begin{smallmatrix} 1 & 0 \\ 12 & 1 \end{smallmatrix} \right), \: \left( \begin{smallmatrix} 5 & 2 \\ 12 & 5 \end{smallmatrix} \right), \: \left( \begin{smallmatrix} 5 & 1 \\ 24 & 5 \end{smallmatrix} \right), \: \left( \begin{smallmatrix} 7 & 2 \\ 24 & 7 \end{smallmatrix} \right) \rangle$, $\gamma_0 = W_{12}$, $\gamma_{-1/3} = W_{12, 4}$, $\gamma_{-1/4} = W_{12, 3}$, $\gamma_{-1/2} = W_{12-, 6}$, and $\gamma_{-1/6} = W_{12-, 2}$.\\

\paragraph{\bf Location of the zeros of the Eisenstein series}
Since $W_{12}^{- 1} (\Gamma_0(12)) W_{12} = W_{12, 4}^{- 1} (\Gamma_0(12)) W_{12, 4} = W_{12, 3}^{- 1} (\Gamma_0(12)) W_{12, 3} = \Gamma_0(12)$, we have
\begin{equation*}
(2 \sqrt{3} z)^{-k} E_{k, 12}^0(W_{12} z) = (6 z - 2)^{-k} E_{k, 12}^{-1/2}(W_{6, 3} z)
 = (4 \sqrt{3} z + \sqrt{3})^{-k} E_{k, 12}^{-1/6}(W_{12, 4} z) = E_{k, 12}^{\infty}(z).
\end{equation*}
Furthermore, we have
\begin{align*}
E_{k, 12}^0 (-1/2 + i / (2 \tan(\theta/2))) &= ((e^{i \theta} -1) / (2 \sqrt{3}))^k E_{k, 12}^{\infty} ((e^{i \theta} - 1) / 12),\\
E_{k, 12}^0 ((e^{i \theta_1} - 5) / 12) &= (- (5 e^{i \theta} - 1) / (4 \sqrt{3}))^k E_{k, 12}^{\infty} ((e^{i \theta} - 5) / 24),\\
E_{k, 12}^0 ((e^{i \theta_2} - 7) / 24) &= (- (7 e^{i \theta} -1) / (4 \sqrt{3}))^k E_{k, 12}^{\infty} ((e^{i \theta} - 7) / 24),\\
E_{k, 12}^0 ((e^{i \theta_1} - 5) / 24) &= (- (5 e^{i \theta} - 1) / (2 \sqrt{3}))^k E_{k, 12}^{\infty} ((e^{i \theta} - 5) / 12),\\
E_{k, 12}^{-1/3} ((e^{i \theta_3} - 5) / 24) &= (- (3 e^{i \theta} + 1) / 2)^k E_{k, 12}^{\infty} ((e^{i \theta} - 1) / 12),\\
E_{k, 12}^{-1/3} ((e^{i \theta_3} - 1) / 12) &= ((3 e^{i \theta} + 1) / 4)^k E_{k, 12}^{\infty} ((e^{i \theta} - 5) / 24),\\
E_{k, 12}^{-1/3} (i \tan(\theta/2) / 3) &= (- (e^{i \theta} + 1) / 4)^k E_{k, 12}^{\infty} ((e^{i \theta} - 7) / 24),\\
E_{k, 12}^{-1/3} (-1/2 + i / (6 \tan(\theta/2))) &= (- (e^{i \theta} - 1) / 2)^k E_{k, 12}^{\infty} ((e^{i \theta} - 5) / 24),\\
E_{k, 12}^{-1/4} ((e^{i \theta_4} - 5) / 24) &= ((e^{i \theta} + 2) / \sqrt{3})^k E_{k, 12}^{\infty} ((e^{i \theta} - 1) / 12),\\
E_{k, 12}^{-1/4} (-1/2 + i \tan(\theta/2) / 4) &= ((e^{i \theta} + 1) / (2 \sqrt{3}))^k E_{k, 12}^{\infty} ((e^{i \theta} - 5) / 24),\\
E_{k, 12}^{-1/4} (i / (4 \tan(\theta/2))) &= ((e^{i \theta} - 1) / (2 \sqrt{3}))^k E_{k, 12}^{\infty} ((e^{i \theta} - 7) / 24),\\
E_{k, 12}^{-1/4} ((e^{i {\theta_4}'} - 1) / 12) &= ((e^{i \theta} - 2) / \sqrt{3})^k E_{k, 12}^{\infty} ((e^{i \theta} - 5) / 12),\\
E_{k, 12}^{-1/2} ((e^{i {\theta_1}'} - 7) / 24) &= ((5 e^{i \theta} + 1) / (2 \sqrt{3}))^k E_{k, 12}^{\infty} ((e^{i \theta} - 1) / 12),\\
E_{k, 12}^{-1/2} ((e^{i {\theta_2}'} - 5) / 24) &= ((7 e^{i \theta} + 1) / (4 \sqrt{3}))^k E_{k, 12}^{\infty} ((e^{i \theta} - 5) / 24),\\
E_{k, 12}^{-1/2} ((e^{i {\theta_1}'} - 1) / 12) &= ((5 e^{i \theta} + 1) / (4 \sqrt{3}))^k E_{k, 12}^{\infty} ((e^{i \theta} - 7) / 24),\\
E_{k, 12}^{-1/2} (i \tan(\theta/2) / 2) &= ((e^{i \theta} + 1) / (2 \sqrt{3}))^k E_{k, 12}^{\infty} ((e^{i \theta} - 5) / 12),\\
E_{k, 12}^{-1/6} (i \tan(\theta/2) / 6) &= (- (e^{i \theta} + 1) / 2)^k E_{k, 12}^{\infty} ((e^{i \theta} - 1) / 12),\\
E_{k, 12}^{-1/6} (-1/2 + i / (3 \tan(\theta/2))) &= (- (e^{i \theta} - 1) / 4)^k E_{k, 12}^{\infty} ((e^{i \theta} - 5) / 24),\\
E_{k, 12}^{-1/6} ((e^{i {\theta_3}'} - 5) / 12) &= (- (3 e^{i \theta} - 1) / 4)^k E_{k, 12}^{\infty} ((e^{i \theta} - 7) / 24),\\
E_{k, 12}^{-1/6} ((e^{i {\theta_3}'} - 7) / 24) &= (- (3 e^{i \theta} - 1) / 2)^k E_{k, 12}^{\infty} ((e^{i \theta} - 5) / 12),
\end{align*}
where $e^{i \theta_1} = (5 - 13 \cos\theta + 12 i \sin\theta) / (13 - 5 \cos\theta)$, $e^{i {\theta_1}'} = (- 5 - 13 \cos\theta + 12 i \sin\theta) / (13 + 5 \cos\theta)$, $e^{i \theta_2} = (7 - 25 \cos\theta + 24 i \sin\theta) / (25 - 7 \cos\theta)$, $e^{i {\theta_2}'} = (- 7 - 25 \cos\theta + 24 i \sin\theta) / (25 + 7 \cos\theta)$, $e^{i \theta_3} = (- 3 - 5 \cos\theta + 4 i \sin\theta) / (5 + 3 \cos\theta)$, $e^{i {\theta_3}'} = (3 - 5 \cos\theta + 4 i \sin\theta) / (5 - 3 \cos\theta)$, $e^{i \theta_4} = (4 + 5 \cos\theta + 3 i \sin\theta) / (5 + 4 \cos\theta)$, and $e^{i {\theta_4}'} = (- 4 + 5 \cos\theta + 3 i \sin\theta) / (5 - 4 \cos\theta)$.

Now, recall that $E_{k, 12}^{\infty}(z) = E_{k, 6}^{\infty}(2 z)$. Then, for $k \leqslant 500$, since we can prove that all of the zeros of $E_{k, 6}^{\infty}$ lie on the lower arcs of $\partial \mathbb{F}_{6}$ by numerical calculation, we have all of the zeros of $E_{k, 12}^{\infty}$ in the lower arcs of $\partial \mathbb{F}_{12}$.\\

\begin{figure}[hbtp]
\begin{center}
{{$E_{k, 12}^{\infty}$}\includegraphics[width=1.5in]{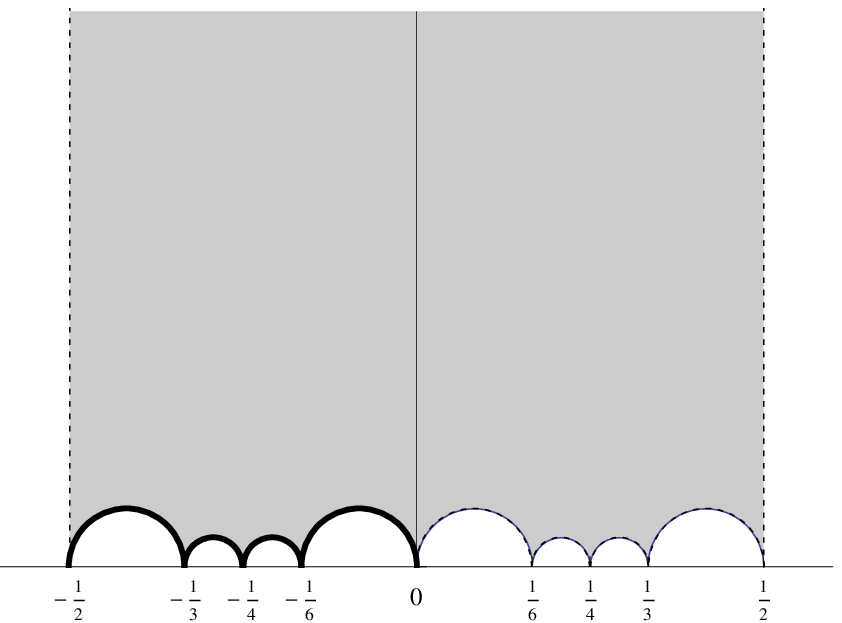}}
\quad \quad
{{$E_{k, 12}^0$}\includegraphics[width=1.5in]{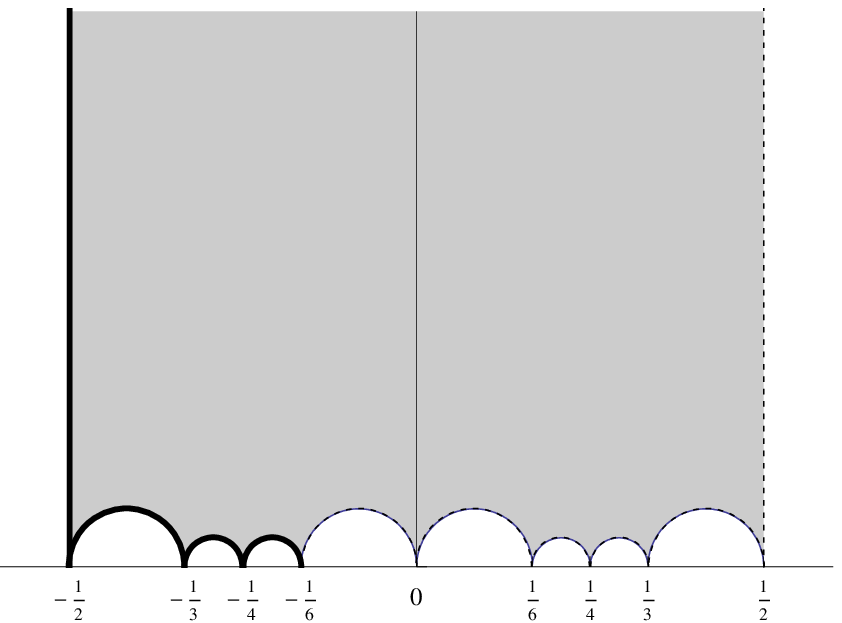}}\\

{{$E_{k, 12}^{-1/3}$}\includegraphics[width=1.5in]{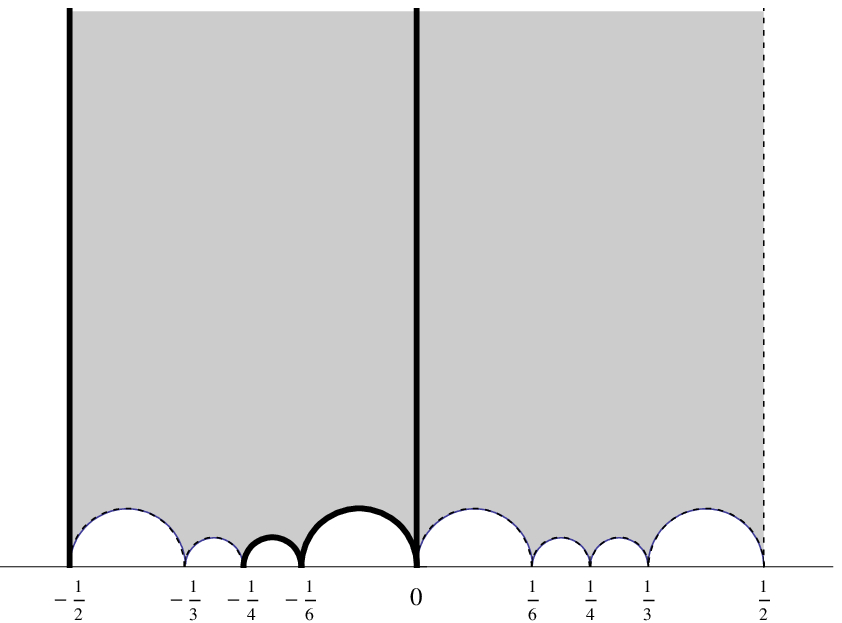}}
\quad \quad
{{$E_{k, 12}^{-1/4}$}\includegraphics[width=1.5in]{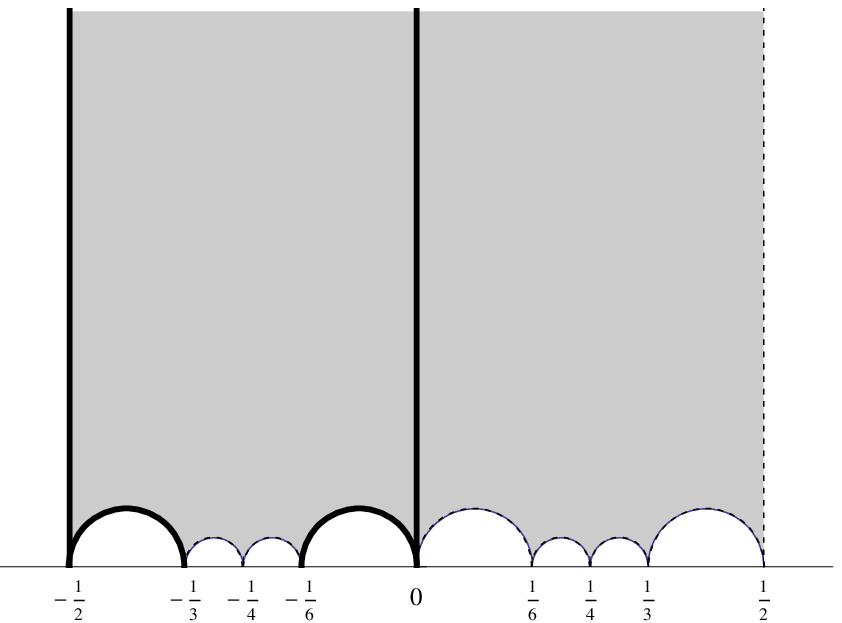}}\\

{{$E_{k, 12}^{-1/2}$}\includegraphics[width=1.5in]{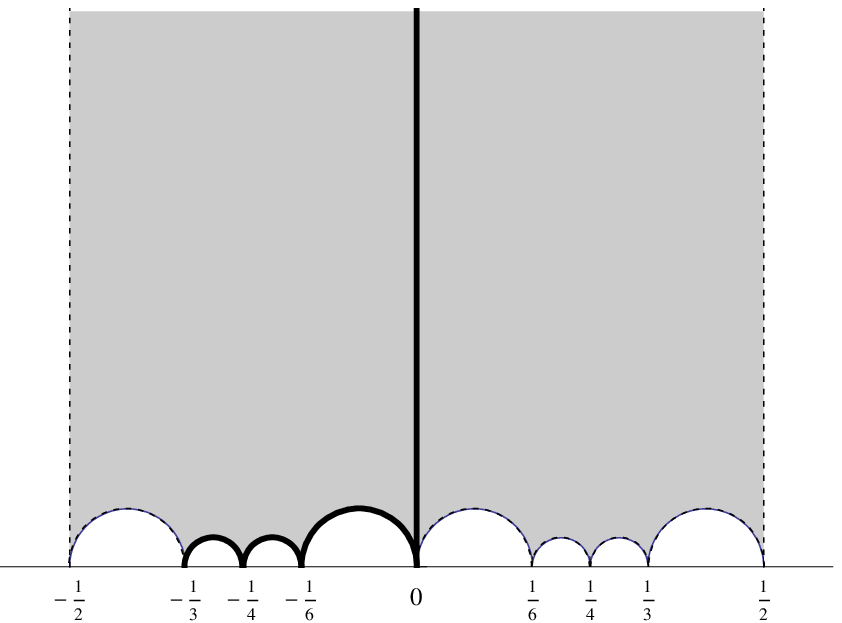}}
\quad \quad
{{$E_{k, 12}^{-1/6}$}\includegraphics[width=1.5in]{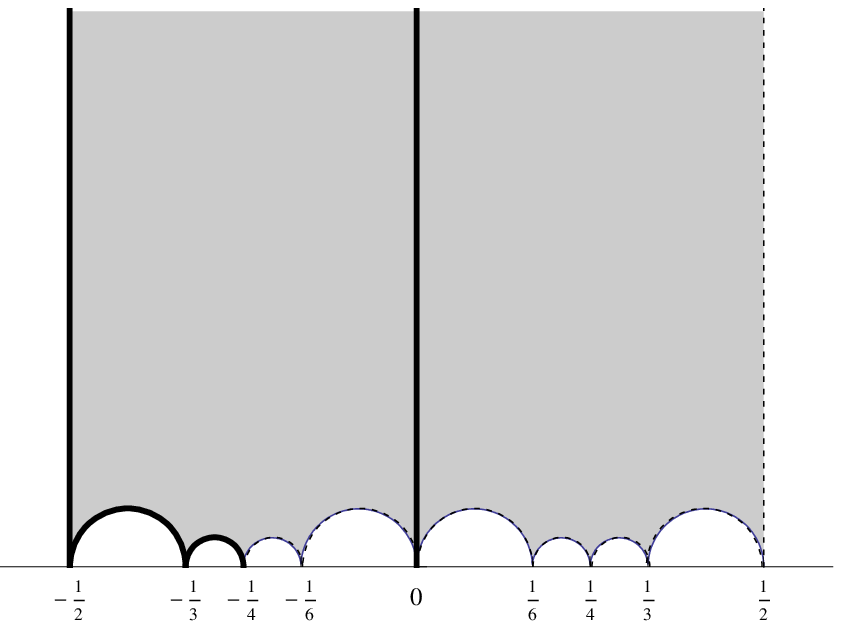}}
\end{center}
\caption{Location of the zeros of the Eisenstein series}
\end{figure}

\paragraph{\bf Location of the zeros of Hecke type Faber Polynomial}
For $m \leqslant 200$, we can prove that all of the zeros of $F_{m, 12}$ lie on the lower arcs of $\partial \mathbb{F}_{12}$ by numerical calculation. 

\markright{\sc On the zeros of Eisenstein Series for the normalizers of congruence subgroups}


\newpage

\begin{thebibliography}{MNS}

\bibitem[ACMS]{ACMS}
D. Alexander, C. Cummins J. McKay, and C. Simons, {\it Completely replicable functions}. In: {\it Groups, combinatorics $\&$ geometry} (M. Liebeck and J. Saxl eds.), 87--98, London Math. Soc. Lecture Note Ser., No. 165, Cambridge Univ. Press, Cambridge, 1992. (Proceedings of the L.M.S. Symposium on Groups and Combinatorics, Durham, 1990.)

\bibitem[AKN]{AKN}
T. Asai, M. Kaneko, and H. Ninomiya, {\it Zeros of certain modular functions and an application}, Comment.
Math. Univ. St. Paul. {\bfseries 46} (1997), 93--101.

\bibitem[BKM]{BKM}
E. Bannai, K. Kojima, and T. Miezaki, {\it On the zeros of Hecke type Faber polynomial}, to appear in Kyushu J. Math.

\bibitem[CN]{CN}
J. H. Conway, S. P. Norton, {\it Monstrous moonshine}, Bull. London Math. Soc., {\bfseries 11}(1979), 308--339.

\bibitem[G]{G}
J. Getz, {\it A generalization of a theorem of Rankin and Swinnerton-Dyer on zeros of modular forms}, Proc. Amer. Math. Soc., {\bfseries 132}(2004), No. 8, 2221--2231.

\bibitem[H]{H}
H. Hahn, {\it On zeros of Eisenstein series for genus zero Fuchsian groups}, Proc. Amer. Math. Soc. {\bfseries 135}(2007), No. 8, 2391--2401.

\bibitem[Ko]{Ko}
N. Koblitz, {\it Introduction to Elliptic Curves and Modular Forms}, Graduate Texts in Mathematics, No. 97, Springer-Verlag, New York, 1984.

\bibitem[Kr]{Kr}
A. Krieg, {\it Modular Forms on the Fricke Group.}, Abh. Math. Sem. Univ. Hamburg, {\bfseries 65}(1995), 293--299.

\bibitem[MNS]{MNS}
T. Miezaki, H. Nozaki, and J. Shigezumi, {\it On the zeros of Eisenstein series for $\Gamma_0^* (2)$ and $\Gamma_0^* (3)$}, J. Math. Soc. Japan, {\bfseries 59}(2007), 693--706.

\bibitem[Q]{Q}
H. -G. Quebbemann, {\it Atkin-Lehner eigenforms and strongly modular lattices}, Enseign. Math. (2), {\bfseries 43}(1997), No. 1-2, 55--65.

\bibitem[RSD]{RSD}
F. K. C. Rankin, H. P. F. Swinnerton-Dyer, {\it On the zeros of Eisenstein Series}, Bull. London Math. Soc., {\bfseries 2}(1970), 169--170.

\bibitem[Se]{Se}
J. -P. Serre, {\it A Course in Arithmetic}, Graduate Texts in Mathematics, No. 7, Springer-Verlag, New York-Heidelberg, 1973. (Translation of {\it Cours d'arithm\'etique $($French$)$}, Presses Univ. France, Paris, 1970.)

\bibitem[SG]{SG}
G. Shimura, {\it On Eisenstein Series}, Duke Math. J., {\bfseries 50}(1983), No. 2, 417--476.

\bibitem[SH]{SH}
H. Shimizu, {\it Hokei kansu. I-III. $($Japanese$)$ $[$Automorphic functions. I-III$]$}, Iwanami Shoten Kiso Sugaku [Iwanami Lectures on Fundamental Mathematics] 8, Iwanami Shoten Publishers, Tokyo, 1977--1978.

\bibitem[SJ1]{SJ1}
J. Shigezumi, {\it On the zeros of Eisenstein series for $\Gamma_0^{*}(p)$ and $\Gamma_0(p)$ of low levels}, M.S. thesis, Kyushu University, 2006.

\bibitem[SJ2]{SJ2}
J. Shigezumi, {\it On the zeros of the Eisenstein series for $\Gamma_0^{*}(5)$ and $\Gamma_0^{*}(7)$}, Kyushu. J. Math. {\bfseries 61}(2007), 527--549.

\end{thebibliography}
\end{document}